\documentclass[11pt,twoside]{article}

\usepackage{amsbsy,amsfonts,amsmath,amssymb,eucal,mathrsfs,graphicx}
\usepackage{mathabx}
\usepackage[all]{xy}
\usepackage{hyperref}

\renewcommand{\geq}{\geqslant}
\renewcommand{\leq}{\leqslant}
\renewcommand{\ge}{\geqslant}
\renewcommand{\le}{\leqslant}

\def\itemn#1{\item[\hspace{0.6mm} {\rm (#1)}]}
\def\itemm#1{\item[\hspace{5mm} {\rm (#1)}]}
\def\on{\stackrel}
\def\too{\longrightarrow}
\def\mto{\longmapsto}
\def\into{\hookrightarrow}
\def\onto{\twoheadrightarrow}
\def\endto{\righttoleftarrow}

\renewcommand{\tilde}{\widetilde}
\renewcommand{\hat}{\widehat}
\renewcommand{\bar}{\overline}
\renewcommand{\square}{\Box}
\def\lb{\lambda}
\newcommand \bbf[1]{\textbf{\textit{#1}}{}}
\newcounter{cpt}

\def\egaldef{\stackrel{{\rm df}}{=}}

\newtheorem{counter}[subsubsection]{$\!\!$}
\newtheorem{subcounter}[subsection]{$\!\!$}
\newenvironment{defi}{\begin{counter} \rm {\bf Definition.}}{\end{counter}}
\newenvironment{defis}{\begin{counter} \rm {\bf Definitions.}}{\end{counter}}
\newenvironment{nota}{\begin{counter} \rm {\bf Notation.}}{\end{counter}}

\newenvironment{prop}{\begin{counter} {\bf Proposition.}}{\end{counter}}
\newenvironment{lemm}{\begin{counter} {\bf Lemma.}}{\end{counter}}

\newenvironment{coro}{\begin{counter} {\bf Corollary.}}{\end{counter}}

\newenvironment{theo}{\begin{counter} {\bf Theorem.}}{\end{counter}}

\newenvironment{rema}{\begin{counter} \rm {\bf Remark.}}{\end{counter}}
\newenvironment{remas}{\begin{counter} \rm {\bf Remarks.}}{\end{counter}}

\newenvironment{noth}{\begin{counter} \rm}{\end{counter}}
\newenvironment{proo}{{\flushleft \bf Proof:}}{\hfill $\square$ \vspace{5mm}}

\DeclareMathOperator{\Hom}{Hom}
\DeclareMathOperator{\Ext}{Ext}

\DeclareMathOperator{\Gr}{Gr}

\DeclareMathOperator{\Id}{Id}
\DeclareMathOperator{\coker}{coker}
\DeclareMathOperator{\Spec}{Spec}
\DeclareMathOperator{\Spf}{Spf}

\DeclareMathOperator{\Gal}{Gal}
\DeclareMathOperator{\Mod}{Mod}
\DeclareMathOperator{\im}{im}
\DeclareMathOperator{\val}{val}

\DeclareMathOperator{\vol}{vol}

\DeclareMathOperator{\cris}{cris}

 \def\cE{{\mathcal E}} \def\cF{{\mathcal F}}
\def\cG{{\mathcal G}}  
 \def\cK{{\mathcal K}} \def\cL{{\mathcal L}}
\def\cM{{\mathcal M}}  \def\cO{{\mathcal O}}
\def\cP{{\mathcal P}} \def\cQ{{\mathcal Q}} 
\def\cS{{\mathcal S}}  \def\cU{{\mathcal U}}

\renewcommand\AA{\mathbb{A}}

\newcommand\GG{\mathbb{G}}

 \newcommand\NN{\mathbb{N}}
 \newcommand\PP{\mathbb{P}}
\newcommand\QQ{\mathbb{Q}} 
\renewcommand\SS{\mathbb{S}} \newcommand\TT{\mathbb{T}}
\newcommand\UU{\mathbb{U}} 
\newcommand\WW{\mathbb{W}} \newcommand\XX{\mathbb{X}}
\newcommand\YY{\mathbb{Y}} \newcommand\ZZ{\mathbb{Z}}

  \def\fI{{\mathfrak I}}
 \def\fK{{\mathfrak K}} 
\def\fM{{\mathfrak M}} \def\fN{{\mathfrak N}} 
  
\def\fS{{\mathfrak S}}  
  \def\fX{{\mathfrak X}}

  \def\sC{\mathscr{C}}
 \def\sE{\mathscr{E}} 
\def\sG{\mathscr{G}} \def\sH{\mathscr{H}} 
  \def\sL{\mathscr{L}}
\def\sM{\mathscr{M}}

\begin{document}

\begin{center}
{\Large \bf Models of group schemes of roots of unity}
\end{center}

\begin{center}
Ariane M\'ezard, Matthieu Romagny, Dajano Tossici
\end{center}

\begin{center}
\today
\end{center}

\bigskip
\bigskip

\begin{center}
{\bf Abstract}

\bigskip

\begin{minipage}{16cm}
{\small Let $\cO_K$ be a discrete valuation ring of mixed
characteristics $(0,p)$, with residue field $k$. Using work of
Sekiguchi and Suwa, we construct some finite flat $\cO_K$-models of
the group scheme $\mu_{p^n,K}$ of $p^n$-th roots of unity, which we call
Kummer group schemes. We carefully set out the general
framework and algebraic properties of this construction.
When $k$ is perfect and $\cO_K$ is a complete totally ramified
extension of the ring of Witt vectors $W(k)$, we provide a parallel
study of the Breuil-Kisin modules of finite flat models of $\mu_{p^n,K}$,
in such a way that the construction of Kummer groups and Breuil-Kisin
modules can be compared. We compute these objects for $n\le 3$.
This leads us to conjecture that all finite flat models of $\mu_{p^n,K}$
are Kummer group schemes.}
\end{minipage}
\end{center}

\bigskip

\tableofcontents


\section{Introduction}

\noindent {\bf 1.1. Context.}
Let $k$ be a perfect field of characteristic $p$, $W(k)$ its ring of
Witt vectors, $K_0$ the fraction field of $W(k)$, $K/K_0$ a finite
totally ramified field extension, and $\cO_K$ its ring of integers.
The aim of the
present paper is the determination of the models over $\cO_K$ of the group
scheme $\mu_{p^n,K}$ of roots of unity, or what is the same by Cartier
duality, of the cyclic group scheme $(\ZZ /p^n\ZZ)_K$. Apart from the intrinsic
interest of the problem, a first motivation
for doing this lies in the study of the representations of the absolute
Galois group of $K$. Indeed, finite flat group schemes and $p$-divisible groups
are extremely important examples of crystalline representations.
Work of Fontaine,
Breuil and Kisin has culminated into a fairly nice description of these
groups using modules with semilinear Frobenius. This description
remains however very abstract and many arithmetico-geometric properties of
the group schemes do not have an easy translation in terms of modules. Thus, one
is in search of concrete examples witnessing the constructions and conjectures
of the general theory, like the filtrations defined by Abbes-Saito~\cite{AS}
and Fargues~\cite{Fa}. We wish to provide such explicit examples and test
these general constructions.

Another important motivation is to
understand the reduction of Galois covers of $K$-varieties. In the case of
covers of curves, it is visible already for isogenies of elliptic curves
(Katz-Mazur~\cite{KM}) but also in higher genus (Abramovich-Romagny \cite{AR})
that it is necessary to let degenerate, along with the varieties, also the
Galois group of the covers. The existence of such group degenerations is
studied more precisely in \cite{Ro} and \cite{To1}.
In the particular case of cyclic covers,
this leads to the question of understanding the models of $\ZZ /p^n\ZZ $.
Here it is worth emphasizing that whereas in the context of Galois
representations one is by choice sticking to the original field
$K$, in the context of reduction of covers it is natural to allow
finite extensions $K'/K$. This enhances the importance of cyclic
$K$-group schemes, since any finite flat commutative group scheme
becomes isomorphic to a product of such after a finite field
extension.

A third motivation comes from the problem of finding an explicit
description of the Hopf algebras of group schemes over a discrete
valuation ring with prescribed generic fiber $G$, in other words the
Hopf orders of the algebra $K[G]$. The most studied
and best-known case is that of $G=\mu_{p^n,K}$. Going beyond the
Tate-Oort classification~\cite{TO}, work in this trend is due mainly to
Larson~\cite{Lar}, Greither~\cite{Gr}, Byott~\cite{By},
Underwood~\cite{Un} and Childs~\cite{Ch}. As a result, one has a
complete classification of Hopf orders for $n=1,2$ and a wealth of
examples for $n=3$. One difference between our approach and some of
these constructions is that we shall find descriptions which offer
information about the cohomology of the associated group schemes.
Another important feature of our constructions is that they require
no assumption on the discretely valued field $K$, whereas the results
obtained by the above authors are valid for $K$ complete, with perfect
residue field, containing a primitive $p^n$-th root of unity.

\bigskip

\noindent {\bf 1.2. Our approach.}
In this text, building on work of Sekiguchi and Suwa, we present a family of
finite flat models of $\mu_{p^n,K}$ which we call {\em Kummer group schemes}.
For this, we consider models of $(\GG_{m,K})^n$ constructed by successive
extensions of affine, smooth, one-dimensional models of $\GG_{m,K}$ with
connected fibres, called {\em filtered group schemes}.
Kummer group schemes are defined as the kernels~$G$ of some well-chosen
isogenies $\cE\to\cF$ between filtered group schemes, and their name
comes from the fact that the exact sequence $0\to G\to\cE\to\cF\to 0$
is an integral model of the usual Kummer isogeny. This sequence is
especially well-suited for the description of torsors under the
group schemes at hand, which as we said before is one of our motivations.
We also point out that the isogenies are given by
explicit equations, and hence so are the kernels. We formulate the
following conjecture:

\bigskip

\noindent {\bf Conjecture.}
{\em Any model of $\mu_{p^n,K}$ over $\cO_K$ is a Kummer group scheme.}

\bigskip

Our aim is to give strong evidence for this statement. We remark
that  this conjecture is true, without  assuming the discrete
valuation ring complete with perfect residue field, in the
case $n\le 2$ (for $n=1$ see e.g.~\cite{WW}, discussion after Theorem 2.5,
and for $n=2$ see~\cite{To2}). In order to explain why we think it
is true in general and what we actually do, let us first consider
the category of finite flat models of $\mu_{p^n,K}$. Using
scheme-theoretic closures, it is not hard to see that any morphism
$G\to G'$ between finite flat $\cO_K$-group schemes factors as the
composition $G\onto G/N\to H\into G'$
 of the quotient by a finite flat subgroup scheme, a morphism which is
an isomorphism on the generic fibre (a so-called {\em model map}) and the closed
immersion of a finite flat subgroup scheme. Models of $\mu_{p^n,K}$ are
special because they have a unique (finite, flat) subgroup and quotient
of a given order. Thus the category of models of $\mu_{p^n,K}$ may be completely
described by the subcategory of groups with model maps as morphisms, which
is just a partially ordered set $(\sC_n,\ge)$, and the two families
of functors $Q_i,S_i:\sC_n\to\sC_i$
given by the finite flat quotient of degree $p^{n+1-i}$, and
the finite flat subgroup of degree~$p^{n+1-i}$, for $i\in\{1,\dots,n+1\}$.

Now let us describe what we do. As we said, we take up and extend a
construction of Sekiguchi and Suwa, and use it to produce models of
$\mu_{p^n}$, the Kummer group schemes.
These group schemes are parametrized by matrices with coefficients
in the ring of Witt vectors $W(\cO_K)$. The choice of a
uniformizer $\pi$ for $\cO_K$ allows to single out a certain set
$\sM_n$ of upper triangular matrices with an interesting structure:
it is embedded in a bigger set of matrices endowed with a
non-associative product, giving rise to a natural order $\succ$.
This set has also operators $\cU^i$ and $\cL^i$ that take a matrix
to its "upper left" and "lower right" square submatrices.

Then, we study Breuil-Kisin modules of models of $\mu_{p^n,K}$. They
can be identified with $u$-integral lattices of the ring of
Laurent series $W_n(k)((u))$, where $k$ is the residue field of
$K$. The set $\sL_n$ of lattices is ordered by inclusion and is
endowed with functors $K_i,I_i:\sL_n\to\sL_i$ given by the kernel
and image of the endomorphisms $p^{n+1-i}$ and $p^{i-1}$ of a
given lattice. The lattices have unique distinguished systems of
generators whose $p$-adic coefficients can be put together into an
upper triangular matrix. In this way, we obtain a set $\sG_n$ of
matrices with coefficients in $k((u))$, with a non-associative
product very similar to that of $\sM_n$ and giving rise to a
natural order $\succ$. This set also has functors $\cU^i$
and~$\cL^i$.

Although not quite "isomorphic", the partially ordered
sets $(\sC_n,\ge)$, $(\sM_n,\succ)$
and $(\sG_n,\succ)$ with their pairs of functors have strong analogies.
There are natural functors $\sM_n\to\sC_n\overset{\sim}{\to}\sG_n$ given
by mapping a matrix to the Kummer group scheme it defines, and then to the
Breuil-Kisin module of that group. The second functor is an equivalence,
constructed in \cite{Ki1}. The basic idea to prove the conjecture above is
to compute the Breuil-Kisin modules of Kummer groups and check that all
modules can be obtained in this way. Unfortunately, there is no direct
way to compute Breuil-Kisin modules. However, computations for $n=2$
(done by Caruso in \cite{To2} Appendix A) and $n=3$ (done in the present
article) show a surprising
phenomenon: it seems that if we replace $\pi$ by $u$ in the matrix of a
Kummer group, we obtain the matrix of its Breuil-Kisin lattice.
In fact, we set up a precise, nontrivial dictionary that indicates how to
translate the congruences in a discrete valuation ring of characteristic~$0$
on one side, into congruences in a discrete valuation ring of
characteristic~$p$ on the other side. The reader may be inspired by a
look at the tables in \ref{Comp_SST_BKT_n=2} (comparison for $n=2$) and
\ref{Comp_SST_BKT_n=3} (comparison for $n=3$ under a simplifying assumption).
She/he will see for her/himself how striking the correspondence is.
However, we wish to say that writing the dictionary already for $n=3$ in the
general case seems challenging.

Finally we observe that in particular we prove that Breuil-Kisin
modules of models of $\mu_{p^n}$ are classified by
$\frac{n(n+1)}{2}$ parameters, as conjectured, in \cite{GrCh} and
\cite{ChUn}, more generally for all models of a fixed group scheme
of order $p^n$. Moreover the Kummer group schemes we constructed
form  a family with exactly this number of parameters.


\bigskip

\noindent {\bf 1.3. Summary of contents.} Here is a short overview
of the contents of the article; each section starts with a more
detailed introduction. The article is divided in two parts
written to be readable independently (to a reasonable extent).
The first part (\S\S \ref{secBK}-\ref{BKn3}) is devoted to
Breuil-Kisin Theory over a complete discrete valuation ring with
perfect residue field. We apply that theory in order
to parametrize the models of $\mu_{p^n}$ in terms of Breuil-Kisin
modules (\S \ref{secBK}). Then we explain the algebraic structure
(called a {\em loop}) of a certain set of matrices (\S
\ref{secloop}) allowing us to rewrite Breuil-Kisin modules in
matricial terms~(\S \ref{seclat}). The main result of this first
part is Theorem \ref{classification_BK} which is a computable
interpretation of Breuil-Kisin Theory (\S \ref{MatrixBK}).
The second part (\S\S \ref{secSS}-\ref{SSn3}) is devoted to
Sekiguchi-Suwa Theory over a general discrete valuation ring of
unequal characteristics. We recall the construction of {\em
filtered group schemes} and formalize it in matricial terms (\S
\ref{secSS}). Then we describe the conditions for certain model
maps of filtered group schemes to be isogenies, whose kernels are
by definition the Kummer group schemes (\S \ref{MatrixSS}).
Finally we proceed with the explicit computation of models of
$\mu_{p^3}$ (\S \ref{SSn3}) with a comparison of the congruences
coming from Breuil-Kisin Theory and Sekiguchi-Suwa Theory
(\ref{Comp_SST_BKT_n=2} an \ref{Comp_SST_BKT_n=3}).

\bigskip

\noindent {\bf 1.4. Acknowledgements.} We thank Xavier Caruso,
Marco Garuti, Noriyuki Suwa and Angelo Vistoli for interesting
conversations related to this article. We are also grateful to
Christophe Breuil for valuable comments on the genesis of the classification
of finite flat group schemes and to Lindsay Childs who kindly sent us
a version of the paper~\cite{GrCh}. We thank the referee
for his careful reading which allowed to correct several
inaccuracies. The first and second
authors especially enjoyed a stay in the Scuola Normale Superiore
di Pisa where part of this work was done. The third author had
fruitful stays at the MPIM in Bonn, at the IHES in
Bures-sur-Yvette, and spent some time in Paris to work on this
project invited by the University Paris 6, the University of
Versailles Saint-Quentin and the IHP, during the Galois Trimester.
The three authors also spent a very nice week in the CIRM in
Luminy. We thank all these institutions for their support and
hospitality.


\section{Breuil-Kisin modules and $\mu$-lattices} \label{secBK}

In this section, we recall the description due to Breuil and Kisin
of the category of finite flat group schemes (understood commutative,
of $p$-power order) in terms of modules with Frobenius. Then, we specialize
to the subcategory of models of the group scheme $\mu_{p^n}$ of roots of unity.

We fix the following notations. Let $k$ be a perfect field of characteristic $p$,
$W=W(k)$ the ring of Witt vectors with coefficients in $k$, and $\fS=W[[u]]$.
We write $W_n=W_n(k)$ the ring of Witt vectors of length $n$ and $\fS_n=W_n[[u]]$.
The rings $\fS$ and $\fS_n$ are endowed with a ring endomorphism~$\phi$
which is continuous for the $u$-adic topology, defined as the
usual Frobenius on $W_n(k)$ and by $\phi(u)=u^p$. Let $K_0$ be the fraction
field of $W(k)$, let $K/K_0$ be a totally ramified extension of degree
$e$ and $\cO_K$ its ring of integers. We fix a uniformizer $\pi$ of~$K$ and
denote by $E(u)$ its minimal polynomial over $K_0$ and $v$ the $p$-adic
valuation with $v(\pi)=1$.
We always use the phrase {\em finite flat group scheme} as a shortcut for
{\em commutative finite flat group scheme of $p$-power order}. We
denote by $(\Gr\!/\cO_K)$ the corresponding category.

\subsection{Breuil-Kisin modules of finite flat group schemes}
\label{sec:Kisin equivalence}

\begin{noth} {\bf The Breuil-Kisin Theorem.}
In recent papers, Breuil and Kisin have proven a classification theorem
for finite flat $\cO_K$-group schemes, in terms of
the category $(\Mod/\fS)$ described as follows:
\begin{itemize}
\item the objects of $(\Mod/\fS)$ are the finitely generated
$\fS$-modules $\fM$ of
projective dimension 1, killed by some power of $p$, and endowed with a
$\phi$-semilinear map $\phi_\fM:\fM\to\fM$ such that $E(u)\fM$ is contained
in the $\fS$-module generated by $\phi_\fM(\fM)$.
\item the morphisms in $(\Mod/\fS)$ are the $\fS$-linear maps compatible
with $\phi$.
\end{itemize}
For any $\fM\in (\Mod/\fS)$, the map $\phi_\fM$ is called the {\em Frobenius}
and most often written simply $\phi$. Note that to $\phi$ is
associated a {\em linear} map $\phi^*\fM\to\fM$, where
$\phi^*\fM:=\fM\otimes_{\fS,\phi}\fS$.
The classification of Breuil and Kisin is the following:

\begin{theo} \label{classif_BK}
There is a contravariant exact equivalence of categories
$(\Gr\!/\cO_K)\to (\Mod/\fS)$.
\end{theo}

One may compose with Cartier duality to get a covariant equivalence, and in
this paper this is what we will do.

The category $(\Mod/\fS)$ was introduced in~\cite{Br1}. To be more
precise, Breuil required moreover that the underlying $\fS$-module
of an object $(\Mod/\fS)$ should be a finite direct sum of modules
$\fS/p^{n_i}\fS$. He conjectured the existence of an equivalence between
$(\Mod/\fS)$ and the category of finite flat group schemes whose
$p^m$-kernels are finite flat for all $m$, and he proved the conjecture
for group schemes killed by $p$, when $p>2$. After that, Kisin realized
that arbitrary finite flat group schemes could be taken into account
by requiring the underlying modules merely to have projective dimension~$1$,
and he proceeded to prove the conjecture in general
(see~\cite{Ki1}, Thm~0.5) for $p>2$. Later Lau~\cite{Lau} and
Liu~\cite{Li} independently proved that the statement also holds for
$p=2$, which in fact was the original motivation of the note of
Breuil~\cite{Br1}.

For the convenience of the reader,
here is a very rough sketch of how the equivalence of the theorem works.
Let $S$ denote the $p$-adic completion of the divided power envelope of
$W[u]$ with respect to the ideal generated by $E(u)$. There is a natural
inclusion $\fS\to S$
but one has to notice that the ring $S$ is
much more complicated than $\fS$. Breuil introduces a category $(\Mod/S)$
whose objects are $S$-modules with a $1$-step filtration and a semi-linear
Frobenius. On the syntomic site of the formal scheme $\Spf(\cO_K)$, all
finite flat group schemes define abelian sheaves. Breuil constructs another
abelian sheaf $\cO^{\cris}_{\infty,\pi}$. This sheaf plays
the role of a sort of dualizing object: Breuil shows that there is a
contravariant equivalence $(\Gr\!/\cO_K)\to (\Mod/S)$ that takes a group
scheme $G$
to the module $\Hom(G,\cO^{\cris}_{\infty,\pi})$, with a quasi-inverse
that takes a module $\cM$ to the group scheme that represents the
syntomic sheaf $\fX\mapsto\Hom(\cM,\cO^{\cris}_{\infty,\pi}(\fX))$. Now there
is a covariant functor $(\Mod/\fS)\to (\Mod/S)$ given by tensoring with the map
$\phi:\fS\to S$. Kisin proves that for any $\cM\in (\Mod/S)$ there
is a unique sub-$\fS$-module $\fM\subset\cM$ such that
$E(u)\fM\subset \langle \phi(\fM)\rangle\subset \fM$. Moreover we can recover
$\cM$ from this submodule in the sense that $\cM\simeq \fM\otimes_{\fS,\phi} S$
so that $(\Mod/\fS)\to (\Mod/S)$ is an equivalence.
\end{noth}

\begin{noth} {\bf Group schemes killed by $p^n$.}
The modules killed by $p^n$ correspond to the group schemes killed
by $p^n$. We will use a somewhat different description of the full
subcategory of $(\Mod/\fS)$ of modules killed by $p^n$, based on
the following lemma.

\begin{lemm}
Let $\fM$ be an $\fS$-module endowed with a $\phi$-semilinear map
$\phi:\fM\to\fM$ such that $\coker(\phi^*\fM\to\fM)$ is killed by $E(u)$.
Assume that $\fM$ is  killed by $p^n$. Then $\fM$ is an $\fS$-module
of projective dimension 1 if and only if $\fM$ is a finite
$\fS_n$-module without $u$-torsion.
\end{lemm}

\begin{proo}
It follows from \cite{Ki1}, Lemma~2.3.2 that $\fM$ has projective
dimension $1$ if and only if it is an iterated extension of finite
free $\fS/p\fS$-modules. By induction, it is immediate that this
is equivalent to the fact that $\fM$
is a finite $\fS_n$-module without $u$-torsion.
\end{proo}

Therefore, the full subcategory of $(\Mod/\fS)$ of modules
killed by $p^n$ is the category $(\Mod/\fS)_n$ defined as follows:
\begin{itemize}
\item the objects of $(\Mod/\fS)_n$ are the finite $\fS_n$-modules $\fM$
with no $u$-torsion endowed with a $\phi$-semilinear map $\phi:\fM\to\fM$ such
that $\coker(\phi^*\fM\to\fM)$ is killed by $E(u)$.
\item the morphisms in $(\Mod/\fS)_n$ are the $\fS_n$-linear maps
compatible with $\phi$.
\end{itemize}
We now record some basic facts concerning $(\Mod/\fS)_n$.

\begin{lemm} \label{phi'_injective}
For any object $\fM$ of $(\Mod/\fS)_n$ the map $\phi^*\fM\to\fM$
is injective.
\end{lemm}

\begin{proo}
This is \cite{Ki2}, Lemma~1.1.9.
\end{proo}

\begin{lemm} \label{lm:kernels_images}
The category $(\Mod/\fS)_n$ has kernels, cokernels, images
and coimages. Kernels and images are given as the kernels and
images in the category of $\fS$-modules.
\end{lemm}

\begin{proo}
Let us prove first that $(\Mod/\fS)_n$ has kernels and images.
For a morphism $f:\fM\to\fN$, let $\fK$ and $\fI$ be the kernel
and the image in the category of $\fS_n$-modules. It is easy to see
that $\fK$ and $\fI$ are finite $\fS_n$-modules, stable under $\phi$,
with no $u$-torsion. Also note that the map $f':=\phi^*f:\phi^*\fM\to\phi^*\fN$
has kernel $\phi^*\fK$ (since $\phi$ is flat) and image $\phi^*\fI$.
The main point is to see that $E(u)$ kills the cokernels of the maps
$\phi^*\fK\to \fK$ and $\phi^*\fI\to \fI$. We start with the kernel.
For any $x\in \fK$ we have $x\in \fM$ and since the cokernel of $\phi^*\fM\to\fM$
is killed by $E(u)$ there exists $y\in\fM'$ such that $E(u)x=\phi(y)$.
Then $f'(y)$ maps to $0$ in $\fN$ and hence is $0$ in $\phi^*\fN$.
It follows that $y\in \phi^*\fK$, as desired. We come to the image.
Let $x\in \fI$ so that $x=f(y)$ for some $y\in\fN$. Then there
exists $z\in \phi^*\fN$ such that $E(u)y=\phi(z)$. Therefore
$E(u)x=\phi(f'(z))$ with $f'(z)\in \phi^*\fI$, as desired.

By Theorem~\ref{classif_BK}, there is on $(\Mod/\fS)_n$ a
contravariant exact involutive equivalence given by Cartier duality.
It follows that $(\Mod/\fS)_n$ has cokernels and coimages.
\end{proo}

\begin{rema}
In general, for a morphism $f:\fM\to\fN$ the objects
$\coker(\ker(f))$ and $\ker(\coker(f))$ are not isomorphic. In the
category $(\Mod/\fS)_n$ this is not so easy to see, because we
have not worked out the description of cokernels. Things are a
little easier in the category of finite flat group schemes. There,
the kernel of a map $u:G\to H$ is the scheme-theoretic closure of
the kernel of the generic fibre $u_K:G_K\to H_K$ inside $G$, and
the cokernel is the Cartier dual of the kernel of the dual of $u$.
For example, if $R$ contains a primitive $p$-th root of unity and
$u:(\ZZ /p\ZZ )_R\to\mu_{p,R}$ is an isomorphism on the generic
fibre, then $\ker(u)=\coker(u)=0$ even though $u$ is not an
isomorphism.
\end{rema}
\end{noth}

\subsection{Lattices of $W_n((u))$} \label{ss:lattices_of_Wn}

We shall see in~\ref{ss:mu_lattices} that the Breuil-Kisin modules of
models of $\mu_{p^n}$ can be identified with lattices in the
$W_n[[u]]$-module $W_n((u))$. For this reason, it is useful to collect
some basic facts on these lattices; knowing their generating systems will
be particularly important in Section~\ref{seclat}. Since the lattices we are
interested in are Breuil-Kisin modules, for simplicity we keep the letters
$\fM,\fN$ (etc.) to denote them.

\begin{defi} \label{df:lattices}
A {\em lattice} $\fM$ is a finitely generated sub-$W_n[[u]]$-module
of $W_n((u))$ such that $\fM[1/u]=W_n((u))$.
We denote by $\sL_n$ the partially ordered set of lattices with inclusions
between them. If a lattice $\fM$ is contained in $W_n[[u]]$, we say that
it is {\em positive} and we write $\fM\ge 0$.
\end{defi}

Note that it is simpler here not to follow british mathematical usage,
so we say {\em positive} instead of {\em non-negative}.

For any two lattices $\fM,\fN$, there exists $\alpha\in\NN$ such that
$u^\alpha\fN\subset \fM$. We define the
{\em volume (or index) of $\fM$ with respect to $\fN$} as
$$
\vol(\fM,\fN)=u^{n\alpha-\lg(\fM/u^\alpha\fN)}
$$
where $\lg$ denotes the length as a $W_n[[u]]$-module. Using the fact that
$\lg(\fN/u^\alpha\fN)=n\alpha$, one sees that the definition is indeed
independent of $\alpha$. Although our base ring is not a Dedekind ring,
this is the analogue of the symbol $\chi(\fM,\fN)$ of \cite{Fa}, d\'ef.~5
and \cite{Se}, chap.~III, no.~1. The {\em volume (or index) of $\fM$} is defined by
$\vol(\fM)=\vol(\fM,W_n[[u]])$, and we have $\vol(\fM,\fN)=\vol(\fM)/\vol(\fN)$.

\begin{noth} {\bf Kernels and images of $p$.} \label{ker_im_p}
For any lattice $\fM$ and any integer $i$ with $1\le i\le n+1$,
we define $\fM[i]=\ker(p^{n+1-i}:\fM\to\fM)$
and $\fM(i)=\im(p^{i-1}:\fM\to\fM)$. We have $\fM(i)\subset\fM[i]$
and these submodules fit into compatible decreasing filtrations:
$$
\xymatrix@R=.1pc@C=.4pc{
\fM= & \fM[1] & \supsetneq & \dots & \supsetneq & \fM[n]
& \supsetneq & \fM[n+1] & **[l] = 0 \ \ \\
& \cup & & & & \cup & & \cup \\
\fM= & \fM(1) & \supsetneq & \dots & \supsetneq & \fM(n)
& \supsetneq & \fM(n+1) & **[l] = 0 \ . \\}
$$
For any two submodules $\fN,\fN'$ of $W_n((u))$,
consider the ideal
$$
(\fN:\fN')=\{x\in W_n[[u]]\,,\,x\fN'\subset\fN\} .
$$
Let $1\le i\le j\le n+1$ be integers. One can see easily, by inverting $u$, that
$$
(\fM[j]:\fM[i])=(\fM(j):\fM(i))=p^{j-i}W_n[[u]].
$$
Besides, since $\fM$ has no $u$-torsion then $\fM[i]\cap\fM[j][1/u]=\fM[j]$
and the map
$$
\fM[i]/\fM[j] \too \fM[i][1/u]/\fM[j][1/u]
$$
is injective. Since $\fM[i][1/u]=p^{i-1}W_n((u))$ and
$p^{i-1}W_n((u))/p^{j-1}W_n((u))=W_{j-i}((u))$, this proves that $\fM[i]/\fM[j]$
is canonically a lattice of $W_{j-i}((u))$. Exactly the same arguments show that
$\fM(i)/\fM(j)$ is canonically a lattice of $W_{j-i}((u))$. In particular,
for $j=n+1$ this says that $\fM[i]$ and $\fM(i)$ are lattices of $W_{n+1-i}((u))$.
\end{noth}

\begin{noth} {\bf Generating sets.} \label{generating_sets}
For each $x\in k$, let $[x]\in W(k)$ be its {\em Teichm\"uller representative}
(see~\ref{ss:p_adic_exp} for a reminder on this notion).
The map $x\mapsto [x]$ is the unique multiplicative section of the projection
onto the residue field.
If $e_1,\dots,e_n$ is a set of generators for $\fM$, then we will call
{\em $T$-combination} a linear combination $t_1e_1+\dots+t_ne_n$ where
$t_1,\dots,t_n$ are Teichm\"uller representatives. In the following result,
and in other places of the paper, we use the same letter for the valuation
of a discrete valuation ring and for the induced function on its artinian
quotients.

\begin{lemm} \label{lemma_T_basis}
Let $\fM$ be a lattice of $W_n((u))$ and let $e_1,\dots,e_n$ be a system
of generators. Let $v_p$ denote the $p$-adic valuation on $W_n$.
Then the following conditions are equivalent:
\begin{trivlist}
\itemn{1} For $1\le i\le n$, we have $v_p(e_i)=i-1$ and
$pe_i\in\langle e_{i+1},\dots,e_n\rangle$.
\itemn{2} For $1\le i\le n$, we have $\fM[i]=\langle e_i,\dots,e_n\rangle$.
\itemn{3} For $1\le i\le n$, we have $v_p(e_i)=i-1$ and 
each element $x\in\fM$ can be written in a unique way as a $T$-combination
$x=[x_1]e_1+\dots+[x_n]e_n$ with $x_i\in k[[u]]$.
\end{trivlist}
\end{lemm}

\begin{proo}
(1) $\Rightarrow$ (2). Set $\fN_i=\langle e_i,\dots,e_n\rangle$.
It is obvious that $\fN_i\subset\fM[i]$,
so we only prove the opposite inclusion. Since $v_p(e_i)=i-1$,
we have $\fN_i[1/u]=p^{i-1}W_n((u))$. Let $x\in\fM[i]$
and write
$$x=x'_1e_1+\dots+x'_ne_n$$
for some coefficients $x'_i\in W_n[[u]]$. The fact that
$pe_i\in\fN_{i+1}$ implies that this linear combination may be
transformed into a $T$-combination $x=[x_1]e_1+\dots+[x_n]e_n$. If
$x\ne 0$ there exists $\nu$ minimal such that $x_\nu\ne 0$. Then
the assumption that $x\in\fM[i]$ gives $[x_\nu]e_\nu\in\fM[i]+\fN_{\nu+1}$.
After tensoring with $W_n((u))$ we obtain
$$
p^{\nu-1}W_n((u))\subset p^{i-1}W_n((u))+p^{\nu}W_n((u))
=p^{\min(i-1,\nu)}W_n((u))
$$
hence $\nu\ge i$, so that $x\in\fN_i$.

\smallskip

\noindent (2) $\Rightarrow$ (3). From $\fM[i][1/u]=p^{i-1}W_n((u))$
we deduce by decreasing induction on $i$ that $v_p(e_i)=i-1$. Now
fix $x\in\fM$. Since $p\fM[i]\subset\fM[i+1]$, we have
$pe_i\in\langle e_{i+1},\dots,e_n\rangle$ for all $i$.
Using this, we may as above write $x$ as a $T$-combination
$x=[x_1]e_1+\dots+[x_n]e_n$ with $x_i\in k[[u]]$. Moreover,
if $[x_1]e_1+\dots+[x_n]e_n=[x'_1]e_1+\dots+[x'_n]e_n$ are two
expressions for $x$, then $([x_1]-[x'_1])e_1\in\fM[2]$. From the
fact that $(\fM[2]:\fM[1])=pW_n[[u]]$ it follows that
$[x_1]-[x'_1]\in pW_{n+1}[[u]]$ and hence $x_1-x'_1=0$. By
induction we get similarly $x_i=x'_i$ for all~$i$.

\smallskip

\noindent (3) $\Rightarrow$ (1). Since $v_p(e_i)=i-1$, the $p$-valuation of a
nonzero element $[x_1]e_1+\dots+[x_n]e_n$ is equal to $\nu-1$ where
$\nu$ is the least integer such that $x_\nu\ne 0$. For $x=pe_i$ we
find $\nu=i+1$, so that $pe_i\in\langle e_{i+1},\dots,e_n\rangle$.
\end{proo}

\begin{defi} \label{df:T_basis}
A set of generators $e_1,\dots,e_n$ of a lattice $\fM$ satisfying
the equivalent conditions of Lemma~\ref{lemma_T_basis} is called
a {\em Teichm\"uller basis}, or a {\em $T$-basis} for short.
\end{defi}

\begin{rema} \label{rm:T_basis_and_li}
Let $e_1,\dots,e_n$ be a $T$-basis of $\fM$ and for each $i$, let $l_i$
be the $u$-adic valuation of the class of $e_i$ in $\fM[i]/\fM[i+1]$ which
is a lattice of $k((u))$. Then, we have $l_1\ge l_2\ge\dots\ge l_n$.
Indeed, by the definition of $l_i$, we have $e_i=\alpha_ip^{i-1} \mod p^i$
with $\val_u(\alpha_i)=l_i$. Therefore $pe_i=\alpha_ip^i \mod p^{i+1}$
and $e_{i+1}=\alpha_{i+1}p^i \mod p^{i+1}$. Since
$pe_i\in\langle e_{i+1},\dots,e_n\rangle$, it follows at once that
$l_i\ge l_{i+1}$.
\end{rema}

\begin{prop} \label{dist_generators}
Let $\fM$ be a lattice of $W_n((u))$. Then there exists a unique
$T$-basis $e_1,\dots,e_n$ of the form:
$$
e_i=u^{l_i}p^{i-1}+[a_{i,i+1}]\,p^i+[a_{i,i+2}]\,p^{i+1}+\dots+[a_{in}]\,p^{n-1}
$$
where $a_{ij}\in k[u,u^{-1}]$ is such that $\deg_u(a_{ij})<l_j$ for all
$i,j$. Moreover, we have $l_1\ge l_2\ge\dots\ge l_n$. Finally $\fM$ is
positive if and only if $l_n\ge 0$ and $a_{ij}\in k[u]$ for all $i,j$.
\end{prop}

\begin{proo}
Existence: we construct the $e_i$ by decreasing induction on $i$,
starting from $i=n$. The module $\fM[n]$ is isomorphic via a
canonical isomorphism to a lattice of $W_1((u))=k((u))$, hence generated
by $u^{l_n}$ for a unique $l_n\in \ZZ$.
The preimage via this isomorphism of this generator is $e_n=u^{l_n}p^{n-1}$.
For $i<n$, assume by induction that $e_{i+1},\dots,e_n$ have been
constructed. The module $\fM[i]/\fM[i+1]$ is again canonically a
lattice of $k((u))$, generated by $u^{l_i}$ for a unique $l_i\in\ZZ $.
Since $\fM[i]\subset p^iW_n((u))$, a lift in $\fM[i]$ of this generator
may be written in the form
$$
e_i=u^{l_i}p^{i-1}+[a_{i,i+1}]\,p^i+[a_{i,i+2}]\,p^{i+1}+\dots+[a_{in}]\,p^{n-1}
$$
for some Laurent series $a_{ij}\in k((u))$. Now write
$a_{i,i+1}=a'_{i,i+1}+u^{l_{i+1}}a''_{i,i+1}$ where
$a'_{i,i+1}\in k[u,u^{-1}]$ is the truncation of $a_{i,i+1}$ in degrees
$\ge l_{i+1}$. Replacing $e_i$ by $e_i-[a''_{i,i+1}]e_{i+1}$, and
rewriting the $p$-adic expansion of the tail $e_i-[a''_{i,i+1}]\,p^i$,
we can fulfill the condition $\deg_u(a_{i,i+1})<l_{i+1}$. Applying the
same process to $a_{i,i+s}$ for $s=1,\dots,n-i$ we can fulfill the
conditions $\deg_u(a_{ij})<l_j$ for all $j$. This finishes the construction
of $e_i$, and by induction, of $e_1,\dots,e_n$. The elements $e_i$ are
such that $\fM[i]=\langle e_i,\dots,e_n\rangle$ by construction.

Uniqueness: the choice of the generator of $\fM[i]/\fM[i+1]$ in the previous
induction is normalized by the fact that we are looking for generators $e_i$
with leading coefficients $u^{l_i}p^{i-1}$. The choice of the remaining
coefficients of $e_i$ is imposed by the condition on the degrees.
This proves that the system $e_1,\dots,e_n$ is unique. Finally the inequalities
between the $l_i$ are given by Remark~\ref{rm:T_basis_and_li} and the statement
about positivity is obvious.
\end{proo}

\begin{defi} \label{df:dist_basis}
The $T$-basis of Lemma~\ref{dist_generators} is called the
{\em distinguished basis} of $\fM$.
\end{defi}

\begin{rema}
Let $\fM$ be a lattice with distinguished basis $e_1,\dots,e_n$.
Then there exist series $b_{ij}\in k[[u]]$ and a set of equalities
$$
R_i : \quad pe_i=[b_{ii}]\,e_{i+1}+\dots+[b_{i,n-1}]\,e_n
$$
for $1\le i\le n$. It can be proven that in fact
$$
\langle\: e_1,\dots,e_n\:|\:R_1,\dots,R_n\:\rangle
$$
is a presentation by generators and relations of $\fM$ as a
$W_n[[u]]$-module. We will not need this.
\end{rema}
\end{noth}

\subsection{Breuil-Kisin modules of models of $\mu_{p^n,K}$}
\label{ss:mu_lattices}

We finally specialize to our main object of interest, namely, the
finite flat models of $\mu_{p^n,K}$.

\begin{noth} {\bf Models and $\mu$-lattices.}
The natural morphisms between
models are the {\em model maps}, which are by definition morphisms
of $R$-group schemes inducing an isomorphism on the generic fibre.
Let us see how the category of models of $\mu_{p^n,K}$ with model maps
can be described concretely in terms of Breuil-Kisin modules.

Let $\bar K$ be an algebraic closure of $K$. For any two finite flat
group schemes $G,G'$ with associated Breuil-Kisin modules $\fM,\fM'$, we have:
$$
G_K\simeq G'_K \iff G(\bar K)\simeq G'(\bar K) \iff \fM[1/u]\simeq
\fM'[1/u]
$$
where $G(\bar K)$ and $G'(\bar K)$ are viewed as representations
of the absolute Galois group $\Gal(\bar{K}/K)$. The first equivalence
is clear, let us explain briefly the second. If we introduce the Kummer
extension $K_\infty=\cup_{n\ge 0}\,K(\sqrt[p^n]{\pi})$, then a result
of Fontaine says that the module $\fM[1/u]$ determines the
$\Gal(\bar{K}/K_\infty)$-representation associated to $G$ (see \cite{Fo},
Remark A.3.4.1). By a result of Breuil (\cite{Br2}, Theorem~3.4.3),
this representation in turn determines the crystalline
$\Gal(\bar{K}/K)$-representation $G(\bar{K})$.

Recall that we are using the covariant equivalence
$(\Gr\!/\cO_K)\to (\Mod/\fS)$ given by~\ref{classif_BK} and Cartier
duality. Thus the $\fS$-module associated to the group scheme
$\mu_{p^n,R}$ is $\fM=\fS_n$ with its usual Frobenius. From this, we deduce that
$\fM$ is the module associated to a model of $\mu_{p^n,K}$ if and
only if $\fM[1/u]$ is isomorphic to $\fS_n[1/u]=W_n(k)((u))$ with its
Frobenius. Since $\fM$ has no $u$-torsion, we may then see it as a
submodule of $W_n(k)((u))$. As far as the morphisms are concerned, the
model maps correspond to {\em inclusions} between
submodules of $W_n(k)((u))$. We are lead to the following notions.

\begin{defis} \label{df:mulattices}
A {\em $\mu$-lattice} is a lattice $\fM\subset W_n((u))$ such that
$E(u)\fM\subset\langle \phi(\fM)\rangle\subset \fM$, where $\phi$ is
the Frobenius of $W_n((u))$.
We denote by $\sL_n^\mu$ the partially ordered set of
$\mu$-lattices with inclusions between them.
\end{defis}

The letter '$\mu$' reminds us of $\mu_{p^n}$. Note that since a $\mu$-lattice
$\fM$ is stable under Frobenius, it is positive, for otherwise
there would exist an element $x\in \fM$ with negative $u$-valuation
and then the valuation of $\phi^n(x)$ would tend to $-\infty$,
in contradiction with the finite generation of $\fM$.
What has been said before means that the Breuil-Kisin classification gives
an equivalence of categories between $\sL_n^\mu$ and the category
of models of $\mu_{p^n}$ with model maps.
\end{noth}

\begin{noth} {\bf Kernels and images of $p$.} \label{k_i_p}
Let $G$ be a model of $\mu_{p^n,K}$. For $1\le i\le n+1$, define:
\begin{itemize}
\item $G[i]$ the scheme-theoretic closure of $\ker(p^{n+1-i}:G_K\to G_K)$ in $G$,
\item $G(i)$ the scheme-theoretic closure of $\im(p^{i-1}:G_K\to G_K)$ in $G$.
\end{itemize}
These are finite flat models of $\mu_{p^{n+1-i}}$.
By definition, there are exact sequences:
$$
\xymatrix{
0 \ar[r] & G[n\!+\!2\!-\!i] \ar[r] & G \ar[r]^{} \ar[d]_>>>>>>>>{p^{i-1}}
& G(i) \ar[r] \ar@{.>}[lld] & 0 \\
0 \ar[r] & G[i] \ar[r] & G \ar[r]^>>>>{} & G(n\!+\!2\!-\!i) \ar[r] & 0 \\}
$$
On the generic fibre, the vertical map $p^{i-1}:G\to G$ vanishes
on $G[n+2-i]$ and its image is a subscheme of $G[i]$. By taking closures,
the same is true everywhere. Therefore, this map induces a morphism of
$R$-group schemes $G(i)\to G[i]$ which is a model map.

Let $\fM$ be the $\mu$-lattice associated to $G$. Starting from the exact
sequences above and using the fact that the
Breuil-Kisin equivalence is exact, we see that $\fM[i]$ is the $\mu$-lattice
of $G[i]$ and $\fM(i)$ is the $\mu$-lattice of $G(i)$. Moreover, the
inclusion $\fM(i)\subset \fM[i]$ and the model map $G(i)\to G[i]$ correspond
to each other.
\end{noth}


\section{The loop of $\mu$-matrices} \label{secloop}\label{section:loop}

In~\ref{ss:lattices_of_Wn}, we have seen that lattices have "nice"
systems of generators. The $p$-adic coefficients of such systems of
generators may be put together into "nice" matrices, called $\mu$-matrices.
We will come back to this in more detail in Section~\ref{seclat}.
In the present section, we focus on the abstract algebra of the set of
$\mu$-matrices. This set has a natural operation $(A,B)\mapsto A*B$ whose
meaning is that if $i:\fM\to\fN$ is an inclusion of lattices, if $A$ is a
matrix associated with a generating system of $\fM$ and if $B$ is a matrix
associated with the inclusion $i$, then $A*B$ is a matrix associated with a
generating system of $\fN$. The operation $*$ is unfortunately neither
associative nor commutative. Still, a good surprise is that $\mu$-matrices
all lie naturally in a set where the operation $*$ becomes
invertible on the left and on the right; this set plays the same role as the
symmetrization of a commutative monoid. The structure that we obtain,
called a {\em loop}, was considered by Manin \cite{Ma} in his study of rational
points on cubic hypersurfaces, essentially because the analogue of the addition
of elliptic curves in higher dimensions fails to be associative.

The key to everything in this section is the use of $p$-adic expansions, which
exist as soon as the coefficient ring of the Witt vectors is a perfect ring of
characteristic $p$. Thus we fix such a perfect ring throughout Section~\ref{secloop}.
For simplicity we denote it by the letter $k$, but note that it need not
be a field. As before, we set $W=W(k)$ and $W_n=W_n(k)$.

Finally we point out that the role of Witt vectors will be very different
in Sections \ref{secSS} to \ref{SSn3}, where we will consider
arbitrary $\ZZ_{(p)}$-algebras as coefficient rings. We will emphasize
this in due time.

\subsection{$p$-adic expansions} \label{ss:p_adic_exp}

\begin{noth} {\bf $p$-adic expansions of Witt vectors.}
Recall that the ring structure of $W$ is given by universal
polynomials with coefficients in $\ZZ$ in countably many variables
$X_0,X_1,X_2,\dots$ For example, there are polynomials
$S_i=S_i(X_0,\dots,X_i)$ and $P_i=P_i(X_0,\dots,X_i)$, for $i\ge 0$,
giving the addition and the multiplication of two vectors
$a=(a_0,a_1,a_2,\dots)$ and $b=(b_0,b_1,b_2,\dots)$ by the rules:
$$
\begin{array}{c}
a+b=(S_0(a,b),S_1(a,b),S_2(a,b),\dots), \medskip \\
ab=(P_0(a,b),P_1(a,b),P_2(a,b),\dots) .
\end{array}
$$
Moreover, since $k$ is perfect all elements have $p$-adic expansions:
$$
a=(a_0,a_1,a_2,\dots)=[a_0]+[a_1^{1/p}]\,p+[a_2^{1/p^2}]\,p^2+\dots
$$
where $[x]:=(x,0,0,\dots)$ is the Teichm\"uller lift of $x\in k$.
Hence the functions $\SS_i$ and $\PP_i$ defined by
$\SS_i(a,b):=S_i(a,b)^{1/p^i}$ and $\PP_i(a,b):=P_i(a,b)^{1/p^i}$
satisfy
$$
\begin{array}{c}
a+b=[\SS_0(a,b)]+[\SS_1(a,b)]\,p+[\SS_2(a,b)]\,p^2+\dots, \medskip \\
ab=[\PP_0(a,b)]+[\PP_1(a,b)]\,p+[\PP_2(a,b)]\,p^2+\dots
\end{array}
$$
In fact, we can define functions $\SS_i$ and $\PP_i$ in any number $r$
of variables by the identities
$$
\begin{array}{c}
a_1+\dots+a_r=[\SS_0(a_1,\dots,a_r)]+[\SS_1(a_1,\dots,a_r)]\,p
+[\SS_2(a_1,\dots,a_r)]\,p^2+\dots, \medskip \\
a_1\dots a_r=[\PP_0(a_1,\dots,a_r)]+[\PP_1(a_1,\dots,a_r)]\,p
+[\PP_2(a_1,\dots,a_r)]\,p^2+\dots
\end{array}
$$
\end{noth}

\begin{noth} {\bf $p$-adic expansions of series.}
We wish to extend the formalism of $p$-adic expansions to the ring
of Laurent series $W((u))$. For this, we extend the definition of
Teichm\"uller lifts to elements $x\in k((u))$ as follows:
if $x=\sum_{j\gg -\infty}\,x_ju^j$ with $x_j\in k$, we set
$$
[x]=\sum_{j\gg -\infty}\,[x_j]\,u^j \ .
$$
Then it is easy to see that for a Laurent series
$a=\sum_{j\gg -\infty}\,a_ju^j$ in $W((u))$, by writing down $p$-adic
expansions of its coefficients one obtains a $p$-adic expansion
$$
a=[a_0]+[a_1]\,p+[a_2]\,p^2+\dots
$$
Let $a=\sum_{j\gg -\infty} a_ju^j$ and $b=\sum_{j\gg -\infty} b_ju^j$
be Laurent series with coefficients in $W$. We extend the definition
of $\SS_i$ by setting
$$
\SS_i(a,b)=\sum_{j\gg -\infty} \SS_i(a_j,b_j)u^j
=\left(\sum_{j\gg -\infty} S_i(a_j,b_j)u^{jp^i}\right)^{1/p^i}
$$
and one verifies immediately that the formula
$a+b=\sum_{i\ge 0} [\SS_i(a,b)]\,p^i$ remains valid. Similarly,
one extends the definition of $\SS_i(a_1,\dots,a_r)$ for Laurent
series $a_s\in W((u))$ in an obvious way. We now come to products.
There are functions $\PP_i$ such that for any $r$ Laurent series
$a_s=\sum_{i\gg -\infty}\,a_{s,i}u^i$ with coefficients in $W$ we have
$$
a_1\dots a_r=[\PP_0(a_1,\dots,a_r)]+[\PP_1(a_1,\dots,a_r)]\,p+
[\PP_2(a_1,\dots,a_r)]\,p^2+\dots
$$
It is a simple exercise to verify that
$$
\PP_i(a_1,\dots,a_r)
=\sum_j\SS_i(\dots\,,\,a_{1,j_1}\cdots a_{r,j_r}\,,\,\dots)\,u^i
$$
where the arguments of $\SS_j$ are all the finitely many possible
products $a_{1,j_1}\cdots a_{r,j_r}$ indexed by $r$-tuples
$(j_1,\dots,j_r)$ such that $j_1+\dots+j_r=j$. For example, if $a$
and $b$ are power series (i.e. Laurent series with nonnegative $u$-valuation)
we have:
$$
\PP_i(a,b)=\sum_j\,\SS_i(a_0b_j,\dots,a_jb_0)\,u^j \ .
$$
\end{noth}

\begin{noth} {\bf A warning on the use of $\SS_i$ and $\PP_i$.}
In the sequel, we will most often use $\SS_i$ and $\PP_i$ for Teichm\"uller
elements $a_i=[x_i]$. In this case, we will usually write
$\SS_i(x_1,\dots,x_r)$ and $\PP_i(x_1,\dots,x_r)$ instead of
$\SS_i([x_1],\dots,[x_r])$ and $\PP_i([x_1],\dots,[x_r])$. This is not
dangerous, but for $x,y\in k((u))$ one must be careful to distinguish between
the sum $x+y$ in $k((u))$ and the sum $[x]+[y]$ of their Teichm\"uller
representatives in $W((u))$. For example, the associativity of the sum of
Witt vectors gives for any elements $a,b,c\in W((u))$ the formula
$\SS_1(a,b,c)=\SS_1(a+b,c)$, and here the sum $a+b$ takes place in $W((u))$.
The reader is invited to compare with formula~\ref{lm:easy_formulas}(1) below.
Among the many formulas relating the $\SS_i$ and the $\PP_i$, most of
them coming from associativity and distributivity of the sum and product
of Witt vectors, we give a few examples:

\begin{lemm} \label{lm:easy_formulas}
Let $a,b,c\in k((u))$ and let $\val$ denote the $u$-valuation. We have:
\begin{trivlist}
\itemn{1} $\SS_1(a,b,c)=\SS_1(a,b)+\SS_1(a+b,c)$.
\itemn{2} $\SS_1(a,b-a)=\SS_1(a,-b)$.
\itemn{3} $\val\left(\SS_i(a,b)\right)\ge \max(\val(a),\val(b))$
for all $i\ge 1$.
\itemn{4} $[a][b]=[ab]$ if $a$ or $b$ is a monomial.
\end{trivlist}
\end{lemm}

Note that the multiplicativity formula $[a][b]=[ab]$ for $a,b\in k$
does not hold in full generality if $a,b\in k((u))$.

\begin{proo}
(1) This comes from the associativity of the sum of Witt vectors.

\smallskip

\noindent (2) It is enough to prove that $S_1(a,b-a)=S_1(a,-b)$.
This can be proven over $\ZZ $, where it follows from
the formula $S_1(x,y)=\frac{1}{p}(x^p+y^p-(x+y)^p)$.

\smallskip

\noindent (3) This comes from the fact that if we write
$a=\sum a_ju^j$ and $b=\sum b_ju^j$, then $\SS_i(a_j,b_j)=0$
as soon as $a_j=0$ or $b_j=0$.

\smallskip

\noindent (4) This is clear.
\end{proo}
\end{noth}

\begin{noth} {\bf $p$-adic expansions of vectors and matrices.}
\label{Teichmueller_reps_for_matrices}
For the computations inside lattices, we will use the notations of
linear algebra. The vectors are all column vectors. If $A$ is a
rectangular matrix with entries $a_{ij}$ in $k((u))$ (for example
$A$ could be a column vector), we will denote by $[A]$ the matrix whose
entries are the Teichm\"uller representatives $[a_{ij}]$. Thus the
entries of $[A]$ are (possibly truncated) Witt vectors. We may as above
consider $p$-adic expansions of matrices with entries in $W((u))$,
but we will have no need for this. For us, the most important vector
will be
$$
p^\star=
\left(
\begin{array}{c}
1 \\ p \\ p^2 \\ \vdots \\
\end{array}
\right)
$$
which for convenience may denote a vector with finitely, or infinitely
many, coefficients. Thus if $x\in W_n((u))^n$ is a vector with
components $x_1,\dots,x_n$ we have:
$$
{}^txp^\star=x_1+x_2p+\dots+x_np^{n-1} \ .
$$
If the $x_i$ are Teichm\"uller representatives, then this linear
combination is called a {\em $T$-combination}. Of course, any linear
combination can be transformed into a $T$-combination:

\begin{lemm} \label{map_rho}
For any rectangular matrix $A$ with entries in $W_n((u))$ with $n$
columns, there is a unique matrix $\rho(A)$ of
the same size with entries in $k((u))$ such that
$$Ap^\star=[\rho(A)]\,p^\star \ .$$
If the entries of $A$ are power series in $u$, or Laurent polynomials,
or polynomials, then so are the entries of $\rho(A)$.
If $A$ is upper triangular (resp. with Teichm\"uller diagonal entries),
then so is $\rho(A)$.
\end{lemm}

\begin{proo}
The equality $Ap^\star=[\rho(A)]\,p^\star$ is equivalent to
finitely many equalities, one for each line of $A$. Thus it is
enough to consider the case where $A$ has
only one line $A=(a_1 \dots a_n)$. Write the
$p$-adic expansion
$$
a_1+a_2p+\dots+a_np^{n-1}=[a'_1]+[a'_2]\,p+\dots+[a'_n]\,p^{n-1} \ .
$$
Obviously the desired matrix is $\rho(A)=(a'_1 \dots a'_n)$.
The remaining assertions are clear.
\end{proo}

There is an algorithmic point of view on the computation of $\rho(A)$ that
will be useful.
In order to explain this, for a coefficient in position $(i,j)$ in an
upper triangular square matrix, let us call
the difference $j-i$ the {\em distance to the diagonal}.

\begin{lemm} \label{algorithm}
Let $\sE$ be the set of upper triangular square matrices of size $n$
with entries in $W((u))$ with Teichm\"uller diagonal entries.
Define a function $F:\sE\to \sE$ as follows. Given a matrix $A$, for $i=1$
to $n$ apply the following rule to the $i$-th line:
\begin{itemize}
\item Find the first non-Teichm\"uller coefficient $a_{i,\nu}$.
\item Write the truncated $p$-adic expansion
$a_{i,\nu}p^{\nu-1}=[a'_{i,\nu}]\,p^{\nu-1}+\dots+[a'_{i,n}]\,p^{n-1}
\; \mod p^n$.
\item Replace $a_{i\nu}$ by $[a'_{i\nu}]$ and for $j>i$ replace
$a_{ij}$ by $a_{ij}+[a'_{ij}]$.
\end{itemize}
After the step $i=n$ has been completed, call the result $F(A)$.
Then, for all $k\ge 0$ we have:
\begin{itemize}
\item the coefficients with distance to the diagonal $\le k$ of the
matrix $F^k(A)$ are Teichm\"uller, where $F^k$ is the $k$-th iterarate of $F$.
\item $F^k(A)\,p^\star=A\,p^\star$.
\end{itemize}
In particular $F^{n-1}(A)=\rho(A)$.
\end{lemm}

\begin{proo}
This is obvious.
\end{proo}
\end{noth}

\subsection{The loop of $\mu$-matrices}

\begin{noth} {\bf Quasigroups and loops.} \label{df:quasigroups}
We start with some definitions from quasigroup theory, referring to the
book of Smith \cite{Sm} for more details. A {\em magma} is a set $X$ endowed with a
binary operation $X\times X\to X$, $(x,y)\mapsto xy$ usually called
{\em multiplication}. A {\em submagma} is a subset $Y\subset X$ that is closed
under multiplication. A {\em quasigroup} is a magma where left and right
division are always possible, in the sense that left multiplications $L_x$
and right multiplications $R_y$ are bijections. Given $x,y\in X$, the unique
element $a$ such that $ax=y$ is denoted $y/x$ (read ``$y$ over $x$'')
and the unique element $b$ such that $xb=y$ is denoted
$x\backslash y$ (read ``$x$ into $y$''). A {\em loop} ({\em boucle} in French,
and... {\em loop} in Italian)
is a quasigroup
with an identity element, i.e. an element $e\in X$ such that
$ex=xe=x$ for all $x\in X$. Thus a loop is a group if and only if the
operation is associative. A {\em magma homomorphism} is a map $f:X\to X'$ such
that $f(x_1x_2)=f(x_1)f(x_2)$ for all $x_1,x_2\in X$.
{\em Quasigroup homomorphisms} and {\em loop homomorphisms} are just magma
homomorphisms.
\end{noth}

\begin{noth} {\bf The loop $\sG_n((u))$.}
In Section~\ref{seclat}, to lattices of $W_n((u))$ we will attach
matrices. The matrices coming in this way appear naturally as objects
in a certain loop which we call the {\em loop of $\mu$-matrices} and
denote by $\sG_n((u))$.
As a set, it is composed of the upper triangular matrices of the form
$$
M(\bbf{l},\bbf{a})=\left(
\begin{array}{ccccc}
u^{l_1} & a_{12} & a_{13} & \dots & a_{1n} \\
& u^{l_2} & a_{23} & & a_{2n} \\
& & \ddots & \ddots & \vdots \\
& & & u^{l_{n-1}} & a_{n-1,n} \\
0 & & & & u^{l_n} \\
\end{array}
\right)
$$
with $\bbf{l}=(l_1,\dots,l_n)\in\ZZ ^n$ and $\bbf{a}=(a_{ij})_{1\le i<j\le n}$
where $a_{ij}\in k((u))$. There is a natural subset $\sG_n[u,u^{-1}]$
composed of matrices with coefficients in $k[u,u^{-1}]$.
In order to keep the notation light, we do not specify the coefficient
ring $k$ in the symbols $\sG_n((u))$ and $\sG_n[u,u^{-1}]$.
Note also that as a general rule, we write $a_{ij}$ instead of $a_{i,j}$,
unless this can disturb comprehension, for example when we write $a_{np,n}$.

If $A,B$ are square matrices with entries in $k((u))$, we set
$A*B=\rho([A][B])$ where $\rho$ is the map from Lemma~\ref{map_rho}.
This matrix is characterized by the equality:
$$
[A][B]\,p^\star=[A*B]\,p^\star \ .
$$
By Lemma~\ref{map_rho}, if $A,B$ are in $\sG_n((u))$ resp. in
$\sG_n[u,u^{-1}]$, then $A*B$ also. It is clear that the identity
matrix is a neutral element for this multiplication. Thus the
triple $(\sG_n((u)),*,\Id)$ is a magma with identity, and
$(\sG_n[u,u^{-1}],*,\Id)$ is a submagma. At this point,
the reader may wish to have a look at the shape of the multiplication
$*$ in the examples of~\ref{examples_for_small_n} below.
\end{noth}

We will now prove that $(\sG_n((u)),*,\Id)$ is a loop.

\begin{prop} \label{sG_n_is_a_loop}
Let $A=M(\bbf{l},\bbf{a})$ and $B=M(\bbf{m},\bbf{b})$ be elements of $\sG_n((u))$.
\begin{trivlist}
\itemn{1} Any coefficient in position $(i,j)$ of $A*B$ with distance
to the diagonal $j-i\ge 1$ has the form:
$$
u^{m_j}a_{ij}+u^{l_i}b_{ij}+
\left(
\begin{array}{c}
\mbox{terms depending on coefficients $a_{i'j'}$ and $b_{i'j'}$} \\
\mbox{whose distance to the diagonal is $j'-i'<j-i$}.
\end{array}
\right).
$$
\itemn{2} The maps $L_A:B\mapsto A*B$ and $R_B:A\mapsto A*B$
are bijections.
\end{trivlist}
Thus, the triple $(\sG_n((u)),*,\Id)$ is a loop.
\end{prop}

\begin{proo}
(1) The entry of $[A][B]$ in position $(i,j)$ is
$$
u^{l_i}[b_{ij}]+\,\left(\sum_{k=i+1}^{j-1}\,[a_{ik}][b_{kj}]\right)
\, + [a_{ij}]u^{m_j} \ .
$$
The coefficients $[a_{ik}]$ and $[b_{kj}]$ in the middle sum have
distance to the diagonal strictly less than $j-i$. When applying the algorithm
of Lemma~\ref{algorithm} to compute $A*B$, at each step the entry $(i,j)$
is replaced by itself plus some terms involving coefficients $a_{st}$
and $b_{st}$ of distance to the diagonal $t-s<j-i$. This proves the claim.

\smallskip

\noindent (2)
The argument is the same for $L_A$ and $R_B$ so we do only the case of
$L_A$. Assume that $A*B=C$ with $A=M(\bbf{l},\bbf{a})$, $B=M(\bbf{m},\bbf{b})$,
$C=M(\bbf{n},\bbf{c})$. We fix $A$ and $C$ and try to solve for $B$. We
determine its entries by increasing induction on the distance
to the diagonal, called $k$. For $k=0$ it is clear that we have $m_i=n_i-l_i$.
By induction, using point (1), it follows directly that the coefficients
$b_{ij}$ of distance to the diagonal $k$ are determined by the entries of
$A$, $C$ and the coefficients $b_{i'j'}$ of lower distance to the diagonal.
\end{proo}

\begin{noth} {\bf Some subloops. The homomorphisms $\cU$ and $\cL$.} \label{U_and_L}
There are some important examples of subloops and loop homomorphisms.
Of course $\sG_n[u,u^{-1}]$ is a subloop of $\sG_n((u))$. Another
example is the subloop of matrices with diagonal entries equal to $1$.
This is in fact the kernel of the morphism of loops
$\varphi:\sG_n((u))\to\ZZ ^n$ to the additive group $\ZZ ^n$ that maps
$A$ to the tuple of its diagonal exponents.

For any square matrix $A$ of size $n$ with entries in some ring,
we denote by $\cU A$ the upper left square submatrix of size $n-1$, i.e.
the matrix obtained by deleting the last row and the last column of $A$.
Similarly we denote by $\cL A$ the lower right square submatrix of size $n-1$,
obtained by deleting the first row and the first column of $A$.

\begin{lemm} \label{homomorphisms_U_and_L}
The mappings $\cU:\sG_n((u))\to \sG_{n-1}((u))$ and $\cL:\sG_n((u))\to \sG_{n-1}((u))$
are commuting loop homomorphisms.
\end{lemm}

\begin{proo}
Let $\tau_{\cU}$ be the truncation map that takes a vector $v$ with $n$
components to the vector whose components are the first $n-1$
components of $v$. Thus $\tau_{\cU} p^\star$ is the vector analogous to
$p^\star$ in dimension one less. Then simple matrix formulas yield:
\begin{align*}
[\cU(A*B)]\,\tau_{\cU} p^\star & = (\cU[A*B])\,\tau_{\cU} p^\star=\tau_{\cU}([A*B]\,p^\star)
=\tau_{\cU}([A][B]\,p^\star) \\
& = \cU[A]\cdot\cU[B]\,\tau_{\cU} p^\star=[\cU A][\cU B]\,\tau_{\cU} p^\star
=[\cU A*\cU B]\,\tau_{\cU} p^\star \ .
\end{align*}
It follows that $\cU(A*B)=\cU A*\cU B$, that is, $\cU$ is a loop homomorphism.

Let $\tau_{\cL}$ be the truncation taking a vector $v$ with $n$ components
to the vector whose components are the last $n-1$ components of $v$.
Thus $\tau_{\cL} p^\star$ is the column vector with components
$p,p^2,\dots,p^{n-1}$. It is still true that if two square matrices
$A,B$ of size $n-1$ with coefficients in $k((u))$ satisfy
$[A]\,\tau_{\cL} p^\star=[B]\,\tau_L p^\star$ then $A=B$. Then a similar
computation as before shows that
$[\cL(A*B)]\,\tau_{\cL} p^\star=[\cL A*\cL B]\,\tau_{\cL} p^\star$, so $\cL$ is a loop
homomorphism.

Finally, the fact that $\cU$ and $\cL$ commute is clear.
\end{proo}
\end{noth}

\begin{noth} {\bf Positive matrices.} \label{def:positive elements}
We say that a matrix $A\in \sG_n((u))$ is {\em positive}, and we write
$A\ge 0$, if its entries are in $k[[u]]$. (Here, as in \ref{df:lattices},
we say {\em positive} instead of {\em non-negative} for simplicity.)
We denote by $\sG_n[[u]]$
the subset of positive elements of $\sG_n((u))$. It is a submagma,
but not a subloop. Similarly $\sG_n[u,u^{-1}]$ has a submagma
$\sG_n[u]=\sG_n[u,u^{-1}]\cap \sG_n[[u]]$.

\end{noth}

\subsection{Examples} \label{examples_for_small_n}

Here is what the operation $*$ looks like for $n=4$.
The product $P=A*B$ is given by
$$
P=\left(
\begin{array}{cccc}
u^{l_1+m_1} &
u^{l_1}b_{12}+u^{m_2}a_{12} &
p_{13} &
p_{14} \\
0 &
u^{l_2+m_2} &
u^{l_2}b_{23}+u^{m_3}a_{23} &
p_{24} \\
0&0 & u^{l_3+m_3} & u^{l_3}b_{34}+u^{m_4}a_{34} \\
0&0 &0 & u^{l_4+m_4} \\
\end{array}
\right) \ .
$$
with
\begin{align*}
p_{13} = & \ u^{l_1}b_{13}+a_{12}b_{23}+u^{m_3}a_{13}
+\SS_1(u^{l_1}b_{12},u^{m_2}a_{12}) \\
p_{24} = & \ u^{l_2}b_{24}+a_{23}b_{34}+u^{m_4}a_{24}+
\SS_1(u^{l_2}b_{23},u^{m_3}a_{23}) \\
p_{14} = &
\ u^{l_1}b_{14}+a_{12}b_{24}+a_{13}b_{34}+a_{14}u^{m_4}+
\SS_2(u^{l_1}b_{12},u^{m_2}a_{12})+ \\
& \quad + \SS_1(u^{l_1}b_{13},a_{12}b_{23},u^{m_3}a_{13},
\SS_1(u^{l_1}b_{12},u^{m_2}a_{12}))
+\PP_1(a_{12},b_{23}) \ .
\\
\end{align*}
Applying the homomorphism $\cU$ (Lemma~\ref{homomorphisms_U_and_L}),
these formulas contain also the formulas of multiplication for
$n\le 4$.

\begin{noth} {\bf Failure of associativity.}
For $n=2$, the loop $\sG_2$ is a group: in fact the multiplication $*$ is the
ordinary multiplication of matrices.
For $n\ge 3$, the multiplication $*$ is not associative. Let us check this.
We have
$$
A*B=\left(
\begin{array}{ccc}
u^{l_1+m_1} &
u^{l_1}b_{12}+u^{m_2}a_{12} &
(A*B)_{13} \\
0 &
u^{l_2+m_2} &
u^{l_2}b_{23}+u^{m_3}a_{23} \\
0&0 & u^{l_3+m_3} \\
\end{array}
\right)
$$
with
$$
(A*B)_{13} = u^{l_1}b_{13}+a_{12}b_{23}+u^{m_3}a_{13}
+\SS_1(u^{l_1}b_{12},u^{m_2}a_{12}) \ .
$$
We now examine the coefficients in position $(1,3)$:
\begin{align*}
((A*B)*C)_{13} =
& \ u^{l_1+m_1}c_{13}+(u^{l_1}b_{12}+u^{m_2}a_{12})c_{23} \\
& +u^{n_3}\big(u^{l_1}b_{13}+a_{12}b_{23}+u^{m_3}a_{13}
+\SS_1(u^{l_1}b_{12},u^{m_2}a_{12})\big) \\
& +\SS_1\big(u^{l_1+m_1}c_{12},u^{n_2}(u^{l_1}b_{12}
+u^{m_2}a_{12})\big)
\end{align*}
and
\begin{align*}
(A*(B*C))_{13} =
& u^{l_1}\big(u^{m_1}c_{13}+b_{12}c_{23}+u^{n_3}b_{13}
+\SS_1(u^{m_1}c_{12},u^{n_2}b_{12})\big) \\
& +a_{12}(u^{m_2}c_{23}+u^{n_3}b_{23})
+u^{m_3+n_3}a_{13}
+\SS_1\big(u^{l_1}(u^{m_1}c_{12}
+u^{n_2}b_{12}),u^{m_2+n_2}a_{12}\big) \ .
\end{align*}
Using the formula $\SS_1(x,y,z)=\SS_1(x,y)+\SS_1(x+y,z)$ from
Lemma~\ref{lm:easy_formulas}, we compute the difference:
\begin{align*}
& ((A*B)*C)_{13}-(A*(B*C))_{13} \\
& \qquad\qquad = \ \SS_1\big(u^{l_1+n_3}b_{12},u^{m_2+n_3}a_{12}\big)
+\SS_1\big(u^{l_1+m_1}c_{12},u^{n_2}(u^{l_1}b_{12}+u^{m_2}a_{12})\big) \\
& \qquad\qquad\qquad - \SS_1(u^{l_1+m_1}c_{12},u^{l_1+n_2}b_{12}\big)
- \SS_1\big(u^{l_1}(u^{m_1}c_{12}+u^{n_2}b_{12}),u^{m_2+n_2}a_{12}\big) \\
& \qquad\qquad = \ (u^{n_3}-u^{n_2})
\SS_1\big(u^{l_1}b_{12},u^{m_2}a_{12}\big) \ .
\end{align*}
This is not zero so $*$ is not associative.

However, we see that this is zero on the subloop
$\ker\phi:\sG_3((u))\rightarrow \ZZ^3$, which then is a group.
Let us verify that for $n\ge 4$, the multiplication $*$ is not associative
even if we restrict it to the subloop $\ker\phi:\sG_4((u))\rightarrow \ZZ^4$.
We shall
check this only for $n=4$. We make the following observation: the
multiplication of $\sG_n((u))$ differs from that of the underlying group
of matrices by terms coming from the operations of Witt vectors, i.e.
involving the sum and product functions $\SS_i$ and $\PP_j$. Since the
ordinary multiplication of matrices is associative, the terms of the
entries in $(A*B)*C$ and $A*(B*C)$ that do not involve $\SS_i$ or
$\PP_j$ are equal. Consequently when we question associativity it is
enough to look at the terms that contain $\SS_i$ or $\PP_j$.
Once this is said, let us compare the entries in position $(1,4)$ of
$(A*B)*C$ and $A*(B*C)$. Looking at the above formulas, we see that
among the terms involving $\SS_i$ or $\PP_j$ the coefficient
$c_{34}$ is present in $((A*B)*C)_{14}$ whereas it
is absent from $(A*(B*C))_{14}$. Then one can easily specialize
the parameters to obtain an example where $((A*B)*C)\ne((A*B)*C)$.
We can also see that $*$ is not {\em diassociative} (i.e. the subloops
generated by two elements are not associative), and hence not a
{\em Moufang loop} like the loops considered by Manin in his book
on cubic forms \cite{Ma}.
\end{noth}

\begin{noth} {\bf Formulas for left and right division.}
Finally, we let $C=A*B$ and give the formulas for $A=C/B$ and
$B=A\backslash C$ for $n=3$. We use the notations $A=M(\bbf{l},\bbf{a})$,
$B=M(\bbf{m},\bbf{b})$, $C=M(\bbf{n},\bbf{c})$.

\noindent The matrix $A=C/B$ is determined by $l_i=n_i-m_i$ for
$i=1,2,3$ and:
\begin{align*}
a_{12} & = \ \frac{c_{12}-u^{n_1-m_1}b_{12}}{u^{m_2}} \\
a_{23} & = \ \frac{c_{23}-u^{n_2-m_2}b_{23}}{u^{m_3}} \\
a_{13} & = \ \frac{c_{13}-u^{n_1-m_1}b_{13}
-\frac{c_{12}-u^{n_1-m_1}b_{12}}{u^{m_2}}b_{23}
-\SS_1(u^{n_1-m_1}b_{12},c_{12}-u^{n_1-m_1}b_{12})}{u^{m_3}}
\end{align*}
\noindent The matrix $B=A\backslash C$ is determined by 
$m_i=n_i-l_i$ for $i=1,2,3$ and:
\begin{align*}
b_{12} & = \ \frac{c_{12}-u^{n_2-l_2}a_{12}}{u^{l_1}} \\
b_{23} & = \ \frac{c_{23}-u^{n_3-l_3}a_{23}}{u^{l_2}} \\
b_{13} & = \ \frac{c_{13}-a_{12}\frac{c_{23}-u^{n_3-l_3}a_{23}}{u^{l_2}}
-u^{n_3-l_3}a_{13}-\SS_1(c_{12}-u^{n_2-l_2}a_{12},u^{n_2-l_2}a_{12})}{u^{l_1}}
\end{align*}
When $C$ is the identity matrix, we see that left inverse and right
inverse coincide.
\end{noth}


\section{Relating lattices and matrices} \label{seclat}

In this section, we consider matrices adapted to well-chosen systems
of generators of lattices. More precisely, we define subsets
$$
\sG_n^\mu((u))\subset \sG_n^d((u))\subset\sG_n^T((u))\subset\sG_n((u))
$$
whose relation to lattices is the following. The set $\sG_n^T((u))$
of {\em $T$-matrices} corresponds to the nice systems of generators of lattices
which we called $T$-bases. The set $\sG_n^d((u))$ of {\em distinguished matrices}
corresponds to the distinguished $T$-bases, that is, to the lattices themselves.
Finally the set $\sG_n^\mu((u))$
of {\em $\mu$-matrices} corresponds to the $\mu$-lattices.
The final result is Theorem~\ref{classification_BK} which formulates
the classification of models of $\mu_{p^n,K}$ in terms of matrices,
well-suited for computations.

From now on, the ring of coefficients $k$ is a perfect field
and $W=W(k)$, $W_n=W_n(k)$.

\subsection{Matrices and lattices}

Recall that lattices, $T$-bases and distinguished bases are defined
in~\ref{ss:lattices_of_Wn}.

\begin{defi} \label{df:T_matrix}
For each $A\in\sG_n((u))$ we consider the column vector
$e_\star=[A]\,p^\star$, its components $e_1,\dots,e_n$, and the
lattice $\fM=\fM(A)$ they generate.
\begin{trivlist}
\itemn{1} We say that $A$ is a {\em $T$-matrix} if
$e_1,\dots,e_n$ is a $T$-basis of $\fM$.
\itemn{2} We say that $A$ is {\em distinguished} if
$e_1,\dots,e_n$ is the distinguished basis of $\fM$.
\end{trivlist}
We denote by $\sG_n^T((u))$, resp. $\sG_n^d((u))$, the set of $T$-matrices,
resp. distinguished matrices, in $\sG_n((u))$.
 We have similar subsets $\sG_n^*[[u]]\subset \sG_n[[u]]$,
$\sG_n^*[u,u^{-1}]\subset \sG_n[u,u^{-1}]$, $\sG_n^*[u]\subset \sG_n[u]$
with $*\in \{T,d\}$.
\end{defi}

Let $\sL_n$ be the set of lattices of $W_n((u))$. We have a
well-defined map
$$
\sG_n((u))\to\sL_n \quad, \quad A\mapsto \fM(A).
$$
Denote by $A(\fM)$ the matrix whose coefficients are the $p$-adic
coefficients of the distinguished basis of $\fM$. Then we have
a section
$$
\sL_n\to\sG^T_n[u,u^{-1}]\subset \sG_n((u)) \quad,\quad
\fM\mapsto A(\fM).
$$

\begin{lemm} \label{lm:T_matrix_T_basis}
Let $A\in\sG_n((u))$ and $\fM=\fM(A)$. Then:
\begin{trivlist}
\itemn{1} $A$ is a $T$-matrix if and only if $\cU A/\cL A\ge 0$,
i.e. $\cU A=B*\cL A$ for some $B\in\sG_n[[u]]$.
\itemn{2} $A\ge 0$ if and only if $\fM\ge 0$.
\end{trivlist}
\end{lemm}

\begin{proo}
(1) Set $e_\star=[A]\,p^\star$.
Due to the shape of matrices in $\sG_n((u))$, we have $v_p(e_i)=i-1$.
It follows from (1) of Lemma~\ref{lemma_T_basis} that $e_\star$ is a
$T$-basis if and only if $pe_i\in\langle e_{i+1},\dots,e_n\rangle$ for
all $i$. This is in turn equivalent to the existence of elements
$b_{ij}\in k[[u]]$ such that
$$
pe_i=[b_{ii}]\,e_{i+1}+\dots+[b_{i,n-1}]\,e_n
$$
for all $i$. Let $B$ be the upper triangular matrix with diagonal
entries $u^{l_i-l_{i+1}}$ and other
entries $b_{ij}\in k[[u]]$. It is simple to see that the set of
equalities above is equivalent to $\cU A=B*\cL A$.

\smallskip

\noindent (2) We have $\fM\ge 0$ if and only if $e_i\in W_n[[u]]$ for all $i$.
Since $e_i=u^{l_i}p^{i-1}+[a_{i,i+1}]\,p^i+\dots+[a_{in}]\,p^{n-1}$,
this means that $u^{l_i}$ and $a_{ij}$ belong to $k[[u]]$ for all $i,j$.
\end{proo}

The construction of the distinguished basis in Lemma~\ref{dist_generators}
shows that the volume of a lattice (defined in \ref{ss:mu_lattices})
can be computed from a $T$-matrix giving rise to it:

\begin{lemm}
For $A\in\sG_n^T((u))$ and $\fM=\fM(A)$, we have $\vol(\fM)=\det(A)$.
\end{lemm}

\begin{proo}
Let $\alpha$ be an integer such that $u^\alpha\fM\subset W_n[[u]]$.
Replacing $\fM$ by $u^\alpha\fM$ and $A$ by $u^\alpha A$, we may assume
that $\alpha=0$. To simplify the notation, we write
$\fM^+=W_n[[u]]/\fM$. Write $A=M(\bbf{l},\bbf{a})$.
We have the following diagram with exact rows and columns:
$$
\xymatrix@R=1.5pc{& 0 \ar[d] & 0 \ar[d] & 0 \ar[d] & \\
0 \ar[r] & \fM[i+1] \ar[r] \ar[d] & \fM[i] \ar[r] \ar[d]
& \fM[i]/\fM[i+1] \ar[r] \ar[d] & 0 \\
0 \ar[r] & W_{n-i}[[u]] \ar[r] \ar[d] & W_{n-i+1}[[u]] \ar[r] \ar[d]
& W_1[[u]] \ar[r] \ar[d] & 0 \\
0 \ar[r] & \fM[i+1]^+ \ar[r] \ar[d] & \fM[i]^+ \ar[r] \ar[d]
& (\fM[i]/\fM[i+1])^+ \ar[r] \ar[d] & 0 \\
& 0 & 0 & 0 & \\}
$$
Since $A\in\sG_n[[u]]$, we have $\fM[i]/\fM[i+1]\simeq u^{l_i}k[[u]]$
and $(\fM[i]/\fM[i+1])^+\simeq k[u]/(u^{l_i})$ of length~$l_i$.
Then the result follows by induction, using the additivity of the
length.
\end{proo}

Let us now look at some natural lattices associated to a lattice $\fM$.
We defined the kernel $\fM[i]$ and the image $\fM(i)$ in \ref{ker_im_p}.
The ring
$W_n((u))$ is endowed with a Frobenius endomorphism $\phi$ whose
restriction to $W_n$ is the Frobenius of the Witt vectors, and such
that $\phi(u)=u^p$. This gives rise to another interesting lattice,
namely the lattice generated by $\phi(\fM)$. Also, for a polynomial
$E(u)\in W_n[u]$ we can consider the lattice $E(u)\fM$.
If $\fM=\fM(A)$, we wish to express the matrices associated to these
lattices in terms of $A$. We will shortly give the result, but we first
need a bit of notation.

\begin{nota} \label{nota:E(u)A}
We denote by $\cP$ the matrix operator taking a square matrix $M$
of size $r$ to the square matrix of size $r+1$ whose upper right block
of size $r$ is $M$ and whose other entries are zero. In pictures,
$$
\cP M=
\left(
\begin{array}{cc}
0 &
\begin{array}{|ccc}
& & \\
& M & \\
& & \\
\hline
\end{array} \\
0 &  0 \\
\end{array}
\right) \ .
$$
The operator $\cP^i$ takes a matrix $M$ of size $r$ to the matrix of size
$r+i$ whose upper right block is $M$ and whose other blocks are zero.
\end{nota}

\begin{defi}
Let $A\in\sG_n((u))$ be a matrix and $E(u)\in W_n[u]$ a polynomial, with
$p$-adic expansion $E(u)=[E_0(u)]+[E_1(u)]\,p+\dots+[E_{n-1}(u)]\,p^{n-1}$.
With the notation of~\ref{nota:E(u)A}, we define:
\begin{trivlist}
\itemn{1} $E(u)\!\diamond\!A=\rho\big(\sum_{i=0}^{n-1}\,[E_i\Id*\cP^i\cU^iA]\big)$,
where $\rho$ is the map from Lemma~\ref{map_rho}.
\itemn{2} $\phi(A)$ is the matrix obtained by applying Frobenius to all the
entries of $A$.
\end{trivlist}
\end{defi}

The two operations $\phi(-)$ and $E(u)\!\diamond\!-$ are compatible
with $\cU$ and $\cL$ in the following sense.

\begin{lemm} \label{lm:U_and_L_commute_with_phi_and_E(u)}
For all matrices $A\in\sG_n((u))$ and polynomials $E(u)\in W_n[u]$, we have:
\begin{trivlist}
\itemn{1} $\cU(\phi(A))=\phi(\cU(A))$ and $\cL(\phi(A))=\phi(\cL(A))$.
\itemn{2} $\cU(E(u)\!\diamond\!A)=E(u)\!\diamond\!\cU(A)$,
and $\cL(E(u)\!\diamond\!A)=E(u)\!\diamond\!\cL(A)$,
\end{trivlist}
where in $E(u)\!\diamond\!\cU(A)$ and $E(u)\!\diamond\!\cL(A)$ it is
the image of $E(u)$ in $W_{n-1}[u]$ that is involved.
\end{lemm}

\begin{proo}
(1) is obvious and we only prove (2). Let $\tau_{\cU}$ be the truncation
map that takes a vector $v$ with $n$ components to the vector whose components
are the first $n-1$ components of $v$, so $\tau_{\cU} p^\star$ is the vector
analogous to $p^\star$ in dimension one less, as in the proof of
Lemma~\ref{homomorphisms_U_and_L}. Since $\cP\cU A$ is the matrix
obtained from $\cU\cP A$ by replacing the last line by $0$, we have:
$[\cP\cU A]\,\tau_{\cU}p^\star = \tau_{\cU} [\cP A]\,p^\star$. It follows
that
$$
[E(u)\!\diamond\!\cU(A)]\,\tau_{\cU}p^\star
=\sum_{i=0}^{n-2} [E_i][\cP^i\cU^{i+1} A]\,\tau_{\cU}p^\star
=\tau_{\cU}\left(\sum_{i=0}^{n-1} [E_i][\cP^i\cU^i A]\,p^\star\right)
=\tau_{\cU}\Big([E(u)\!\diamond\!A]\,p^\star\Big).
$$
But it is exactly the defining property of $M=\cU(E(u)\!\diamond\!A)$ that
$[M]\,\tau_{\cU}p^\star=\tau_{\cU}([E(u)\!\diamond\!A]\,p^\star)$.
This proves that $\cU(E(u)\!\diamond\!A)=E(u)\!\diamond\!\cU(A)$.
The proof for the commutation with $\cL$ is similar: $\cP\cL A$ is the
matrix obtained from $\cL\cP A$ by replacing the first line by $0$, etc.
\end{proo}

\begin{lemm} \label{lemma_matrix_of_M[i]}
Let $A\in\sG_n((u))$ and $\fM\in\sL_n$.
\begin{trivlist}
\itemn{1} If $\fM=\fM(A)$ then:
\begin{trivlist}
\itemm{a} $\fM(i)=\fM(\cU^{i-1}A)$,
\itemm{b} $\fM[i]=\fM(\cL^{i-1}A)$,
\itemm{c} $\langle\phi(\fM)\rangle=\fM(\phi(A))$,
\itemm{d} $E(u)\fM=\fM(E(u)\!\diamond\!A)$.
\end{trivlist}
\itemn{2} If $A$ is a $T$-matrix then
$\cU^{i-1}A$, $\cL^{i-1}A$, $\phi(A)$, $E(u)\!\diamond\!A$ are also $T$-matrices.
\itemn{3} If $A$ is distinguished then $\cU^{i-1}A$, $\cL^{i-1}A$,
$\phi(A)$ are also distinguished.
\end{trivlist}
\end{lemm}

It is not true in general that if $A$ is distinguished then
$E(u)\!\diamond\!A$ is distinguished. There are obvious counter-examples
for $n=2$ as soon as $l_1\ge l_2+1$.

\begin{proo}
(1) Let $e_\star=[A]\,p^\star$. Let us fix $i\in\{1,\dots,n\}$ and define:
\begin{trivlist}
\itemm{a} $f_j=p^{i-1}e_j$ for $1\le j\le n+1-i$,
\itemm{b} $g_j=e_{j+i-1}$ for $1\le j\le n+1-i$,
\itemm{c} $h_j=\phi(e_j)$ for $1\le j\le n$,
\itemm{d} $\ell_j=E(u)e_j$  for $1\le j\le n$,
\end{trivlist}
The elements $f_j$ generate $\fM(i)$ and we have $f_\star=[\cU^{i-1}A]\,p^\star$,
hence $\fM(i)=\fM(\cU^{i-1}A)$. The elements $g_j$ generate $\fM[i]$ and satisfy
$g_\star=[\cL^{i-1}A]\,p^\star$, so that $\fM[i]=\fM(\cL^{i-1}A)$. The elements $h_j$
generate $\langle\phi(\fM)\rangle$ and satisfy
$h_\star=[\phi(A)]\,p^\star$ so
$\langle\phi(\fM)\rangle=\fM(\phi(A))$. Finally the elements $\ell_j$
generate $E(u)\fM$ and moreover a simple matrix computation shows that
$p^i[A]\,p^\star=[\cP^i\cU^iA]\,p^\star$ so
$$
E(u)e_\star=\big(\sum\,[E_i]p^i\big)[A]\,p^\star
=\sum [E_i][\cP^i\cU^iA]\,p^\star
=\big[\rho\big(\sum\,E_i\Id*\cP^i\cU^iA\big)\big]\,p^\star
=[E(u)\!\diamond\!A]\,p^\star \ .
$$
It follows that $E(u)\fM=\fM(E(u)\!\diamond\!A)$.

\smallskip

\noindent (2) Using the characterization 1) in Lemma~\ref{lemma_T_basis}, it
is very easy to prove that $f_\star$, $g_\star$, $h_\star$, $\ell_\star$
are $T$-bases.

\smallskip

\noindent (3) It is immediate that the matrices $ \cU^{i-1}A$, $\cL^{i-1}A$ and
$\phi(A)$ have Laurent polynomial entries and satisfy the condition
on the degrees required to be distinguished.
\end{proo}

\begin{lemm} \label{lemma_inclusion_of_lattices}
Let $A,A'$ be in $\sG_n((u))$ and $\fM=\fM(A)$, $\fM'=\fM(A')$.
\begin{trivlist}
\itemn{1} Assume that $A'\in\sG_n^T((u))$. Then $\fM\subset\fM'$
if and only if $A/A'\ge 0$.
\itemn{2} In particular, the $T$-matrices are the minimal elements
among the matrices $A\in\sG_n((u))$ such that $\fM(A)=\fM$, in the sense
that for any two matrices $A,A'$ with $\fM(A)=\fM(A')=\fM$,
if $A'$ is a $T$-matrix then $A/A'\ge 0$.
\itemn{3} Assume that $A,A'\in\sG_n^T((u))$. Then $\fM=\fM'$ if and only
if $A/A'$ is positive and unipotent.
\end{trivlist}
\end{lemm}

\begin{proo}
Let $e_\star=[A]\,p^\star$ and $e'_\star=[A']\,p^\star$ be the associated
generating sets. Then $\fM\subset\fM'$ if and only if for each $i$
we have $e_i\in\fM'[i]$. This means that there exist scalars
$b_{ij}\in k[[u]]$ such that
$$
e_i=[b_{ij}]\,e'_i+[b_{i,i+1}]\,e'_{i+1}+\dots+[b_{i,n}]\,e_n \ .
$$
Let $B$ be the upper triangular matrix with coefficients $b_{ij}$.
These equalities amount to $e_\star=[B]\,e'_\star$, in other words
$[A]\,p^\star=[B][A']\,p^\star=[B*A']\,p^\star$. Thus $A/A'=B\ge 0$
and this proves (1). Now (2) and (3) follow immediately.
\end{proo}

\begin{rema} \label{rm:order}
It follows from this lemma that the relation $\succ$ on $\sG_n^T((u))$
defined by $A\succ B$ if and only if $A/B\ge 0$
is reflexive and transitive.
\end{rema}

\subsection{Matricial description of Breuil-Kisin modules}\label{MatrixBK}

Finally we arrive at the description in terms of matrices of
the Breuil-Kisin modules corresponding to a group scheme which is a
model of $\mu_{p^n}$. We recall that $K$ is a finite totally ramified
field extension of $K_0$, the fraction field of the Witt ring $W=W(k)$ of
a perfect field $k$ of characteristic $p>0$, and that $E(u)$ is the
Eisenstein polynomial of a fixed uniformizer $\pi\in \cO_K$.

\begin{defi} \label{df:distinguished_matrix}
We say that $A=M(\bbf{l},\bbf{a})\in\sG_n((u))$ is a {\em $\mu$-matrix} if it is
distinguished and if
\begin{trivlist}
\itemn{1} $\phi(A)/A\ge 0$,
\itemn{2} $(E(u)\!\diamond\!A)/\phi(A)\ge 0$.
\end{trivlist}
We denote by $\sG_n^\mu((u))$ the set of $\mu$-matrices in $\sG_n((u))$.
\end{defi}

With the induced order of $\sG_n^T((u))$ ({\em cf} Remark~\ref{rm:order}),
the set $\sG_n^\mu((u))$ is an ordered set.
Since $\cU$ and $\cL$ are loop homomorphisms (\ref{homomorphisms_U_and_L}),
take positive matrices to positive matrices (obvious), and commute with
$\phi$ and $E(u)\!\diamond\!-$ (\ref{lm:U_and_L_commute_with_phi_and_E(u)}),
one sees that if $A\in\sG_n((u))$ is a $\mu$-matrix, then $\cU A$ and
$\cL A$ are also $\mu$-matrices.

\begin{theo} \label{classification_BK}
The maps $G\mapsto \fM(G)$ and $\fM\mapsto A(\fM)$ give
bijections between:
\begin{itemize}
\item the set of isomorphism classes of $R$-models of $\mu_{p^n,K}$,
\item the set $\sL_n^\mu$ of $\mu$-lattices, i.e. finitely generated
sub-$W_n[[u]]$-modules of $W_n((u))$ satisfying
$E(u)\fM\subset \langle\phi(\fM)\rangle\subset\fM$,
\item the set $\sG_n^\mu((u))$ of $\mu$-matrices, i.e. matrices
\end{itemize}
$$
A=\left(
\begin{array}{ccccc}
u^{l_1} & a_{12} & a_{13} & \dots & a_{1n} \\
& u^{l_2} & a_{23} & & a_{2n} \\
& & \ddots & \ddots & \vdots \\
& & & u^{l_{n-1}} & a_{n-1,n} \\
0 & & & & u^{l_n} \\
\end{array}
\right)
$$
where $l=(l_1,\dots,l_n)\in\NN ^n$ and $a_{ij}\in k[u]$ for all $i,j$,
such that:
\begin{trivlist}
\itemm{1} $\deg_u(a_{ij})<l_j$ whenever $1\le i<j\le n$,
\itemm{2} $\cU A/\cL A\ge 0$,
\itemm{3} $\phi(A)/A\ge 0$,
\itemm{4} $(E(u)\!\diamond\!A)/\phi(A)\ge 0$.
\end{trivlist}
These bijections are increasing: if $G,G'$ are models of $\mu_{p^n,K}$
with associated lattices $\fM,\fM'$ and distinguished matrices $A,A'$,
then the following conditions are equivalent:
\begin{itemize}
\item there exists a model map $G\to G'$,
\item $\fM\subset \fM'$,
\item $A/A'\ge 0$.
\end{itemize}
Finally, these bijections are "compatible with quotients and kernels":
\begin{itemize}
\item $G(i)$, $\fM(i)$ and $\cU^{i-1}A$ correspond to each other, and
\item $G[i]$, $\fM[i]$ and $\cL^{i-1}A$ correspond to each other.
\end{itemize}
\end{theo}

\begin{proo}
The increasing bijection between models of $\mu_{p^n}$ and $\mu$-lattices
is the Breuil-Kisin equivalence. The map $\fM\mapsto A(\fM)$ is the map
taking a lattice to its distinguished matrix, so that $A=A(\fM)=M(\bbf{l},\bbf{a})$
satisfies the conditions (1) and (2).
It remains to prove that the additional conditions satisfied by a $\mu$-lattice
translate into the additional conditions (3) and (4) in the theorem.
Indeed, the condition (3) is a translation of the fact that
$\langle\phi(\fM)\rangle\subset\fM$ and the condition (4) is a translation
of the fact that $E(u)\fM\subset \fM$. Moreover, since $\fM$ is positive
(see~\ref{ss:mu_lattices}), then so is~$A$ and hence $l_i\ge 0$. This gives
the refinement in the statement of the theorem. The fact that the bijection
between $\mu$-lattices and $\mu$-matrices is increasing is
Lemma~\ref{lemma_inclusion_of_lattices}. The fact that the bijections are
compatible with quotients and kernels comes from Lemma~\ref{lemma_matrix_of_M[i]}
and the fact that $\cU$ and $\cL$ preserve $\mu$-matrices.
\end{proo}

\begin{rema} \label{rm:recap}
Let us recapitulate some of the information we have on the parameters.
\begin{trivlist}
\itemn{1} We have $l_1\ge l_2\ge\cdots\ge l_n$ since $\cU A/\cL A\ge 0$
($A$ is a $T$-matrix). In fact, the positivity of $\cU A/\cL A$ corresponds
to the existence of the model maps $G(i)\to G[i]$ of~\ref{k_i_p}, for all $i$.
\itemn{2} We have $l_i\ge 0$ and $\val_u(a_{i,i+1})\ge l_{i+1}/p$, for all $i$.
Indeed, since $\phi(A)/A\ge 0$ there exists a positive matrix $B=M(\bbf{m},\bbf{b})$
such that $\phi(A)=B*A$. Comparing the diagonal entries, we get
$(p-1)l_i=m_i\ge 0$ thus $l_i\ge 0$. Comparing the entries at distance $1$
from the diagonal, we get
$(a_{i,i+1})^p=u^{(p-1)l_i}a_{i,i+1}+u^{l_{i+1}}b_{i,i+1}$.
Thus $(a_{i,i+1})^p\equiv 0 \mod u^{l_{i+1}}$.
\itemn{3} We have $e/(p-1)\ge l_1$ since $(E(u)\!\diamond\!A)/\phi(A)\ge 0$.
Indeed, the upper left entry of $E(u)\!\diamond\!A$ is $u^{e+l_1}$ and
the upper left entry of $\phi(A)$ is $u^{pl_1}$. The result follows.
\end{trivlist}
\end{rema}

Theorem \ref{classification_BK} gives already very precise information
on the structure of the set of models of $\mu_{p^n}$. In a naive way,
it is parametrized by $n$ integers
$0\leq l_n\leq \cdots\leq l_1\leq e/(p-1)$ and at most
$\sum_{i=1}^{n-1}il_i$ elements of $k$ (the coefficients $a_{ij}$),
as follows from condition~(1) in Theorem~\ref{classification_BK}.
%

\begin{defi}
The parameters $(l_1,\ldots,l_n)$ of a model of $\mu_{p^n}$ are called
the {\em type} of the model.
\end{defi}

The geometric interpretation of the type of a model of $\mu_{p^n}$
is quite clear. Theorem~\ref{classification_BK} gives a precise
geometric interpretation for the other (somehow more mysterious)
parameters of the Breuil-Kisin modules: some of
them parametrize flat subgroup schemes or quotients, and some others
parametrize extensions between such subquotients, models of $\mu_{p^s}$ and
$\mu_{p^r}$ for $1\leq r,s\leq {n-1}$.

\begin{rema}
A remaining open question is the structure of this set of parameters.
The explicit computation of relations is completed for $n=3$ in
Section~\ref{BKn3}. Since the functions $\SS_i$ and $\PP_i$ involved
in the operation $A*B$ are defined by exponentiation with respect to
negative powers of $p$, a high enough power of Frobenius transforms
the constraints defining $\mu$-matrices into polynomial relations between
the coefficients of the $a_{ij}$. Hence up to Frobenius, we can easily
define the variety of models of $\mu_{p^n}$. The study of the dimension
and irreducible components of this variety has to be compared to the
works of Imai and Caruso (\cite{Car}, \cite{Im})
on Kisin's moduli space of models of $\mu_{p^n}$ (\cite{Ki2}).
\end{rema}


\section{Computation of $\mu$-matrices for $n=3$} \label{BKn3}

Since the bijections in Theorem~\ref{classification_BK} are compatible
with quotients and kernels, the matricial formulas for the models of
$\mu_{p^n}$ contain the matricial formulas for the models of $\mu_{p^i}$
for all $i\le n$. In this section, we work out the conditions in
Theorem~\ref{classification_BK} for $n= 3$ and $p\geq 3$. We stress that they include also the case $n=1,2$. And in these cases one gets the formulas obtained by Caruso in \cite{To2}, Appendix A.

\subsection{Computation of the matrices}
We have
$$
A=\left(
\renewcommand{\arraystretch}{1.5}
\begin{array}{ccc}
u^{l_1} & a_{12} & a_{13} \\
   0   & u^{l_2} & a_{23} \\
   0   &   0     & u^{l_3} \\
\end{array}
\right)
\quad,\quad
\cU A=\left(
\renewcommand{\arraystretch}{1.5}
\begin{array}{cc}
u^{l_1} & a_{12} \\
  0    & u^{l_2} \\
\end{array}
\right)
\quad,\quad
\cL A=\left(
\renewcommand{\arraystretch}{1.5}
\begin{array}{cc}
u^{l_2} & a_{23} \\
  0    & u^{l_3} \\
\end{array}
\right) .
$$
Using Examples \ref{examples_for_small_n} we find
$$
\cU A/\cL A=\left(
\renewcommand{\arraystretch}{1.5}
\begin{array}{cc}
u^{l_1-l_2} & \frac{a_{12}-u^{l_1-l_2}a_{23}}{u^{l_3}} \\
  0    & u^{l_2-l_3} \\
\end{array}
\right) .
$$
Moreover we have
$$
\phi(A)=\left(
\renewcommand{\arraystretch}{1.5}
\begin{array}{ccc}
u^{pl_1} & a_{12}^p & a_{13}^p \\
   0   & u^{pl_2} & a_{23}^p \\
   0   &   0     & u^{pl_3} \\
\end{array}
\right)
$$
and
$$
\phi(A)/A=\left(
\renewcommand{\arraystretch}{1.5}
\begin{array}{ccc}
u^{(p-1)l_1} & \frac{a_{12}^p-u^{(p-1)l_1}a_{12}}{u^{l_2}} & p_{13} \\
   0   & u^{(p-1)l_2} & \frac{a_{23}^p-u^{(p-1)l_2}a_{23}}{u^{l_3}} \\
   0   &   0     & u^{(p-1)l_3} \\
\end{array}
\right)
$$
where
$$
p_{13}=\frac{a_{13}^p-u^{(p-1)l_1}a_{13}
-\frac{a_{12}^p-u^{(p-1)l_1}a_{12}}{u^{l_2}}a_{23}-\SS_1(u^{(p-1)l_1}a_{12},-a_{12}^p)}{u^{l_3}} .
$$
Finally we compute $E(u)\!\diamond\!A$. Note that for $n=3$ we have
$E_i\Id*\cP^i\cU^iA=E_i\cP^i\cU^iA$ for all~$i$, but this is false already
for $n=4$, because of the failure of multiplicativity of Teichm\"uller
representatives of polynomials (see Lemma~\ref{lm:easy_formulas}). Thus
$$
E(u)\!\diamond\!A=\left(
\renewcommand{\arraystretch}{1.5}
\begin{array}{ccc}
u^{e+l_1} & u^ea_{12}+u^{l_1}E_1
& u^ea_{13}+a_{12}E_1+\SS_1(u^ea_{12},u^{l_1}E_1)+u^{l_1}E_2 \\
   0   & u^{e+l_2} & u^ea_{23}+u^{l_2}E_1 \\
   0   &   0     & u^{e+l_3} \\
\end{array}
\right) .
$$
and
$$
(E(u)\!\diamond\!A)/\phi(A)=\left(
\renewcommand{\arraystretch}{1.5}
\begin{array}{ccc}
u^{e-(p-1)l_1} & \frac{u^ea_{12}+u^{l_1}E_1-u^{e-(p-1)l_1}a_{12}^p}{u^{pl_2}}
& q_{13} \\
0   & u^{e-(p-1)l_2}
& \frac{u^ea_{23}+u^{l_2}E_1-u^{e-(p-1)l_2}a_{23}^p}{u^{pl_3}} \\
0   &   0     & u^{e-(p-1)l_3} \\
\end{array}
\right)
$$
where
\begin{align*}
q_{13} = & \ \frac{u^ea_{13}+a_{12}E_1+\SS_1(u^ea_{12},u^{l_1}E_1)+u^{l_1}E_2
-u^{e-(p-1)l_1}a_{13}^p}{u^{pl_3}} \\
& \\
& \ -\frac{\frac{u^ea_{12}+u^{l_1}E_1-u^{e-(p-1)l_1}a_{12}^p}{u^{pl_2}}a_{23}^p-
\SS_1(u^{e-(p-1)l_1}a_{12}^p,u^ea_{12}+u^{l_1}E_1-u^{e-(p-1)l_1}a_{12}^p)}{u^{pl_3}} \ .
\end{align*}

\subsection{Translation of the conditions of the theorem}

\medskip

\noindent\textbullet\ Condition (1) yields:
$$
\boxed{\deg_u(a_{12})\le l_2-1}
\quad,\quad
\boxed{\deg_u(a_{13})\le l_3-1}
\quad\mbox{and}\quad
\boxed{\deg_u(a_{23})\le l_3-1} \ .
$$

\medskip

\noindent\textbullet\ Condition (2) yields:
$$
\boxed{l_1\ge l_2\ge l_3}
\quad\mbox{and}\quad
\boxed{a_{12}-u^{l_1-l_2}a_{23}\equiv 0 \mod u^{l_3}} \ .
$$

\medskip

\noindent\textbullet\ Condition (3) yields:
$$a_{12}^p-u^{(p-1)l_1}a_{12}\equiv 0 \mod u^{l_2} \ , $$

$$a_{23}^p-u^{(p-1)l_2}a_{23}\equiv 0 \mod u^{l_3} $$

\noindent and

$$
a_{13}^p-u^{(p-1)l_1}a_{13}-\frac{a_{12}^p-u^{(p-1)l_1}a_{12}}{u^{l_2}}a_{23}
-\SS_1(u^{(p-1)l_1}a_{12},-a_{12}^p) \equiv 0 \mod u^{l_3} \ .
$$
Since $(p-1)l_1\ge l_2$, the first two are equivalent to:
$$
\boxed{a_{12}^p\equiv 0 \mod u^{l_2}}
\quad\mbox{and}\quad
\boxed{a_{23}^p\equiv 0 \mod u^{l_3}} \ .
$$
Concerning the third, observe that since $(p-1)l_1\ge l_3$
the term $u^{(p-1)l_1}a_{13}$ can be neglected. Also since
$\val_u(\SS_1(x,y))\ge \max(\val_u(x),\val_u(y))$ by
Lemma~\ref{lm:easy_formulas}, we see that the $\SS_1$ term can be neglected.
Finally the term $u^{(p-1)l_1-l_2}a_{12}a_{23}$ can also be neglected:
indeed $pl_1\ge 2l_1\ge l_2+l_3$ implies that its valuation is at
least
$$
\big((p-1)l_1-l_2\big)+\frac{l_2}{p}+\frac{l_3}{p}
=\frac{1}{p}\big((p-1)(pl_1-l_2)+l_3\big)
\ge\frac{1}{p}\big((p-1)l_3+l_3\big)=l_3 \ .
$$
So the third condition is equivalent to:
$$
\boxed{a_{13}^p-u^{-l_2}a_{12}^pa_{23}\equiv 0 \mod u^{l_3}} \ .
$$

\medskip

\noindent\textbullet\ Condition (4) yields:
$$
\boxed{u^ea_{12}+u^{l_1}E_1-u^{e-(p-1)l_1}a_{12}^p \equiv 0 \mod u^{pl_2}} \ ,
$$

$$
\boxed{u^ea_{23}+u^{l_2}E_1-u^{e-(p-1)l_2}a_{23}^p \equiv 0 \mod u^{pl_3}}
$$

\noindent and

$$
\begin{array}{cl}
u^ea_{13}+a_{12}E_1+\SS_1(u^ea_{12},u^{l_1}E_1)+u^{l_1}E_2
-u^{e-(p-1)l_1}a_{13}^p & \\ & \\
-\frac{u^ea_{12}+u^{l_1}E_1-u^{e-(p-1)l_1}a_{12}^p}{u^{pl_2}}a_{23}^p
-\SS_1(u^{e-(p-1)l_1}a_{12}^p,-u^ea_{12}-u^{l_1}E_1) & \equiv 0 \mod u^{pl_3} \ . \\
\end{array}
$$
Finally the last but one boxed congruence implies that
$$
\SS_1(u^{e-(p-1)l_1}a_{12}^p,u^ea_{12}+u^{l_1}E_1-u^{e-(p-1)l_1}a_{12}^p)
\equiv 0 \mod u^{pl_2}.
$$
Hence it vanishes also modulo $u^{pl_3}$ and we obtain:
$$
\boxed{
\begin{array}{rr}
u^ea_{13}+a_{12}E_1+\SS_1(u^ea_{12},u^{l_1}E_1)+u^{l_1}E_2
-u^{e-(p-1)l_1}a_{13}^p & \\ & \\
-\frac{u^ea_{12}+u^{l_1}E_1-u^{e-(p-1)l_1}a_{12}^p}{u^{pl_2}}a_{23}^p
& \equiv 0 \mod u^{pl_3} \\
\end{array}} \ .
$$

\begin{coro}\label{coro1}
Let $p\geq 3$. Let $\fM\in(\Mod/\fS)_3$ be the Breuil-Kisin module
of a finite flat $R$-model of $\mu_{p^3,K}$.
Then there exists a unique family of parameters
$(l_1,l_2,l_3,a_{12},a_{13},a_{23})$ composed of three integers
$0\leq l_3\leq l_2\leq l_1\leq e/(p-1)$ and three polynomials
$a_{12},a_{13},a_{23}\in k[u]$ satisfying:
\begin{trivlist}
\itemn{i}
$\deg _u a_{12}\leq l_2-1$, $\deg _u a_{13}\leq l_3-1$, $\deg _u a_{23}\leq l_3-1$,
\medskip
\itemn{ii}
$a_{12}-u^{l_1-l_2}a_{23}\equiv 0\mod u^{l_3}$, $a_{12}^p
\equiv 0\mod u^{l_2}$, $a_{23}^p\equiv 0\mod u^{l_3}$,
\medskip
\itemn{iii}
$a_{13}^p-u^{-l_2}a_{12}^pa_{23}\equiv 0\mod u^{l_3}$,
\medskip
\itemn{iv}
$u^ea_{12}+u^{l_1}E_1-u^{e-(p-1)l_1}a_{12}^p \equiv 0 \mod u^{pl_2}$ and
$u^ea_{23}+u^{l_2}E_1-u^{e-(p-1)l_2}a_{23}^p \equiv 0 \mod u^{pl_3}$,
\medskip
\itemn{v}
$\begin{array}[t]{l}
u^ea_{13}+a_{12}E_1+\SS_1(u^ea_{12},u^{l_1}E_1)+ \\
+ u^{l_1}E_2
-u^{e-(p-1)l_1}a_{13}^p
-\frac{u^ea_{12}+u^{l_1}E_1-u^{e-(p-1)l_1}a_{12}^p}{u^{pl_2}}a_{23}^p
\equiv 0 \mod u^{pl_3},
\end{array}$
\end{trivlist}
such that $\fM=\fM(A)$ with
$$A=\left(
\renewcommand{\arraystretch}{1.5}
\begin{array}{ccccc}
u^{l_1} & [a_{12}] & [a_{13}] \\
0& u^{l_2} & [a_{23}] \\
0&0 & u^{l_3} \\
\end{array}
\right).
$$
\end{coro}

\subsection{The tamely ramified case}

In the tamely ramified case $(e,p)=1$, some of these congruences can
be simplified. To begin with, let us prove that
$$
\boxed{l_1\ge pl_2} \quad\mbox{and}\quad \boxed{l_2\ge pl_3} \ .
$$
Let us prove the first inequality. If $l_2=0$ there is nothing to show.
Otherwise we have $l_2>0$ and we claim that the only monomial of
degree $l_1$ in the polynomial
$$
u^ea_{12}+u^{l_1}E_1-u^{e-(p-1)l_1}a_{12}^p
$$
is $u^{l_1}E_1(0)$. Indeed the first term has valuation
$$
\val(u^ea_{12})\ge e+l_2/p>e\ge (p-1)l_1\ge l_1 \ .
$$
Moreover since $a_{12}^p$ is a $p$-th power, the degrees of the monomials
of $u^{e-(p-1)l_1}a_{12}^p$ are of the form
$$
e-(p-1)l_1+ip=e+l_1-p(l_1-i)
$$
for some integer $i$. Since $(e,p)=1$, this degree is not congruent to
$l_1$ modulo $p$. This proves that $u^{l_1}E_1(0)$ is the only monomial
of degree $l_1$ and then the congruence
$$
u^ea_{12}+u^{l_1}E_1-u^{e-(p-1)l_1}a_{12}^p \equiv 0 \mod u^{pl_2}
$$
forces $l_1\ge pl_2$. The proof that $l_2\ge pl_3$ is similar.

It follows that the condition given by the congruence
$a_{12}-u^{l_1-l_2}a_{23}\equiv 0 \mod u^{l_3}$ is empty since we
already know that both terms have valuation at least $l_3$.

It follows also that the congruences implied by condition (4) become:
$$
\boxed{u^{e-(p-1)l_1}a_{12}^p \equiv 0 \mod u^{pl_2}}
\quad,\quad
\boxed{u^{e-(p-1)l_2}a_{23}^p \equiv 0 \mod u^{pl_3}}
$$
and
$$
\boxed{
a_{12}E_1-u^{e-(p-1)l_1}a_{13}^p
-\frac{u^ea_{12}+u^{l_1}E_1-u^{e-(p-1)l_1}a_{12}^p}{u^{pl_2}}a_{23}^p
\equiv 0 \mod u^{pl_3}} \ .
$$
Then, in the tamely ramified case, the parametrisation of models of
$\mu_{p^3,K}$ is much easier:

\begin{coro}
In the tamely ramified case $(e,p)=1$, the models of $\mu_{p^3}$ over $\cO_K$ are
classified by three integers $0\leq p^2l_3\leq pl_2\leq l_1\leq e/(p-1)$ and three
polynomials $a_{12},a_{13},a_{23}\in k[u]$ satisfying:
\begin{trivlist}
\itemn{i}
$\deg _u a_{12}\leq l_2-1$, $\deg _u a_{13}\leq l_3-1$, $\deg _u a_{23}\leq l_3-1$,
\medskip
\itemn{ii}
$u^{e-(p-1)l_1}a_{12}^p \equiv 0 \mod u^{pl_2}$, $u^{e-(p-1)l_2}a_{23}^p
\equiv 0 \mod u^{pl_3}$,
\medskip
\itemn{iii}
$a_{12}E_1-u^{e-(p-1)l_1}a_{13}^p
-\frac{u^ea_{12}+u^{l_1}E_1-u^{e-(p-1)l_1}a_{12}^p}{u^{pl_2}}a_{23}^p
\equiv 0 \mod u^{pl_3}$.
\end{trivlist}
\end{coro}


\begin{rema}
The tamely ramified case seems to be easy to compute in higher dimension.
In Corollary \ref{coro1}, even for $n=3$, we can see that ramification
intervenes in the computation in a very delicate way: not only through
the coefficient $E_2$ of the Eisenstein polynomial
but also through the lifting modulo $p^2$ of the parameters via
$\SS_1(u^ea_{12},u^{l_1}E_1)$.
\end{rema}


\section{Sekiguchi-Suwa Theory}\label{secSS}

In this section, we recall and complement some aspects of Sekiguchi-Suwa
Theory. The main definitions and results are given in Subsections
\ref{ss:some_defs}, \ref{ss:DAHe}, \ref{ss:main_thms}. For an extended
version, see \cite{MRT}. We also give an interpretation of these results
from a matricial point of view: we introduce the set $\sM_n$ of
matrices parametrizing filtered group schemes, and study its basic
properties. This is the topic of Subsection~\ref{sub:SS with matrices}.

\subsection{Some definitions about Witt vectors} \label{ss:some_defs}

\begin{noth} {\bf The maps $V,F,T$.}
We recall here some definitions  about Witt vectors. We emphasize that in
contrast with Sections \S~\ref{secBK} to \S~\ref{BKn3}, we consider Witt vectors
with coefficients in an arbitrary ring,
not necessarily perfect of characteristic $p$. In particular, we need to
consider quotients of a discrete valuation ring of unequal characteristics.
For $r\ge 0$, we recall the definition of the $r$-th Witt polynomial:
$$
\Phi_r(X_0,\dots,X_r)=X_0^{p^r}+pX_{1}^{p^{r-1}}+\dots+p^rX_r.
$$
Then for each ring $A$ the following maps are defined:
\begin{itemize}
    \item[-]Verschiebung:
    \begin{align*}
    V:W(A)&\too W(A)\\
(a_0,a_1,a_2,\dots)&\longmapsto (0,a_0,a_1,a_2,\dots)
    \end{align*}
    \item[-]Frobenius:
 \begin{align*}
    F:W(A)&\too W(A)\\
\bbf{a}=(a_0,a_1,a_2,\dots)&\longmapsto
(F_0(\bbf{a}),F_1(\bbf{a}),F_2(\bbf{a}),\dots)
    \end{align*}
where the polynomials $F_r(\XX)=F_r(X_0,\dots,X_r)\in
\ZZ[X_0,\dots,X_{r+1}]$ are defined inductively by
$$
\Phi_r(F_0(\XX),F_1(\XX),\dots,F_r(\XX))=\Phi_{r+1}(X_0,\dots,X_{r+1}).
$$
    \item[-]$T$ map:
 \begin{align*}
    T:W(A)\times W(A) & \too W(A)\\
(\bbf{a},\bbf{x})&\longmapsto
T_{\bbf{a}}\bbf{x}=
(T_0(\bbf{a},\bbf{x}),T_1(\bbf{a},\bbf{x}),T_2(\bbf{a},\bbf{x}),\dots)
 \end{align*}
where the polynomials
$T_r(\AA,\XX)=T_r(A_0,\dots,A_r,X_0,\dots,X_r)
\in\ZZ[A_0,\dots,A_r,X_0,\dots,X_r]$ are defined inductively by
$$
\Phi_n(T_0,\dots,T_n)=\sum_{i=0}^{n}p^{n-i}(A_{n-i})^{p^i}
\Phi_i(X_0,\dots,X_i).
$$
\end{itemize}
Since $\Phi_n(T_0,\dots,T_n)$ is linear in the variables $\Phi_i(X_0,\dots,X_i)$,
we see that for fixed $\bbf{a}$ the map $T_{\bbf{a}}$ is a morphism of additive
groups. Moreover, it is easy to see that for any ring $A$ and Witt vectors
$\bbf{a},\bbf{x}\in W(A)$ with $\bbf{a}=(a_0,\dots,a_n,\dots)$
we have explicitly
$T_{\bbf{a}}\bbf{x}=\sum_{k=0}^{\infty}V^k([a_k]\bbf{x})$ (see
\cite{SS1}, Lemma 4.2). For instance if $\bbf{a}=[a_0]$ is a Teichm\"uller element
then $T_{\bbf{a}}$ is nothing else than left multiplication by $[a_0]$, and
in particular $T_1$ is the identity. If $\bbf{x}=[x_0]$ is Teichm\"uller then
$T_{\bbf{a}}([x_0])=(a_0x_0,a_1x_0,a_2x_0,\dots)$.
For each ring $A$ an element $\lambda\in A$, we set
$$
\lambda.\bbf{a}\egaldef(\lambda a_0,\lambda a_1,\lambda a_2,\dots)
=T_{\bbf{a}}([\lambda]).
$$
Clearly $\lambda_1.(\lambda_2.\bbf{a})=(\lambda_1\lambda_2).\bbf{a}$
which will usually be written $\lambda_1\lambda_2.\bbf{a}$. The ideal
$\lambda.W(A)$ is the kernel of the morphism of
rings $W(A)\to W(A/\lambda A)$. If two vectors $\bbf{a},\bbf{b}$
are congruent modulo $\lambda.W(A)$, we sometimes write simply
$\bbf{a}\equiv\bbf{b}\mod \lambda$. We will also have to consider the
following type of Witt vectors with coefficients in the ring $A[1/\lambda]$:
$$
\frac{\bbf{a}}{\lambda}\egaldef
\left(\frac{a_0}{\lambda},\frac{a_1}{\lambda},\frac{a_2}{\lambda},\dots\right).
$$
The notations $\frac{1}{\lambda}\bbf{a}$ or $\bbf{a}/\lambda$ may also
be used when it is convenient.

\begin{defi}
For any ring $A$, we define the {\em subfunctor of finite Witt vectors} by
$$
W^f(A)=\big\{(a_0,a_1,a_2,\dots)\in W(A)\;;
\text{ $a_i=0$ for $i\gg 0$}\big\}
$$
and the {\em completion} of $W(A)$ by
$$
\widehat{W}(A)=\big\{(a_0,a_1,a_2,\dots)\in W(A)\;;
\text{ $a_i=0$ for  $i\gg 0$ and $a_i$ is nilpotent for all $i$}\big\}.
$$
\end{defi}

Note that $W^f(A)$ is not a subgroup of $W(A)$, but $\widehat{W}(A)$ is an ideal
in $W(A)$ which is stable under $F$ and $V$
(see \cite{MRT}, 2.2.1, 2.2.3, 2.2.4).
\end{noth}

\begin{noth} {\bf The $T$-multiplication.}
We shall define a new product between matrices whose entries are Witt
vectors. We need to start with some elementary properties of the map $T$
when one of the variables is fixed.

\begin{lemm}\label{lem:Ta injective}
Let $A$ be a ring and $\bbf{a}=(a_0,a_1,a_2,\dots)\in W(A)$  with $a_0$ not
a zero divisor. Then the morphism $T_{\bbf{a}}$ is injective. If $a_0$ is
invertible then it is an isomorphism.
\end{lemm}

\begin{proo}
Let us suppose that $a_0$ is not a zero divisor and that
$T_{\bbf{a}}\bbf{x}=0$ with $\bbf{x}=(x_0,x_1,\dots)\in W(A)$.
We prove, by induction, that $x_n=0$ for any $n$. Since
$\Phi_0(T_{\bbf{a}}\bbf{x})=a_0x_0=0$ and since $a_0$ is not a
zero divisor then  $x_0=0$. We now suppose that $x_i=0$
for $i\le n$. This means that $\bbf{x}=V^{n+1}\bbf{y}$ with
$\bbf{y}=(x_{n+1},\dots,x_r,\dots)\in W(A)$. Therefore
$$
T_\bbf{a}\bbf{x}=\sum_{k=0}^{\infty}V^k([a_k]V^{n+1}(\bbf{y}))=
V^{n+1}(\sum_{k=0}^{\infty}V^k([a_k^{p^{n+1}}]\bbf{y}))=0.
$$
In particular we have $a_0^{p^{n+1}} x_{n+1}=0$.
Since $a_0$ is not a zero divisor then
 $x_{n+1}=0$.

Let us now suppose that $a_0$ is invertible. Let $\bbf{y}\in W(A)$.  Let
$p_n:W(A)\to W_n(A)$ and $p_{n,k}: W_{n}(A)\to W_k(A)$, if $n\ge k$, the
natural projections. We now prove that for any $n\in \NN$ there exist
$\bbf{x}_n\in W_n(A)$ such that $T_{p_n(\bbf{a})}\bbf{x}_n=p_n(\bbf{y})$ and
$p_{n,n-1}(\bbf{x}_n)=\bbf{x}_{n-1}$. This clearly implies that there
exists $\bbf{x}\in W(A)$ such that $T_{\bbf{a}}\bbf{x}=\bbf{y}$.

We prove the above statement by induction. Clearly
$\bbf{x}_0=(x_0)=(\frac{y_0}{a_0})\in A$. Let us suppose that there exists
$\bbf{x}_{n}=(x_0,\dots,x_n)$ such that $T_{p_n(\bbf{a})}\bbf{x}_n=p_n(\bbf{y})$.
The required $\bbf{x}_{n+1}$ is given by $(x_0,\dots,x_{n+1})$ with $x_{n+1}$
such that
$$
V^{n+1}[a_0x_{n+1}]=p_{n+1}(\bbf{y})-
\sum_{i=1}^{n+1}V^i([a_i](x_0,\dots,x_n,0))-[a_0](x_0,\dots,x_n,0).
$$
The existence of $x_{n+1}$ is ensured by the fact that $a_0$ is invertible and by
the fact that the projection of the right hand side on $W_n$ is zero by induction.
\end{proo}

\begin{lemm}\label{lem:cancellation property Ta}
For any $\bbf{x}=(x_0,x_1,x_2,\dots)\in W(A)$ with $x_0$ not a
zero divisor, the map $T_{\bullet}\bbf{x}$ is injective.
If $x_0$ is invertible then it is bijective.
\end{lemm}

\begin{proo}
Let $\bbf{a}=(a_0,a_1,a_2,\dots)$ and
$\bbf{b}=(b_0,b_1,b_2,\dots)$ as above. We will prove by
induction that $a_n=b_n$ for any $n$. If
$T_{\bbf{a}}\bbf{x}=T_{\bbf{b}}\bbf{x}$ in particular
$a_0x_0=b_0x_0$. Since $x_0$ is not a zero divisor them
$a_0=b_0$. Now let us suppose that $a_i=b_i$ for $i\le n$. We
prove $a_n=b_n$. By hypothesis, we have
$$
T_{\bbf{a}}\bbf{x}-T_{\bbf{b}}\bbf{x}=T_{\bbf{a}}\bbf{x}
=\sum_{k=0}^{\infty}V^k(([a_k]-[b_k])\bbf{x})
=\sum_{k={n+1}}^{\infty}V^k(([a_k]-[b_k])\bbf{x})=0
$$
In particular we have $a_{n+1} x_0=b_{n+1}x_0$ which implies
$a_{n+1}=b_{n+1}$ since $x_0$ is not a zero divisor.
To prove the surjectivity when $x_0$ is invertible one proceeds in a similar
way as in the previous lemma and it is even simpler.
\end{proo}

We now introduce a new, nonassociative product between matrices with
Witt vector entries.

\begin{defi} \label{df:star_T}
Let $M=(\bbf{m}_{i}^j)$ and $N$ be two matrices belonging to $M_n(W(A))$.
We define the $T$-multiplication by
$$
M\star_T N:=  T_M(N)
$$
where $T_M$ is the matrix of operators $(T_{\bbf{m}_i^j})_{1\le i,j\le n}$.
\end{defi}

Endowed with this composition law, $M_n(W(A))$ is a magma and the identity matrix
is a two-sided unit element. We will now consider the set
$\sH_n(W(A))\subseteq M_n(W(A))$ of upper triangular matrices of the form
$$
\left(
\renewcommand{\arraystretch}{1.5}
\begin{array}{ccccc}
\bbf{a}_1^1 & \bbf{a}_{1}^{2}  & \bbf{a}_1^3 & \dots & \bbf{a}_1^n \\
    0   & \bbf{a}_2^2 & \bbf{a}_2^3 & \dots & \bbf{a}_{2}^{n} \\
\vdots & & & & \vdots \\
& & 0& \bbf{a}_{n-1}^{n-1} & \bbf{a}_{n-1}^n \\
& & & 0& \bbf{a}_n^n \\
\end{array}
\right)
$$
with  $\bbf{a}_i^j=(a_{i0}^j,a_{i1}^j,a_{i2}^j,\dots)$ and
$a_{i0}^i$ not a zero divisor. We refer to \ref{U_and_L} for the definition
of the operators $\cU$ and $\cL$, taking a square matrix to its
upper left and lower right codimension $1$ submatrices.

\begin{lemm}\label{lem:cancellation property starT}
The set  $\sH_n(W(A))$ is a submagma of $M_n(W(A))$ and the cancellation laws
hold, i.e. if $M\star_T N=M'\star_T N$ then $M=M'$ and if $M\star_T N=M\star_T N'$
then $N=N'$. Moreover if $A$ is a field then $\sH_n(W(A))$ is a loop.
\end{lemm}

\begin{proo}
It is easy to prove that $\sH_n(W(A))$ is stable under $\star_T$.
We now prove that the cancellation laws hold by induction on $n$.
For $n=1$ this is just lemmas~\ref{lem:Ta injective} and
\ref{lem:cancellation property Ta}. Let us suppose that the cancellation
laws hold in $\sH_n(W(A))$ and prove them for $\sH_{n+1}(W(A))$.
We observe that for any $M,N\in \sH_{n+1}(W(A))$ we have
$$
\cU(M\star_T N)= \cU(M)\star_T \cU(N)
$$
and
$$
\cL(M\star_T N)= \cL(M)\star_T \cL(N).
$$
Therefore if  $M\star_T N=M'\star_T N$ we have, by induction, that
$\cU(M)=\cU(M')$ and $\cL(M)=\cL(M')$. Similarly if $M\star_T N=M\star_T N'$
then $\cU(N)=\cU(N')$ and $\cL(N)=\cL(N')$. It remains to
prove that $\bbf{m}_1^{n+1}={\bbf{m}'}_1^{n+1}$ and
$\bbf{n}_1^{n+1}={\bbf{n}'}_1^{n+1}$.  We begin with the first.
From
$$
(M\star_T N)_1^{n+1}=(M'\star_T N)_1^{n+1}
$$
it follows that
$$
\sum_{j=1}^{n+1}T_{\bbf{m}_1^j}\bbf{n}_j^{n+1}=\sum_{j=1}^{n+1}T_{{\bbf{m}'}_1^j}\bbf{n}_j^{n+1}.
$$
Since $\bbf{m}_1^j={\bbf{m}'}_1^j$ for $j=1,\dots,n$ it follows that
$$
T_{\bbf{m}_1^{n+1}}\bbf{n}_{n+1}^{n+1}=T_{{\bbf{m}'}_1^{n+1}}\bbf{n}_{n+1}^{n+1}
$$
which implies $\bbf{m}_1^{n+1}={{\bbf{m}'}_1}^ {n+1}$ by Lemma
\ref{lem:Ta injective}. Now from
$$
(M\star_T N)_1^{n+1}=(M\star_T N')_1^{n+1}
$$
it follows that
$$
\sum_{j=1}^{n+1}T_{\bbf{m}_1^j}\bbf{n}_j^{n+1}=\sum_{j=1}^{n+1}T_{\bbf{m}_1^j}{\bbf{n}'}_j^{n+1}.
$$
Since $\bbf{n}_1^j={\bbf{n}'}_1^j$ for $j=1,\dots,n$ it follows that
$$
T_{\bbf{m}_1^1}\bbf{n}_1^{n+1}=T_{\bbf{m}_1^1}{\bbf{n}'}_1^{n+1}
$$
which implies $\bbf{n}_1^{n+1}={\bbf{n}'}_1^ {n+1}$ by Lemma
\ref{lem:cancellation property Ta}.

To prove the fact that $\sH_n(W_n(A))$ is a loop if $A$ is a field one proceeds
similarly, using the second part of Lemmas \ref{lem:Ta injective} and
\ref{lem:cancellation property Ta}.
\end{proo}
\end{noth}





\subsection{Deformed Artin-Hasse exponentials} \label{ss:DAHe}

In this section we introduce some deformations of Artin-Hasse exponentials
which we will need in the following.

\begin{defi}
Given indeterminates $\Lambda$, $U$ and $T$, we define a formal power series
in $T$ with coefficients in $\QQ[\Lambda,U]$ by
$$
E_p(U,\Lambda,T)=
(1+\Lambda T)^{\frac{U}{\Lambda}}\prod_{k=1}^\infty\,
(1+\Lambda^{p^k}T^{p^k})^{\frac{1}{p^k}\left(\left(\frac{U}{\Lambda}\right)^{p^k}
-\left(\frac{U}{\Lambda}\right)^{p^{k-1}}\right)} \ .
$$
\end{defi}

It satisfies basic properties such as $E_p(0,\Lambda,T)=1$ and
$E_p(MU,M\Lambda,T)=E_p(U,\Lambda,MT)$, where $M$ is another indeterminate.
It is a deformation of the classical Artin-Hasse exponential
$E_p(T)=\prod_{k=0}^\infty\,\exp(T^{p^k}/p^k)$ in the sense that $E_p(1,0,T)=E_p(T)$.
To see this it is sufficient to observe that, for any $k$, the series
$(1+\Lambda^{p^k}T^{p^k})^{\frac{1}{p^k}\big(\frac{1}{\Lambda^{p^k}}-\frac{1}{\Lambda^{p^{k-1}}}\big)}$
is equal to
  $\bigg({(1+\Lambda^{p^k}T^{p^k})^{\frac{1}{\Lambda^{p^k}}}}\bigg)^{\frac{1-\Lambda^p}{p^k}}$,
  and this gives
 $\exp({T^{p^k}}/{p^k})$
  for $\Lambda=0$.
\begin{defi} \label{df:DAHE}
Given a vector of indeterminates $\UU=(U_0,U_1,\dots)$, we define
a power series in~$T$ with coefficients in $\QQ[\Lambda,U_0,U_1,\dots]$ by
\begin{equation*}
E_p(\UU,\Lambda,T)=\prod_{\ell=0}^\infty\,E_p(U_\ell,\Lambda^{p^\ell},T^{p^\ell}).
\end{equation*}
\end{defi}

We have the following fundamental lemma.

\begin{lemm}
The series $E_p(U,\Lambda,T)$
and $E_p(\UU,\Lambda,T)$ are integral at $p$, that is, they have their coefficients
in $\ZZ_{(p)}[\Lambda,U]$ and $\ZZ_{(p)}[\Lambda,U_0,U_1,\dots]$ respectively.
\end{lemm}
\begin{proo}
See \cite{SS2}, Corollary 2.5.
\end{proo}

It follows from this lemma that given a $\ZZ_{(p)}$-algebra $A$, elements
$\lambda,a\in A$ and $\bbf{a}=(a_0,a_1,\dots)\in A^\NN$, we have specializations
$E_p(a,\lambda,T)$ and $E_p(\bbf{a},\lambda,T)$ which are power series in
$T$ with coefficients in~$A$. We usually consider $\bbf{a}$ as a Witt vector,
i.e. as an element in $W(A)$.


\begin{rema} \label{rm:W_and_Lambda}
Let $\AA^1=\Spec(\ZZ_{(p)}[\Lambda])$ be the affine line over the ring of
$p$-integers $\ZZ_{(p)}$, with coordinate $\Lambda$, and write $W_{\AA^1}$
for the scheme of Witt vectors over $\AA^1$.
We remark (see \cite{SS2}, Corollary 2.9.1) that, generalizing what happens
for the Artin-Hasse exponential, the deformed exponential of
Definition~\ref{df:DAHE} gives a homomorphism
$$
W_{\AA^1}\too \mathbf{\Lambda}_{\AA^1}
$$
where $\mathbf{\Lambda}_{\AA^1}=\Spec(\ZZ_{(p)}[\Lambda,X_1,\dots,X_n,\dots])$
is the $\AA^1$-group scheme  whose group of
$R$-points, for any $\ZZ_{(p)}[\Lambda]$-algebra $R$, is the abelian
multiplicative group $1+TR[[T]]$.
(We hope that the difference between the symbols $\Lambda$ and
$\mathbf{\Lambda}$ is visible enough.)
The above homomorphism is in fact a closed immersion. We also note that
there is an isomorphism:
$$
\prod_{p\nmid k}W_{\AA^1} \simeq \mathbf{\Lambda}_{\AA^1}
$$
which works as follows. With any
$\ZZ_{(p)}$-algebra $R$, any element $\lambda\in R$, and any family of Witt
vectors $\bbf{a}_k=(a_{k0},a_{k1},a_{k2},\dots)\in W(A)$ indexed by the
prime-to-$p$ integers~$k$, this isomorphism associates the series
$F(T)=\prod_{p\nmid k}\,E_p(\bbf{a}_k,\lambda,T^k)$
(see \cite{MRT}, Lemma 3.1.2). \hfill $\square$
\end{rema}


Here are a couple more definitions which will be useful in the sequel.
We set
$$
\tilde{p}E_p(\UU,\Lambda,T)=E_p(V(U_0^p,U_1^p,\ldots),\Lambda,T).
$$
where $V$ is the Verschiebung. Using the isomorphism
$\prod_{k\nmid p}W_{\AA^1}\simeq\mathbf{\Lambda}_{\AA^1}$
described above, one extends this definition to any element of
$1+T\ZZ_{(p)}[U_1,\dots,U_n,\Lambda][[T]]$. The result is a
group scheme endomorphism
$$
\tilde p:\mathbf{\Lambda}_{\AA^1}\too\mathbf{\Lambda}_{\AA^1}.
$$
In \cite{SS1} this operator is called $[p]$, but we prefer $\tilde p$ to
avoid confusion with Teichm\"uller representatives. Also, we define an
additive endomorphism $F^\Lambda:=F-[\Lambda^{p-1}]:W_{\AA^1}\to W_{\AA^1}$.
For each element $\lambda$ in a $\ZZ_{(p)}$-algebra $R$, this gives
an endomorphism $F^\lambda:W_R\to W_R$. When $R$ is a discrete
valuation ring with uniformizer $\pi$ and $\lambda=\pi^l$ for some
$l>0$, we will sometimes write $F^{(l)}$ instead of $F^{\pi^l}$
(see e.g. the statement of Theorem~\ref{thm 2 SS}).

\begin{defi}\label{eq: E tilde}
Let $\Lambda_2$ be another indeterminate.
For any $H$ in  $1+T\ZZ_{(p)}[U_1,\dots,U_n,\Lambda][[T]]$ we define
the series
\begin{equation}
\tilde{E}_p(\WW,\Lambda_2, H)=H^{\frac{W_0}{\Lambda_2}}
\prod_{r=1}^{\infty}\big(\tilde{p}^rH\big)
^{\frac{1}{p^r\Lambda_2^{p^r}}\Phi_{r-1}(F^{\Lambda_2}(\WW))}.
\end{equation}
\end{defi}

From the definition, one sees that $\tilde{E}_p(\WW,\Lambda, H)$
gives a bilinear group scheme homomorphism
$$
W_{\AA^1}\times \mathbf{\Lambda}_{\AA^1}\to \mathbf{\Lambda}_{\AA^1}.
$$
With some quite simple computations  one shows  the following lemma.

\begin{lemm}\label{eq: E tilde avec E et operation T}
In the group
$1+T\ZZ_{(p)}[\WW,\frac{\UU}{\Lambda_2},\Lambda,\Lambda_2][[T]]$, we have
$$
\tilde{E}_p(\WW,\Lambda_2, E_{p}(\UU,\Lambda;T))
=E_p(T_{\UU/\Lambda_2}\WW,\Lambda;T).
$$
\end{lemm}

\begin{proo}
See \cite{SS1}, Proposition 4.11.
\end{proo}

In particular, we have
$\tilde{E}_p(\WW,\Lambda, 1+\Lambda T)=E_p(\WW,\Lambda;T)$.
Finally we define the following series.

\begin{defi}\label{G_p}
For any $H$ as above, we define
\begin{equation*}
G_p(\WW,\Lambda_2,H)
=\prod_{r=1}^{\infty}\bigg(\frac{1+(H-1)^{p^r}}{\tilde{p}^rH}\bigg)
^{\frac{1}{p^r\Lambda_2^{p^r}}
{\Phi_{r-1}(\WW)}}\in 1+T\QQ[\WW,\UU,\Lambda,\Lambda_2,\frac{1}{\Lambda_2}][[T]].
\end{equation*}
\end{defi}

Using \cite{SS2}, Lemma 2.8, one sees immediately that
\begin{equation}\label{eq: G_p vs E_p}
G_p(F^{\Lambda_2}(\WW),\Lambda_2, H)
=\frac{E_p(\WW,\Lambda_2;\frac{H-1}{\Lambda_2})}{\tilde{E}_p(\WW,\Lambda_2, H)}.
\end{equation}
We remark that for any $H$ as above we have
\begin{equation}\label{eq:G_p additive}
G_p(\WW,\Lambda_2, H)\,G_p(\WW',\Lambda_2, H)
=G_p(\WW+\WW',\Lambda_2, H)
\in 1+T\QQ[\WW, \WW', \UU,\Lambda,\Lambda_2\frac{1}{\Lambda_2}][[T]]
\end{equation}
where $\WW+\WW'$ is the sum of Witt vectors. We finally have the following lemma.

\begin{lemm}
We have
$G_p(\WW,\Lambda_2,E_p(\UU,\Lambda_2;T))
\in\ZZ_{(p)}[\WW,\frac{\UU}{\Lambda_2},\Lambda,\Lambda_2][[T]]$.
\end{lemm}


\begin{proo}
See \cite{SS1}, Proposition 4.12.
\end{proo}

It is quite simple to verify the following equality.

\begin{lemm}\label{eq: E tilde avec G et operation T}
We have
$\tilde{E}_p(\WW,\Lambda_3, G_{p}(\UU,\Lambda_2;H))
=G_p(T_{\UU/\Lambda_3}\WW,\Lambda_2;H)$.
\end{lemm}

\begin{proo}
See \cite{SS1}, Proposition 4.13.
\end{proo}

\subsection{Main theorems of Sekiguchi-Suwa Theory} \label{ss:main_thms}

In this section, we briefly recall the main results of Sekiguchi-Suwa Theory,
stated in \cite{SS1}. One can also find a summary of this
theory in wider generality in \cite{MRT}. From now on, we denote by $R$ a
discrete valuation ring of unequal characteristics. We stress that, in contrast
with Sections \S~\ref{secBK} to \S~\ref{BKn3} we do not assume that $R$ is complete
and neither that its residue field is perfect. We will denote by $\pi$ a fixed
uniformizer of $R$ and by $v$ the valuation of $R$.

\begin{defi} \label{df:filtered_gp_schemes}
Let $l,l_1,\dots,l_n$ be integers.
\begin{trivlist}
\itemn{1}
We let $\cG^{(l)}$ be the group scheme $\Spec(R[T,1/(\pi^l T+1)])$ with
group law $T*T'=T+T'+\pi^l TT'$, the unique group law such
that the morphism
$\alpha:\Spec(R[T,1/(\pi^l T+1)])\to \GG_m=\Spec(R[T,1/T])$
given by $T\mapsto 1+\pi^l T$ is  a group scheme homomorphism.
\itemn{2} Let $\cE$ be a flat $R$-group
scheme. If there exist exact sequences of flat $R$-group schemes
$$
0\too \cG^{(l_i)}\too \cE_i\too \cE_{i-1}\too 0
$$
for $1\le i \le n$, with $\cE_0=0$ and $\cE_n=\cE$, we call the
sequence of flat $R$-group schemes
$$
\cE_1=\cG^{(l_1)},\cE_2,\dots,\cE_n=\cE
$$
or, sometimes, simply $\cE$,  a filtered $R$-group scheme of type
$(l_1,\dots,l_n)$.
\end{trivlist}
\end{defi}

\begin{rema}
One can define a group scheme $\cG^{(\lambda)}$ for each $\lb \in R$, in such
a way that $\cG^{(l)}:=\cG^{(\pi^l)}$ is just the group scheme defined
in~\ref{df:filtered_gp_schemes}. In this article, we care only about the
isomorphism class of $\cG^{(\lambda)}$ which depends only on $\lambda$
up to units, so we prefer to adopt the more compact notation.
\end{rema}


\begin{theo}\label{thm 1 SS}
Let $\cE=(\cE_1,\dots,\cE_n)$ be a filtered group scheme of
type $({l_1},\dots,{l_n})$, with $l_i>0$ for each $i$. Then there are
compatible  open immersions of $\cE_i\to \AA^{i}$ and elements
$$
{D}_i\in H^0(\AA^i_R,\cO_{\AA^i_R})=R[T_1,\dots, T_i]
$$
such that,
for each $1\le i\le n$, the Hopf algebra of $\cE_i$ is given by
$$
R[\cE_i]= R\Big[T_1,\dots,T_i,\frac{1}{1+\pi^{l_1}
T_1},\frac{1}{{D}_1(T_1)+\pi^{l_2}
T_2},\dots,\frac{1}{{D}_{i-1}(T_1,\dots,T_{i-1})+\pi^{l_i}
T_i}\Big]
$$
The group law of $\cE_i$ is the one which makes  the morphism
\begin{align*}
\alpha_{\cE_i}: & \ \cE_i\too (\GG_{m,R})^i\\
& (T_1,\dots,T_i) \mto (1+\pi^{l_1} T_1,{D}_1(T_1)+\pi^{l_2}
T_2,\dots, {D}_{i-1}(T_1,\dots,T_{i-1})+\pi^{l_i} T_i)
\end{align*}
a group-scheme homomorphism and the reduction modulo
$\pi^{l_{i+1}}$ of the function ${D}_i:\AA^i_R\to\AA^1_R$ factors
into a group scheme homomorphism
${D_i}_{|\cE_i}:\cE_{i,R/\pi^{l_{i+1}}R}\to
\GG_{m,R/{\pi^{l_{i+1}}}R}\subseteq \AA^1_{R/{\pi^{l_{i+1}}}R}$.

Moreover if $l_{n+1}$ is a positive integer and
${D}_n:\AA^n_R\to\AA^1_R$ is a function whose reduction modulo
$\pi^{l_{n+1}}$ factors into a group scheme homomorphism
$$
{D_{n}}_{|\cE_n}:\cE_{n,R/\pi^{l_{n+1}}R}\to \GG_{m,R/{\pi^{l_{n+1}}}R}
\subseteq \AA^1_{R/{\pi^{l_{n+1}}}R}$$
then
$$
R[\cE_{n+1}]:=
R[\cE_n]\left[T_{n+1},1/({D}_n(T_1,\dots,T_n)+\pi^{l_{n+1}}T_{n+1})\right]
$$
is the Hopf algebra of a filtered group scheme $\cE_{n+1}$ of type
$(l_1,\dots,l_{n+1})$, where the group scheme structure is the
only one which turns into a group scheme homomorphism the
morphism $\alpha_{\cE_{n+1}}: \cE_{n+1}\to (\GG_{m,R})^{n+1}$
which extends $\alpha_{\cE_{n}}$ and sends $T_{n+1}$ to
${D}_n(T_1,\dots,T_n)+\pi^{l_{n+1}}T_{n+1}$.

Finally, a polynomial ${D}'_{n}\in R[T_1,\ldots,T_n]$ with the
same reduction modulo $\pi^{l_{n+1}}$ as ${D}_{n}$ gives the same
filtered group scheme up to isomorphism.
\end{theo}

\begin{proo}
See \cite{SS1}, Theorem 3.2 and Theorem 3.3.
\end{proo}

%

In fact one can describe very explicitly the polynomials which
appear in the above theorem. In the next statement and in the rest
of the article, we sometimes write $f:X\endto$ for a map $f:X\to
X$ from some set to itself.

\begin{theo}\label{thm 2 SS}
Let $\cE$ be a filtered group scheme of type $(l_1,\dots, l_n)$
with $l_i>0$ for each $i$. Then there exist elements
${\bbf{a}}_{i}^j\in W^f(R)$ with $1\le i<j\le n$, whose reductions
modulo $\pi^{l_{j}}$ are in $\widehat{W}(R/\pi^{l_{j}}R)$, such
that
\begin{itemize}
\item one can take,  for any $j=1,\dots,n-1$,
${D}_j(T_1,\dots,T_j)$ as the truncation of
$$E_p(({\bbf{a}}^{j+1}_i)_{1\le i\le j},
(\pi^{l_k})_{1\le k\le j}; T_1,\dots, T_j)$$ in degree $r$, where
$ E_p\big(({\bbf{a}}^{j+1}_i)_{1\le i\le j}, (\pi^{l_k})_{1\le
k\le j}; T_1,\dots, T_j\big)$ is the series defined by induction
$$
\prod_{i=1}^j E_p\left({\bbf{a}}^{j+1}_i, \pi^{l_i};
\frac{T_i}{E_p(({\bbf{a}}^{i}_s)_{1\le s\le i-1},
(\pi^{l_k})_{1\le k\le i-1}; T_1,\dots, T_{i-1})}\right)
$$
and  $r$ is the degree of the reduction of this series
modulo $\pi^{l_{j+1}}$, which is  a polynomial;
\item the reduction modulo $\pi^{l_j}$ of each $({\bbf{a}}^{j}_i)_{1\le i\le j-1}$
is in the kernel of the operator
$$
U^{j-1}:\widehat{W}(R/\pi^{l_{j}}R)^{j-1}\endto
$$
defined as follows:
$U^1$ is defined as $F^{(l_1)}:=F-[\pi^{(p-1)l_1}]$ and
we define
$$
U^{n}=\left(
\renewcommand{\arraystretch}{1.5}
\begin{array}{ccccc}
  F^{(l_1)} & -T_{\bbf{b}_1^2}& -T_{\bbf{b}_1^3} &\dots & -T_{\bbf{b}_1^n}\\
  0 & F^{(l_2)} & -T_{\bbf{b}_2^3} &\dots& -T_{\bbf{b}_2^n} \\
  \vdots & 0 & \ddots & \ddots &\vdots \\
  0 & 0 & \ddots & \ddots & -T_{\bbf{b}_{n-1}^n} \\
  0 & 0 & \dots & 0 & F^{(l_n)} \\
\end{array}%
\right)
$$
where $\cU(U^n)=U^{n-1}$ and $\cL(U^n)$ are defined by induction and
$$
\bbf{b}_1^{n}:=\frac{1}{\pi^{l_n}} \left(F^{(l_1)}{\bbf{a}}_1^{n}
-\sum_{t=2}^{n-1}T_{\bbf{b}_1^t}{\bbf{a}}_t^n\right)
=\frac{U^{n-1}({\bbf{a}}^{n}_i)_{1\le i\le n-1}}{\pi^{l_n}};
$$
\item  for any $l\in \NN$, we have an isomorphism
$$
\ker\left(U^n:\widehat{W}(R/\pi^lR)^n\endto\right)
\too\Hom_{R/\pi^{l}R\mbox{-}\Gr}({i^*\cE,{\GG}_{m,R/\pi^{l}R}}),
$$
given by
$$
\bbf{c}^n\mto E_p(\bbf{c}^{n},(\pi^{l_j})_{1\le j\le
n},T_1,\dots,T_n),
$$
where $i$ is the closed immersion $\Spec(R/\pi^l R)\to \Spec(R)$.
\end{itemize}
\end{theo}

\begin{proo}
See \cite{SS1}, Theorem 5.1 and Theorem 5.2.
\end{proo}

%
%

\subsection{Sekiguchi-Suwa Theory from a matricial point of view}
\label{sub:SS with matrices}

Our purpose here is to introduce "simple" matrices parametrizing filtered
group schemes (\ref{the_set_Mn}) and to translate in matricial terms
the main operations on group schemes: quotients and subgroups
(\ref{U_L_quotients_sbgps}) and model maps (\ref{pos_matrices_model_maps}).
In the following, we always suppose that the parameters $l_i$ of the
filtered group schemes we are considering are positive ($l_i>0$).

Let $\sH_n(W(K))$ be the loop constructed
in~\ref{lem:cancellation property starT}. For matrices $A,B\in\sH_n(W(K))$
we will make use of the notations $A/B$ and $A\backslash B$ as defined
in~\ref{df:quasigroups}.
In a similar way as in~\ref{def:positive elements} we will say that
a matrix $A$ in $\sH_n(W(K))$ is {\em positive}, and we will write
$A\ge 0$, if it belongs to $\sH_n(W(R))$.

\begin{noth} {\bf The set $\sM_n$.} \label{the_set_Mn}
To start with, we need a technical remark allowing to
reformulate the congruences in Theorem~\ref{thm 2 SS}.
We consider an upper triangular matrix of the following form:
$$
A=\left(
\renewcommand{\arraystretch}{1.5}
\begin{array}{ccccc}
[\pi^{l_1}] & {{\bbf{a}}}_{1}^{2}  & {{\bbf{a}}}_1^3
& \dots & {{\bbf{a}}}_1^n \\
& {[\pi^{l_2}]} & {{\bbf{a}}}_2^3
& & \vdots \\
& & \ddots & \ddots & \vdots \\
& & & [\pi^{l_{n-1}}] & {{\bbf{a}}}_{n-1}^n \\
0 & & & & [\pi^{l_n}] \\
\end{array}
\right)\in M_n(W^f(R)) \ , \ \mbox{all } l_i>0.
$$

\begin{lemm} \label{lm:technical_point}
For each matrix $A$ as above, let $F(A)$ be the matrix obtained by applying
Frobenius to all entries. Then the following conditions are equivalent:
\begin{trivlist}
\itemn{1} for each $j\in\{1,\dots,n\}$, the reduction of
${({\bbf{a}}_i^j)}_{1\le i\le j-1}$ belongs to
$\widehat{W}(R/\pi^{l_j}R)^{j-1}$ and
$$
U^{j-1}({{\bbf{a}}_i^j})_{1\le i\le  j-1}\equiv 0 \mod \pi^{l_j}
$$
where $U^{j-1}$ is defined by induction in Theorem \ref{thm 2 SS}.
\itemn{2} $F(A)/A\ge 0$.
\end{trivlist}
\end{lemm}

Note that the operator $U^{j-1}$ in (1) depends only on the
vectors ${\bbf{a}}^{k}\in W(R)^k$ with $1\le k\le j-1$.

\begin{proo}
In fact, we have
$$
F(A)/A=\left(
\renewcommand{\arraystretch}{1.5}
\begin{array}{ccccc}
[\pi^{(p-1)l_1}] & {\bbf{b}}_{1}^{2}  & {\bbf{b}}_1^3 & \dots & {\bbf{b}}_1^n \\
    0   & {[\pi^{(p-1)l_2}]} & {\bbf{b}}_2^3 & \dots & {\bbf{b}}_{2}^{n} \\
\vdots & \ddots & \ddots & & \vdots \\
\vdots & & \ddots & [\pi^{(p-1)l_{n-1}}] & {\bbf{b}}_{n-1}^n \\
0 & \dots & \dots & 0 & [\pi^{(p-1)l_n}] \\
\end{array}
\right)\in \sH_n(W(K))
$$
where the $\bbf{b}_i^j$ are defined as in \ref{thm 2 SS}. By the
definition of $\bbf{b}_i^j$, this matrix is in $\sH_n(W(R))$ if
and only if the congruences in (1) are satisfied. It remains to
prove that if $F(A)/A\ge 0$ then for each $j$ the reduction of
${({\bbf{a}}_i^j)}_{1\le i\le j-1}$ belongs to
$\widehat{W}(R/\pi^{l_j}R)^{j-1}$. We prove this by induction on
$n$. We observe that since the entries of $A$ are in $W^f(A)$,
then this condition simply means that the entries of $A$ are
congruent to $0$ modulo $\pi$, i.e. $A/[\pi]\Id$ is positive. For
$n=1$ there is nothing to prove. Let us suppose the statement true
for $n$ and prove it for $n+1$. Then one has
\begin{equation}\label{eq:F(A)/A}
F(A)/A=
\left(
\renewcommand{\arraystretch}{1.5}
\begin{array}{cccc}
& & & \\
& \cU F(A)/\cU A & &
\frac{F((\bbf{a}_i^{n+1})_{1\le i\le n})
-T_{(\cU F(A)/\cU A)}(\bbf{a}_i^{n+1})_{1\le i\le n}}{\pi^{l_{n+1}}} \\
& & & \\
0 & \dots & 0 & [\pi^{(p-1)l_{n+1}}] \\
\end{array}
\right)
\end{equation}
and
$$
F(A)/A= \left(
\begin{array}{cc}
[\pi^{(p-1)l_{n+1}}] &
\frac{F(({\bbf{a}_1^{j}})_{2\le j\le n+1})
-T_{(\cL F(A)/\cL A)}(\bbf{a}_1^{j})_{2\le j\le n+1}}{\pi^{l_{j}}} \\
0 & \\
\vdots & \cL F(A)/\cL A \\
0 & \\
\end{array}
\right)
$$
where we use the notation $T_M(N)$ from Definition~\ref{df:star_T}.
By the inductive hypothesis, the matrices $\cU A/[\pi]\Id$ and $\cL A/[\pi]\Id$
are positive. Therefore it remains to prove that
 $\bbf{a}_1^{n+1}\equiv 0\mod \pi$.  Since $F(A)/A$ is positive and
$l_{n+1}>0$ then from \eqref{eq:F(A)/A} we derive
$$
F(\bbf{a}_1^{n+1})\equiv
T_{(\cU F(A)/\cU A)}{(\bbf{a}_i^{n+1})_{1\le i\le n}} \mod \pi.
$$
Since by induction $\bbf{a}_i^{n+1}\equiv 0 \mod \pi^{l_{n+1}}$ for
$i=2,\dots,n$ and $\cU^n F(A)/\cU^n A=([\pi^{(p-1)l_1}])$  with $l_1>0$,
then we have
$$
F(\bbf{a}_1^{n+1})\equiv [\pi^{(p-1)l_1}]\bbf{a}_1^{n+1}\equiv 0\mod \pi.
$$
This implies that $\bbf{a}_1^{n+1}\equiv 0\mod \pi$.
\end{proo}

Theorems \ref{thm 1 SS} and \ref{thm 2 SS} imply that to any matrix satisfying
the equivalent conditions of Lemma~\ref{lm:technical_point} one can
attach a unique filtered group scheme $\cE(A)$.
Conversely, for any filtered group scheme $\cE$ one can find a matrix $A$
satisfying these conditions such that $\cE=\cE(A)$. This leads us to introduce
the relevant set of matrices. Note that if $\cE$ is given, then a matrix $A$
such that $\cE=\cE(A)$ is not unique. So we have to identify the
equivalence relation saying that two matrices define the same filtered group;
this will be done in \ref{pos_matrices_model_maps}.

\begin{defi}
Let $n\in\NN$ and $\bbf{l}=(l_1,\dots,l_n)\in (\NN_{>0})^n$. We define
$$
\sM^{\bbf{l}}_n:=\big\{A=(\bbf{a}_{i}^j)\in M_n(W^f(R)),
\text{upper triangular}, \bbf{a}_i^i=[\pi^{l_i}]
\text{ for } 1\le i\le n \text{ and } F(A)/A\ge 0\big\}
$$
and $\sM_n:=\bigcup \sM^{\bbf{l}}_n$, the union being over all
${\bbf{l}\in (\NN_{>0})^n}$.
\end{defi}

\begin{rema}
If $A\in \sM_n$, then it is not necessarily the case that
$F(A)\in \sM_n$. There are  counterexamples already for $n=2$,
with $l_2\gg l_1$.
\end{rema}

By Theorems~\ref{thm 1 SS}, \ref{thm 2 SS} and
Lemma~\ref{lm:technical_point}, one can associate with any $A\in
\sM^{\bbf{l}}_n$ a filtered group scheme $\cE(A)$ of type
$(l_1,\dots,l_n)$. It is constructed by successive extensions
defined by deformed exponentials ${D}_j(T_1,\dots,T_j)$ equal to
the truncation of $E_p(({\bbf{a}}_{i}^{j+1})_{1\le i\le j},
(\pi^{l_k})_{1\le k\le j}; T_1,\dots, T_j)$ in degree $r$, where
$r=r_i$ is the degree of the reduction modulo $\pi^{l_{i+1}}$ of
this series. We call ${D}_j$ the {\em truncated exponential}
associated with $(\bbf{a}_i^{j+1})_{1\le i\le j}$. 
(Note that similar truncated exponentials appear in the work~\cite{GrCh}.)
With the vocabulary introduced in the article \cite{MRT} (see especially
3.2 and 4.3 there), the vectors ${\bbf{a}}^j$ are {\em frames} for
the filtered group $\cE(A)$, and the matrix~$A$ may be called a
{\em matrix of frames}.
\end{noth}

\begin{noth} {\bf Operators $\cU$ and $\cL$ versus quotients and subgroups.}
\label{U_L_quotients_sbgps}
It is clear that for each $A\in\sM_n$ and $i\in\{1,\dots,n\}$
we have $\cU^{n-i}A\in\sM_i$ and $\cL^iA\in\sM_{n-i}$. Here is the precise
meaning of the operators $\cU$ and $\cL$ for filtered group schemes.

\begin{prop}
Let $\cE=(\cE_1,\dots,\cE_n)$ be a filtered group scheme of type
$\bbf{l}=(l_1,\dots,l_n)$, and $A\in\sM^{\bbf{l}}_n$. For
$1\le i\le n-1$ consider the exact sequence
$0\to\cS_{n-i}\to\cE_n\to\cQ_i\to 0$
where $\cQ_i:=\cE_i$ (quotient of dimension $i$) and
$\cS_{n-i}:=\ker(\cE_n\to\cE_i)$ (subgroup of codimension $i$). Then:
\begin{trivlist}
\itemn{1} $\cQ_i$ is a filtered group scheme of type $(l_1,\dots,l_i)$.
If $\cE=\cE(A)$ then $\cQ_i=\cE(\cU^{n-i}A)$.
\itemn{2} $\cS_{n-i}$ is a filtered group scheme of type $(l_{i+1},\dots,l_n)$.
If $\cE=\cE(A)$ then $\cS_{n-i}=\cE(\cL^iA)$.
\end{trivlist}
\end{prop}

\begin{proo}
Assertion (1) comes from the inductive construction of $\cE_n$. For the
proof of (2), we set $\cK_d=\ker(\cE_d\to\cE_i)$ for each $d\ge i+1$.
First, we show by induction on $d$ that $\cK_d$ is a filtered group scheme
of type $(l_{i+1},\dots,l_d)$. The initialization at $d=i+1$ is clear and
the inductive step is verified since the morphism $\nu_d:\cE_{d+1}\to\cE_d$
with kernel $\cG^{(l_{d+1})}$ induces an exact sequence:
$$
0\too\cG^{(l_{d+1})}\too\cK_{d+1}\stackrel{\nu_d}{\too}\cK_d\too 0.
$$
In order to prove that $\cS_{n-i}=\cE(\cL^iA)$ if $\cE=\cE(A)$,
we examine more closely the way these extensions are built. The extension
$\cE_{d+1}$ is constructed from $\cE_d$ using a morphism
$D_d:\cE_d\to i_*\GG_m$ where $i:\Spec(R/\pi^{l_{d+1}}R)\into\Spec(R)$ is
the closed immersion. This morphism is the deformed exponential defined
by the coefficients $\bbf{a}^{d+1}_i$ in the $(d+1)$-th column of $A$.
The extension $\cK_{d+1}$ is constructed from $\cK_d$ using the morphism
$D_d|_{\cK_d}:\cK_d\to i_*\GG_m$. In the coordinates $T_1,\dots,T_d$
of~\ref{thm 2 SS}, the closed subgroup scheme $\cK_d\subset\cE_d$ is
defined by the vanishing of the coordinates $T_1,\dots,T_i$. It follows
that $D_d|_{\cK_d}$ is obtained from the deformed Artin-Hasse exponential
$D_d$ by setting $\bbf{a}^{d+1}_1=\bbf{a}^{d+1}_2=\dots=\bbf{a}^{d+1}_i=0$.
Hence the matrix of coefficients that defines $\cK_d$ is the boxed
middle matrix:
\begin{center}
{\footnotesize $A=\left(
\begin{array}{ccc}
\begin{array}{cc} \ddots & \\ & [\pi^{l_i}] \\ \end{array} & * & * \\
0 &
\begin{array}{|ccc|}
\hline
[\pi^{l_{i+1}}] & & \\
& \ddots & \\
& & [\pi^{l_d}] \\
\hline
\end{array}
& * \\
0 & 0 & \begin{array}{cc} [\pi^{l_{d+1}}] & \\ & \ddots \\ \end{array} \\
\end{array}
\right)$}
\end{center}
In symbols, $\cK_d=\cE(\cU^{n-d}\cL^iA)$. For $d=n$, we get
$\cS_{n-i}=\cK_n=\cE(\cL^iA)$.
\end{proo}

\end{noth}


%



\begin{noth} {\bf Positive matrices versus model maps.}
\label{pos_matrices_model_maps}
We use the word {\em unitriangular} as a synonym for
{\em upper triangular unipotent}.

\begin{prop} \label{pp:pos_mat_and_model_maps}
Let $\cE=\cE(A)$ and $\cE'=\cE(A')$ be two filtered group schemes,
with $A,A'\in \sM_n$.
\begin{itemize}
\item There exists a (unique) model map $\cE\to \cE'$ which commutes
with $\alpha_{\cE}$ and $\alpha_{\cE'}$ if and only if $A/A'\ge 0$.
In particular, the relation $\succ$ in $\sM_n$ given by $A\succ A'$
if and only if $A/A'\ge 0$, is transitive. Moreover $\cE$ and $\cE'$ are
isomorphic if and only if $A/A'$ is positive  and  unitriangular.
\item If $\varphi:\cE\to \cE'$ is  such a model map,
the morphism of groups
$$
\varphi^*:\Hom_{R/\pi^{l}R\mbox{-}\Gr}({i^*\cE',{\GG}_{m,R/\pi^{l}R}})
\too \Hom_{R/\pi^{l}R\mbox{-}\Gr}({i^*\cE,{\GG}_{m,R/\pi^{l}R}})
$$
is given, using the isomorphism of Theorem \ref{thm 2 SS},
by the operator $T_{A/A'}$.
\end{itemize}
\end{prop}

\begin{proo}
We prove by induction on the dimension $n$  the following more precise statements:
\begin{itemize}
\item There exists a (unique) model map $\cE\to \cE'$ which
commutes with $\alpha_{\cE}$ and $\alpha_{\cE'}$ if and only if
$A/A'\ge 0$.
\item Let $D_1,\dots, D_{n-1}$ (resp. $D_1',\dots, D_{n-1}'$) be the
truncated polynomials determined by $A$ (resp. $A'$). If a model map
$\varphi^{n}:\cE\to \cE'$ exists then it is given by
$\varphi^{n}=(\varphi_{i})_{1\le i\le n}$, where
$\varphi_1(T_1)=\pi^{l_1'-l_1}$ and
$$
\varphi_{i}(T_1,\dots,T_{i})
=\frac{D_{i-1}(T_{1},\dots,T_{i-1})-{D'_{i-1}}
(\varphi_1(T_1),\dots,\varphi_{i-1}(T_1,\dots,T_{i-1}))}
{\pi^{{l'_{i}}}}+\pi^{l_{i}-{l'_{i}}}T_{i}
$$
for all $i>1$.
\item Write $A=(\bbf{a}_i^j)$ and
$A'=({\bbf{a}'}_i^j)$ in $\sM_{n}$ with $\bbf{a}_i^i=[\pi^{l_i}]$
and ${\bbf{a}'_i}^i=[\pi^{l'_i}]$ for $1\le i \le n$. Let
$(\bbf{w}_i^n)_{1\le i\le n}\in W(R)^n$ be a vector whose
reduction modulo $\pi^l$ belongs to
$$
\ker(U^n:\hat{W}(R/\pi^lR)^n\endto).
$$
Let $H_1=1+\pi^{l_1}T_1$ and
$$
H_n=\frac{E_p\big((\bbf{a}_i^n)_{1\le i \le n-1},
(\pi^{l_i})_{1\le i \le n-1},\TT\big)
\big(1+\pi^{l_{n}}\frac{T_{n}}{E_p((\bbf{a}_i^n)_{1\le i \le n-1},
(\pi^{l_i})_{1\le i \le n-1},\TT)}\big)}{E_p((
  {\bbf{a}'}_i^n)_{1\le i \le n-1},
(\pi^{l_1})_{1\le i \le n},\varphi^{n}(\TT))}
$$
for $n>1$. Then
\begin{equation}\label{eq:E_p(phi)}
\begin{split}
E_p(\bbf{w}^{n},(\pi^{l_1},\dots,\pi^{l_{n}}),\varphi^{n}(\TT))=
& \ E_p(T_{A/ A'}(\bbf{w}^{n}),(\pi^{l_1},\dots,\pi^{l_n}),\TT) \\
& \times \prod_{2\le r\le
n}G_p(U^n_r(\bbf{w}^n),\pi^{l_r};H_{r-1})
\end{split}
\end{equation}
where $U_r^n$ is the $r$-th row of $U^n$, and
\begin{equation}\label{eq:H_n}
\begin{split}
H_n=&\frac{E_p\left((\bbf{a}_i^n)_{1\le i \le n-1}-T_{\cU A/\cU A'}
({\bbf{a}'}_i^n)_{1\le i \le n-1},(\pi^{l_1},\dots,\pi^{l_{n-1}}); \TT\right) \,}
{\prod_{2\le r\le n} G_p(U^n_r(\bbf{w}^n),\pi^{l'_r};H_{r-1})}\\
&\times E_p\left([\pi^{l_{n}}], \pi^{l_{n}};
\frac{\TT}{E_p(({\bbf{a}'}_i^n)_{1\le i \le n-1},
(\pi^{l_i})_{1\le i \le n-1},\varphi^{n}(\TT))}\right).
\end{split}
\end{equation}
\end{itemize}
Note that~(\ref{eq:E_p(phi)}) implies that $\varphi^*$ is given
by the operator $T_{A/A'}$, as asserted in the statement of the proposition.

If $n=1$ we have $\cE=\cG^{(l_1)}$ and $\cE'=\cG^{(l_1')}$ for some
positive integers $l_1,l_1'$. In this case $A/A'\ge 0$ simply means
$l_1\ge l_1'$ and the above statement is known: see
\cite{SOS}, Proposition 1.4 for the first part and second part and
\cite{SS2}, Remark 3.8 for the third part.

We now suppose that the three statements hold true for some $n\ge 1$
and we prove them for $n+1$.
 We have
$$
A/A'=
\left(
\begin{array}{cccc}
& & & \\
& \cU A/\cU A' & &
\frac{(\bbf{a}_i^{n+1})_{1\le i\le n}
-T_{\cU A/\cU A'}(({\bbf{a}'}_i^{n+1})_{1\le i\le n})}{\pi^{{l'_{n+1}}}} \\
& & & \\
0 & \dots & 0 & [\pi^{l_{n+1}-l'_{n+1}}]\\
\end{array}
\right).
$$
Let $D_1,\dots, D_{n+1}$ (resp. $D_1',\dots, D_n'$) be the truncated
exponentials that define $\cE$ (resp. $\cE'$).
The model map $\varphi_{n+1}$ which we are looking for should commute
with $\alpha_{\cE}$ and $\alpha_{\cE'}$, so if it exists it is unique.
So one sees immediately that it exists if and only if
there exists a model  map  $\cE_n\to \cE'_n$ and, if we write
$\varphi^{n+1}=(\varphi_{i})_{1\le i\le n+1}$,  the polynomial
$$
\varphi_{n+1}(T_1,\dots,T_{n+1})=\frac{D_{n}(T_{1},\dots,T_{n})-{D'_n}(\varphi_1(T_1),\dots,\varphi_{n}(T_1,\dots,T_{n}))}{\pi^{{l'_{n+1}}}}+\pi^{l_{n+1}-{l'_{n+1}}}T_{n+1}
$$
belongs to $R[T_1,\dots,T_n]$.
Therefore, by the inductive hypothesis, there exists a model
map between $\cE$ and $\cE'$ if and only if $\cU A/ \cU A'$ is positive,
$l_{n+1}\ge {l'_{n+1}}$ and
\begin{equation}\label{eq: equation for model map}
D_{n}(T_{1},\dots,T_{n})\equiv {D'_n}(\varphi_1(T_1),
\dots,\varphi_{n}(T_1,\dots,T_{n}))\mod \pi^{l'_{n+1}}.
\end{equation}
By induction we have that
$$
\varphi^*:\Hom_{R/\pi^{l}R\mbox{-}\Gr}({i^*\cE(\cU(A')),{\GG}_{m,R/\pi^{l}R}})
\too \Hom_{R/\pi^{l}R\mbox{-}\Gr}({i^*\cE(\cU(A)),{\GG}_{m,R/\pi^{l}R}})
$$
is given by the operator $T_{\cU A/\cU A'}$. This means that the
equation \eqref{eq: equation for model map} is equivalent to
$$
(\bbf{a}_i^{n+1})_{1\le i\le n}\equiv
T_{\cU A/\cU A'}\big((\bbf{a}_i^{n+1})_{1\le i\le n}\big)\mod \pi^{l'_{n+1}}.
$$
Thus we have proved the first and second part of the statement for $n+1$.

It remains to prove the formulas \eqref{eq:E_p(phi)} and \eqref{eq:H_n} for $n+1$.
But the second one clearly follows from the first one, so we just have
to prove \eqref{eq:E_p(phi)}.
%
Let us suppose that \eqref{eq:E_p(phi)} is true for $n$ and prove it for $n+1$.
We clearly have
\begin{equation}\label{eq:formula to use induction in prop}
\begin{split}
E_p((\bbf{w}_r^{n+1})_{1\le r\le n+1},(\pi^{l_1},\dots,\pi^{l_{n+1}}),
\varphi^{n+1}(\TT))= & \\
E_p\left((\bbf{w}_r^{n+1})_{1\le r\le n},
(\pi^{l_1},\dots,\pi^{l_n}),\varphi^n(\TT)\right)
E_p\left(\bbf{w}_{n+1}^{n+1},\pi^{l_{n+1}}\right.&,
\left.\frac{\varphi_{n+1}(\TT)}{E_p(({\bbf{a}'}_i^{n+1})_{1\le i \le n},
(\pi^{l_i})_{1\le i \le n},\varphi^{n}(\TT))}\right).
\end{split}
\end{equation}
Moreover, by induction
\begin{equation}\label{eq:induction hypothesis in prop}
\begin{split}
E_p\big((\bbf{w}_i^{n+1})_{1\le i\le n},
(\pi^{l_1},\dots,\pi^{l_n}),\varphi^n(\TT)\big)=
& \ E_p\big(T_{\cU A/ \cU A'}(\bbf{w}_i^{n+1})_{1\le i\le n},
(\pi^{l_1},\dots,\pi^{l_n}),\TT\big) \\
& \times
\prod_{2\le r\le n}
G_p\big(U^n_r((\bbf{w}_i^{n+1})_{1\le i \le n}),\pi^{{l'}_r};H_{r-1}\big).
\end{split}
\end{equation}
Now
\renewcommand{\arraystretch}{3}
\begin{align*}
E_p\Big(\bbf{w}_{n+1}^{n+1},\pi^{l'_{n+1}}, &
\frac{\varphi_{n+1}(\TT)}{E_p(({\bbf{a}'}_i^{n+1})_{1\le i \le n},
(\pi^{l_i})_{1\le i \le n},\varphi^{n}(\TT))}\Big) =
\ E_p\Big(\bbf{w}_{n+1}^{n+1},\pi^{l_{n+1}},\frac{H_n-1}{\pi^{l'_{n+1}}}\Big)
\\ & \\
\stackrel{\eqref{eq: G_p vs E_p}}{=}
& \ \tilde{E}_p\Big(\bbf{w}_{n+1}^{n+1},\pi^{l_{n+1}},H_n\Big)
\,G_p\Big(F^{l_{n+1}}\bbf{w}_{n+1}^{n+1},\pi^{l_{n+1}},H_n\Big) \\ & \\
= & \ E_p\Big(T_{\frac{(\bbf{a}_i^{n+1})_{1\le i \le n}- T_{\cU A/\cU A'}
{({\bbf{a}'}_i^{n+1})_{1\le i \le n}}}{\pi^{l_{n+1}}}}\bbf{w}_{n+1}^{n+1},
(\pi^{l_1},\dots,\pi^{l_n});\TT\Big) \\ & \\
& \ \times
E_p\Big(T_{[\pi^{l_{n+1}-l'_{n+1}}]}\bbf{w}_{n+1}^{n+1},
\pi^{l_{n+1}};\frac{T_{n+1}}{E_p((\bbf{a}_i^{n+1})_{1\le i \le n-1},
(\pi^{l_i})_{1\le i \le n},\varphi^{n}(\TT))}\Big) \\ & \\
& \ \times \prod_{r=2}^n
G_p\Big(-T_{\frac{U_r^n(\bbf{a}_i^n)_{1\le i \le n-1}}
{\pi^{l_{n+1}}}}\bbf{w}_{n+1}^{n+1},
\pi^{l_{n+1}};H_{r-1}\Big) \,G_p\Big(F^{l_{n+1}}
\bbf{w}_{n+1}^{n+1},\pi^{l_{n+1}},H_n\Big)
\end{align*}
where in the last equality we have used Equation~\eqref{eq:H_n},
Lemma~\ref{eq: E tilde avec E et operation T},
Lemma~\ref{eq: E tilde avec G et operation T},
Equation~\eqref{eq:G_p additive} and the fact that
$\tilde{E}_p(\WW,\Lambda, H)$ gives
a bilinear group scheme homomorphism
$W_{\AA^1}\times \mathbf{\Lambda}_{\AA^1}\to \mathbf{\Lambda}_{\AA^1}$.
Now using \eqref{eq:formula to use induction in prop} and
\eqref{eq:induction hypothesis in prop} one gets the result.

Finally we remark that if $A/A'\ge 0$ and $A'/A\ge 0$ then
necessarily $A/A'$ and $A'/A$ are unitriangular, as it is very
easy to verify.
\end{proo}

The order relation $\succ$ from the previous proposition induces
an equivalence relation on $\sM_n$:

\begin{defi}
For any $A,A'\in \sM_n$ we write $A\sim A'$ if and only if $A/A'$
is positive and unitriangular.
\end{defi}

This relation characterizes when two matrices in $\sM_n$
define the same filtered group scheme:

\begin{coro} \label{co:bijection}
The map $A\mapsto \cE(A)$ induces an increasing bijection between the set
$\sM_n/\!\sim$ ordered by the relation $\succ$ and the set of isomorphism
classes of filtered group schemes of dimension $n$ ordered by the
existence of a model map.
\end{coro}

\end{noth}


\section{Kummer group schemes} \label{MatrixSS}

In this section, following Sekiguchi and Suwa's approach, we specify
Theorem~\ref{thm 2 SS} for filtered $R$-group schemes containing a
model of $\mu_{p^n}$. The main result (\ref{result3} below) is a
generalization of Theorem~9.4 of \cite{SS1}, which covers the
particular case where the finite flat subgroup is the constant
group scheme $(\ZZ/p^n\ZZ)_R$. As it turns out, the main
difficulty is to find the {\em statement} of the generalized
theorem, for then the proof of \cite{SS1} carries over smoothly.

We point out an important fact: the computation of successive
extensions by groups $\cG^{(l)}$, which is the essence of the existence of
filtered group schemes, proceeds differently when $l>0$ and when $l=0$.
The former case is treated by Theorem~\ref{result3}, and we indicate in
Remark~\ref{rm:l_i=0} how to handle the easier case $l=0$.

\subsection{Finiteness of closures of finite flat subgroups}

Let $\cE=\cE_n$ be a filtered group scheme of type
$(l_1,\dots,l_n)$. Let $\alpha:\cE\to (\GG_m)^n$ be a morphism of
filtered $R$-group schemes which is an isomorphism on the generic fibre.
Let $\Theta^n:(\GG_m)^n\to (\GG_m)^n$ be the morphism defined by
$$
\Theta^n(T_1,\dots,T_n)
=(T_1^p,T_2^pT_1^{-1},\dots,T_{n}^pT_{n-1}^{-1}).
$$
The kernel of $\Theta^n$ is a subgroup isomorphic to $\mu_{p^n,R}$ which
we call the {\em Kummer $\mu_{p^n}$ of $\GG_m^n$}. Via the map
$\alpha$, we can see the Kummer $\mu_{p^n,K}$ as a closed subscheme of
$\cE_K$. We define the {\em pre-Kummer subgroup} $G_n$ as the
scheme-theoretic closure of $\mu_{p^n,K}$ in $\cE$, and we call it the
{\em Kummer subgroup} when it is finite over $R$.
In spite of the notation,  $G_n$ depends on the
choice of $\alpha$ (see Theorem \ref{thm 1 SS}).
If $G_n$ is finite, then the quotient $\cF_n$ is a filtered group
scheme and the quotient map $\Psi^n:\cE_n\to\cF_n$ is an isogeny.
In this case, for each $\mu\in R$ we have a pullback map
$$
(\Psi^n)^*:\Hom_{R/\mu R\mbox{-}\Gr}(\cF_n,\GG_m)
\too \Hom_{R/\mu R\mbox{-}\Gr}(\cE_n,\GG_m).
$$
%
We know by Theorem~\ref{thm 2 SS} that using the deformed Artin-Hasse
exponentials, the groups on both sides may be identified with suitable
kernels of additive operators $U^n$ on Witt vector groups. Once this is done,
Sekiguchi and Suwa express $(\Psi^n)^*$ by a matrix called $\Upsilon^n$.
Let us give some details in the case $n=1$ that initiates the induction.
Then we have $\cE_1\simeq\cG^{(l_1)}$
and the closure of $\mu_{p,K}$ is finite flat if and only
if $v(p)\ge (p-1)l_1$, see e.g. \cite{MRT}, Lemma~5.1.1.
Moreover $\cF_1=\cE_1/G_1\simeq\cG^{(pl_1)}$ and one may check that
the pullback $(\Psi^1)^*$ is expressed by the one-term matrix
$\Upsilon^{1}=(T_{p[\lambda_{1}]/\lambda_{1}^p})$. Note that the operator
$T_{p[\lambda_{1}]/\lambda_{1}^p}$ indeed takes the kernel of $F^{(pl_1)}$
into the kernel of $F^{(l_1)}$, see \cite{MRT}, Lemma~5.2.8. Let us come back to
an arbitrary dimension $n$.

In this setting, we can characterize the situation where the pre-Kummer group
scheme $G_{n+1}$ inside a filtered group scheme $\cE_{n+1}$ is finite and
flat. In the statement below, we will denote by $U^i$ the matrices involved
in the construction of $\cE_n$ like in Theorem~\ref{thm 2 SS}, and
$\bar{U}{}^i$ the matrices involved in the construction of $\cF_n$
(this is the notation of \cite{SS1}). Note that the inductive construction of
$\bar{U}{}^n$ is included in the statement of the theorem via the vectors
$\bbf{u}^n$.

\begin{theo} \label{result3}
Let $n\geq 1$, $\bbf{l}=(l_1,\ldots,l_{n+1})$ with $l_i>0$ for
each $i$, and $A\in\sM^{\bbf{l}}_{n+1}$. Let
$\cE=\cE({A})=(\cE_1,\ldots,\cE_{n+1})$ be the filtered group
scheme of type $\bbf{l}$ defined by $A$. Assume that
$G_n\subset\cE_n$ is finite flat. Then, the following conditions
are equivalent:
\begin{itemize}
\item[{\rm (i)}] $G_{n+1}$ is finite flat,
\item[{\rm (ii)}] $v(p)\ge (p-1)l_{n+1}$ and there exist vectors
$\bbf{u}^{n+1}$ and $\bbf{v}^{n+1}$ in $W^f(R)^n$, with the reduction of
$\bbf{u}^{n+1}$ modulo $\pi^{pl_{n+1}}$ lying in
$\ker(\bar{U}{}^n:\widehat{W}(R/\pi^{pl_{n+1}}R)^n\endto)$,
such that
%
$$
p\bbf{a}^{n+1}_i-\bbf{c}_i^{n+1}-(\Upsilon^n\bbf{u}^{n+1})_i
=\pi^{pl_{n+1}}.\bbf{v}^{n+1}_i
$$
for all $1\leq i\leq n$, where
$\bbf{c}^2=\bbf{c}^2_1=[\pi^{l_1}]\in W(R)$,
$\bbf{c}^{n+1}=(\bbf{a}^{n},[\pi^{l_{n}}])\in W(R)^n$  for $n\geq
2$, with $\bbf{a}^n=(\bbf{a}_i^n)_{1\le i <n}$.
\end{itemize}
In this case,
the filtered group
scheme $\cF_{n+1}=\cE_{n+1}/G_{n+1}$ is obtained  from $\cE_n/G_n$
using the deformed Artin-Hasse exponential defined by $\bbf{u}^{n+1}$.
Moreover, if $\Psi^{n+1}:\cE_{n+1}\to\cF_{n+1}$ is the induced morphism
and if $\mu \in R\setminus\{0\}$ then the morphism
$$
(\Psi^{n+1})^*: \Hom_{R/\mu R\mbox{-}\Gr}(\cF_{n+1},\GG_m)
\to \Hom_{R/\mu R\mbox{-}\Gr}(\cE_{n+1},\GG_m)
$$
is given by
$$
\Upsilon^{n+1}=
\left(
\begin{array}{cccc}
& & & \\
& \Upsilon^n & & T_{\bbf{v}^{n+1}} \\
0 & \dots & 0 & T_{p[\pi^{l_{n+1}}]/\pi^{pl_{n+1}}} \\
\end{array}
\right).$$
\end{theo}

\begin{proo}
We make an induction on $n$. It is convenient to set $\cE_0=\{1\}$ and to
start the induction at $n=0$, in which case the result is known (see e.g.
\cite{MRT}, Lemma~5.1.1). For the last statement we will prove something
more precise. Write $\cE_{n+1}=\cE(A)$ and $\cF_{n+1}=\cE(B)$
with $A=(\bbf{a}_i^j)$ and $B=({\bbf{u}}_i^j)$ in $\sM_{n+1}$
with $\bbf{a}_i^i=[\pi^{l_i}]$ and $\bbf{u}_i^i=[\pi^{pl_i}]$
for $i=1,\dots,n+1$. 
Let
$$K_0:= (1+\pi^{l_1} T_1)^p=E_p(p[\pi^{l_1}],\pi^{l_1},T_1)$$
and for $r\ge 1$ let
$$
K_r:=\frac{(E_p(\bbf{a}^{r+1}, (\pi^{l_i})_{1\le i \le
r},\TT)+\pi^{l_{r+1}}T_{r+1})^p}
{(E_p(\bbf{a}^r,(\pi^{l_i})_{1\le i\le r-1},
\TT)+\pi^{l_{r}}T_r)\,E_p({\bbf{u}}^{r+1},
(\pi^{l_i})_{1\le i\le r}, \Psi^r(\TT))} \in
R[T_1,\dots,T_{r+1}].
$$
Given $\bbf{w}^{n}\in W(R)^{n}$ whose reduction belongs to
$\ker(\bar{U}{}^{n}:\hat{W}(R/\mu R)^{n}\endto)$, we will prove that
\begin{equation}\label{eq:E_p(Psi)}
\begin{split}
E_p(\bbf{w}^{n},(\pi^{pl_1},\dots,\pi^{pl_{n+1}}),\Psi^{n}(\TT)) = & \\
E_p(\Upsilon^n(\bbf{w}^{n}),(\pi^{l_1},\dots,\pi^{l_{n+1}}),\TT)
& \prod_{1\le r\le n}G_p(\bar{U}{}^n_r(\bbf{w}^{n}),\pi^{pl_r},K_{r-1}) ,
\end{split}
\end{equation}
where $\bar{U}{}^n_r$ is the $r$-th row of $\bar{U}{}^n$.
Since $G_p(\bar{U}{}^n_r(\bbf{w}^{n}),\pi^{pl_{r}};K_{r-1})\in 1+TR[[T]]$
for each $r\ge 1$ (\cite{SS1}, Prop.~9.3), Equation (\ref{eq:E_p(Psi)})
implies
\begin{equation}\label{eq:K_n}
K_r= \frac{E_p(p\bbf{a}^r-\bbf{c}^r-\Upsilon^r\bbf{u}^r,
(\pi^{l_1},\dots,\pi^{l_r});\TT)\,E_p(p[\pi^{l_{r+1}}],\pi^{l_{r+1}};
\frac{T_{r+1}}{E_p({\bbf{a}}^{r+1},
(\pi^{l_i})_{1\le i\le r}; \TT)})}
{\prod_{i=1}^rG_p(U_i^r(\bbf{u}^r),\pi^{pl_{r+1}};K_{i-1})}.
\end{equation}
For $n=1$ the formula \eqref{eq:E_p(Psi)} follows
from \ref{eq: E tilde avec E et operation T} and \eqref{eq: G_p vs E_p}.
We now suppose that the theorem and the formula \eqref{eq:E_p(Psi)} are
true for $n-1$ and we prove them for $n$.
We do this in three steps (a)-(b)-(c).

\bigskip

\noindent (a) We prove that (i) is equivalent to (ii). Among the
objects constructed inductively at the same time as the filtered
groups $\cE_n,\cF_n$, we consider the polynomials $D_r,D'_r$
(truncated exponentials associated respectively to $\cE_n$ and
$\cF_n$) and the isogenies $\Psi^r:\cE_r\to \cF_r$, for $1\leq
r\leq n-1$. We also introduce the notation:
$$
C_{n+1}=C_{n+1}(T_1,\dots,T_{n+1})
:=(D_n(T_1,\dots,T_n)+\pi^{l_{n+1}}T_{n+1})^p
(D_{n-1}(T_1,\dots,T_{n-1})+\pi^{l_{n}}T_{n})^{-1}.
$$
We have $K[G_{n+1}]=K[G_n][T_{n+1}]/(C_{n+1}-1)$. Assume
that $G_{n+1}$ is finite over $R$; then it is finite over $G_n$. It follows
(\cite{Ei}, Prop.~4.1) that $C_{n+1}-1\equiv 0 \mod \pi^{pl_{n+1}}$ and
\begin{equation}\label{Kumsub}
R[G_{n+1}]=R[G_n][T_{n+1}]/\Big(\frac{C_{n+1}-1}
{\pi^{pl_{n+1}}}\Big).
\end{equation}
In particular $C_{n+1}$, seen as an element of
$\Hom_{R/\pi^{pl_{n+1}}R}(G_n,\GG_m)$, is the trivial morphism.
If we apply the functor $\Hom_{R/\pi^{pl_{n+1}}R}(-,\GG_m)$ to the
short exact sequence
$$
0\too G_n\stackrel{i_n}{\too} \cE_n \stackrel{\Psi^n}{\too} \cF_n \too 0,
$$
we obtain a long exact sequence
$$
0 \too \Hom_{R/\pi^{pl_{n+1}}R}(\cF_n,\GG_m)\on{(\Psi^n)^*} \too
\Hom_{R/\pi^{pl_{n+1}}R}(\cE_n,\GG_m) \stackrel{i_n^*}{\too}
\Hom_{R/\pi^{pl_{n+1}}R}(G_n,\GG_m) \too \dots
$$
As we noticed before, the element $C_{n+1}$ lives in
$\ker(i_n^*)$ and hence is equal to $(\Psi^n)^*(D'_n)$ for some
$D'_n\in \Hom_{R/\pi^{pl_{n+1}}}(\cF_{n+1},\GG_m)$. Now we use the
description of groups of homomorphisms from a filtered group
scheme to $\GG_m$ in terms of vectors, as given by the third point of
Theorem~\ref{thm 2 SS}. Let $\bbf{u}^{n+1}\in W(R)^n$ be a lift in
$W(R)^n$ of a vector corresponding to $D'_n$. The equality
$$
C_{n+1}=(\Psi^n)^*(D'_n)
$$
translates to $n$ equalities
$p\bbf{a}^{n+1}_i-\bbf{c}_i^{n+1}=(\Upsilon^n\bbf{u}^{n+1})_i$
in $W(R/\pi^{pl_{n+1}}R)^n$, for $1\le i\le n$. Lifting this to $W(R)^n$ shows
that (ii) holds. Conversely, if (ii) holds then $C_{n+1}$ is of
the form $(\Psi^n)^*(D'_n)$ for some $D'_n$, hence it has trivial image under
$i_n^*$. Thus $C_{n+1}\equiv 1 \mod \pi^{pl_{n+1}}$ and the expression
(\ref{Kumsub}) defines a finite flat group scheme $G_{n+1}$ over $R$.

\bigskip

\noindent (b) Let $\cF_{n+1}=\cE_{n+1}/G_{n+1}$. Now we prove that $\cF_{n+1}=\cF'_{n+1}$ where $\cF'_{n+1}$
is the filtered group scheme obtained from $\cF_{n}$ using the vector
$\bbf{u}^{n+1}$. The $R$-algebra of $\cF'_{n+1}$ is
$$
R[\cF_n][T_{n+1},1/(D'_{n}+\pi^{pl_{n+1}}T_{n+1})].
$$
where $D'_n\in R[T_1,\ldots,T_n]$ is the truncation of
$E_p(\bbf{u}^{n+1})\prod_{j=1}^n E_p(\bbf{u}^{n+1}_j,\pi^{l_j},
T_j/E_p(\bbf{u}^{i-1},(\pi^{l_i}),\TT))$ as defined in Theorem
\ref{thm 2 SS}. Let $\Psi_1,\ldots,\Psi_n$ be the polynomials
defining the isogeny $\Psi^n:\cE_n\rightarrow \cF_n$. Let
$$
\Psi_{n+1}(T_1,\ldots,T_n)= \frac{1}{\pi^{pl_{n+1}}}
\left(\frac{(D_n(T_1,\dots,T_n)+\pi^{l_{n+1}}T_{n+1})^p}
{D_{n-1}(T_1,\dots,T_{n-1})+\pi^{l_{n}}T_n}
-D'_n(\Psi_1(\TT),\ldots,\Psi_n(\TT))\right).
$$
Then the morphism $R[\cF_{n+1}]\rightarrow R[\cE_{n+1}]$,
$T_i\mapsto \Psi_i(T_1,\ldots,T_{i})$ defines an
isogeny $\cE_{n+1}\to \cF'_{n+1}$ with kernel $G_{n+1}$. Therefore
$\cF_{n+1}$ is isomorphic to $\cF'_{n+1}$ as a filtered group scheme.

\bigskip

\noindent (c) We now prove the formula \eqref{eq:E_p(Psi)}.
We have
\begin{equation}\label{eq:formula to use induction in Thm}
\begin{split}
E_p\left(\bbf{w}^{n+1},(\pi^{pl_1},\dots,\pi^{pl_{n+1}}),\Psi^{n+1}(\TT)\right) = & \\
E_p\left((\bbf{w}_r^{n+1})_{1\le r\le n},(\pi^{pl_1},\dots,
\pi^{pl_{n}}),\Psi^n(\TT)\right)
& \, E_p\left(\bbf{w}_{n+1}^{n+1},\pi^{pl_{n+1}},
\frac{\Psi_{n+1}(\TT)}{E_p({\bbf{u}}^{n+1},
(\pi^{l_i})_{1\le i\le n}; \Psi^n(\TT))}\right). \\
\end{split}
\end{equation}
By the induction hypothesis we have
\begin{equation}\label{eq:induction hypothesis in Thm}
\begin{split}
E_p\left((\bbf{w}_r^{n+1})_{1\le r\le n},(\pi^{pl_1},\dots,
\pi^{pl_{n}}),\Psi^{n}(\TT)\right)
= & \\
E_p\left(\Upsilon^n\bbf{w}^{n+1},(\pi^{l_1},\dots,\pi^{l_{n}}),\TT\right)
\prod_{1\leq r\leq n}
& G_p\left(U_r^n((\bbf{w}^{n+1}_s)_{1\leq s\leq n}),\pi^{pl_{n+1}},K_{r-1}\right).
\end{split}
\end{equation}
Moreover we have
\begin{align*}
E_p\bigg(\bbf{w}_{n+1}^{n+1},\pi^{pl_{n+1}},
& \frac{\Psi_{n+1}(\TT)}{E_p({\bbf{u}}^{n+1},
(\pi^{l_i})_{1\le i\le n}; \Psi^n(\TT))}\bigg)
= \ E_p\left(\bbf{w}_{n+1}^{n+1},\pi^{pl_{n+1}},\frac{K_n-1}{\pi^{pl_{n+1}}}\right) \\
& \\
\stackrel{\eqref{eq: G_p vs E_p}}{=} &
\ \tilde{E}_p\Big(\bbf{w}_{n+1}^{n+1},\pi^{pl_{n+1}},K_n\Big)
\,G_p\Big(F^{(pl_{n+1})}\bbf{w}_{n+1}^{n+1},\pi^{pl_{n+1}},K_n\Big) \\ & \\
= & \ E_p\Big(T_{\frac{(p\bbf{a}^n-\bbf{c}^{n-1}-\Upsilon^n\bbf{u}^{n+1})}{\pi^{pl_{n+1}}}}\bbf{w}_{n+1}^{n+1},(\pi^{l_1},\dots,\pi^{l_n});\TT\Big) \\ & \\
& \ \times E_p\left(T_{\frac{p[\pi^{l_{n+1}}]}{\pi^{pl_{n+1}}}}\bbf{w}_{n+1}^{n+1},\pi^{l_{n+1}};\frac{T_{n+1}}{E_p(\bbf{a}^{n+1}, (\pi^{l_i})_{1\le i \le
n},\TT)}\right) \\ & \\
& \ \times \prod_{1\le r\le n}
G_p\bigg(-T_{\frac{U_r^n\bbf{u}^{n+1}}{\pi^{pl_{n+1}}}}\bbf{w}^{n+1}_{n+1},\pi^{pl_{n+1}};K_{r-1}\bigg) \\ & \\
& \ \times
G_p\left(F^{(l_{n+1})}\bbf{w}_{n+1}^{n+1},\pi^{pl_{n+1}};K_n\right)
\end{align*}
where in the last equality we have used Equation~\eqref{eq:K_n},
Lemma~\ref{eq: E tilde avec E et operation T},
Lemma~\ref{eq: E tilde avec G et operation T},
Equation~\eqref{eq:G_p additive} and the bilinearity of
$\tilde{E}_p(\WW,\Lambda, H):W_{\AA^1}\times
\mathbf{\Lambda}_{\AA^1}\to \mathbf{\Lambda}_{\AA^1}$,
see~\ref{eq: E tilde}.
Now using Equations~\eqref{eq:formula to use induction in Thm} and
\eqref{eq:induction hypothesis in Thm} one gets the result.
\end{proo}

\begin{defi}
Let $\cE(A)=(\cE_1,\ldots,\cE_n)$ be a filtered group scheme.
We say that {\em $A$ satisfies the integrality conditions} if for any
$1\leq i\leq n$, the upper left square submatrix $\cU^{n-i}A$ of $A$
satisfies the conditions (ii) in Theorem \ref{result3} applied to
an $i$-dimensional matrix.
\end{defi}

In other words, $A$ satisfies the integrality conditions if and only if
the pre-Kummer subgroups $G_i$ are finite flat in $\cE_i$, for $1\leq i\leq n$.
From the proof of Theorem \ref{result3}, we deduce an explicit formula
for these models of $\mu_{p^n,K}$.

\begin{coro}
Let $\cE_n=\cE(A)$ be a filtered group scheme given by a family of
parameters $A=(\bbf{a}_i^j)$ satisfying the integrality conditions.
Let ${D}_1(T_1)\in R[T_1]$ be any lifting of
$E_p(\bbf{a}_1^2,\pi^{l_1},T_1)\mod \pi^{l_2},\ldots,$
and ${D}_{n-1}(T_1,\ldots,T_{n-1})\in R[T_1,\ldots,T_{n-1}]$
be any lifting of
$$
E_p(\bbf{a}^{n},(\pi^{l_1},\ldots,\pi^{l_{n-1}}), T_1,\ldots,T_{n-1})
\mod\pi^{l_n}
$$
as defined in Theorem~\ref{thm 2 SS}.
Then $G_n$ is a finite flat $R$-group scheme, defined in affine $n$-space
in coordinates $T_1,\dots,T_n$ by the $n$ equations:
\begin{align*}
\frac{(1+\pi^{l_1}T_1)^p-1}{\pi^{pl_1}} \ , \
& \frac{({D}_1+\pi^{l_2}T_2)^p(1+\pi^{l_1}T_1)^{-1}-1}{\pi^{pl_2}}
\ , \ \dots \\
& \dots \ , \
\frac{({D}_{n-1}+\pi^{l_n}T_n)^p({D}_{n-2}+\pi^{l_{n-1}}T_{n-1})^{-1}-1}
{\pi^{pl_n}}. \\
\end{align*}
\end{coro}

\begin{proo}
This is a translation of Theorem~\ref{result3}.
\end{proo}

\begin{rema} \label{rm:l_i=0}
In the statement of Theorem~\ref{result3}, it is assumed that $l_i>0$
for all $i$. Here is what to do so as to obtain a description of all
Kummer group schemes, including the case where some $l_i$ vanish.
We make some preliminary observations. First, as it is easy to see from
the case $n=2$, the type of a Kummer group is necessarily ordered:
$l_1\ge l_2\ge\dots\ge l_n$ (see \ref{consequences utiles}). Second,
there are no nontrivial extensions of a filtered group scheme by
$\GG_m$ (see~\cite{SS1}, Prop.~3.1). Third, it is easy to see that
the only Kummer group scheme of $(\GG_m)^n$ with type
$(0,\dots,0)$ is $\mu_{p^n}$. After these preliminaries, it
remains to see how to describe Kummer groups of type $\bbf{l}$
with $l_1\ge\dots\ge l_r>l_{r+1}=\dots=l_n=0$. Such a Kummer group
$G$ lies in a filtered group $\cE=\GG_m^{n-r}\times\cE'(A')$ where
$\cE'(A')$ is filtered of type $\bbf{l}'=(l_1,\dots,l_r)$. We
define $\cE(A):=\cE$ when
$$
A=\left(
\begin{array}{cccc}
& A' & 0 &\\
& & & \\
& 0 & V & \\
\end{array}
\right)
$$
with $V$ unipotent in $M_{n-r}(W^f(R))$ (therefore equivalent to
the identity, since invertible: use Lemma \ref{lem:Ta injective}).
Moreover, $G$ is an extension of a finite flat Kummer group $G'$
of $\cE'$ by $\mu_{p^{n-r}}$. Using the same argument as in
\cite{To2}, Prop.~3.6 there is an exact sequence:
$$
0\too \ZZ/p^r\ZZ\too \Ext^1(G',\mu_{p^{n-r}})\too H^1(S,(G')^\vee)\too 0
$$
where $S=\Spec(R)$ and $(G')^\vee$ is the Cartier dual of $G'$. Then with
the same proof as in \cite{To2}, Cor.~3.20 we see that the Kummer
subgroups of $\cE$ are given by the image of $1\in \ZZ/p^r\ZZ$.
They have the following ring of functions:
$$
R[G]=\frac{R[G'][T_{r+1}]}{(T_{r+1}^{p^{n-r}}(D_{r-1}+\pi^{l_r}T_r)^{-1}-1)}.
$$
\end{rema}

\subsection{Matricial translation of the integrality conditions}
\label{ss:matrices_for_Kummer}

We translate the previous results on Kummer group schemes in terms of
matrices, in order to emphasize the formal similarities with the
classification of models of $\mu_{p^n}$ by their Breuil-Kisin lattices
in~\ref{classification_BK}.
Up to now in Sections \S\ref{secSS} and \S\ref{MatrixSS}, we followed
Sekiguchi and Suwa's notation $a_i^j$ for the entries of matrices.
However, in order to compare the parameters with those of the Breuil-Kisin
classification, we will henceforth write $a_{ij}$. In the statement of the
following result, we use the operator $\cP$ introduced in \ref{nota:E(u)A}.

\begin{theo}\label{result4}
Let $G$ be a Kummer group scheme. Then there exists
$\bbf{l}=(l_1,\dots,l_n)\in(\NN)^n$ and upper triangular matrices
$$
A=\left(
\renewcommand{\arraystretch}{1.5}
\begin{array}{ccccc}
[\pi^{l_1}] & \bbf{a}_{12} & \bbf{a}_{13} & \dots & \bbf{a}_{1n} \\
& [\pi^{l_2}] & \bbf{a}_{23} & & \bbf{a}_{2n} \\
& & \ddots & \ddots & \vdots \\
& & & [\pi^{l_{n-1}}] & \bbf{a}_{n-1,n} \\
0 & & & & [\pi^{l_n}]\\
\end{array}
\right)
\ , \quad
B=\left(
\renewcommand{\arraystretch}{1.5}
\begin{array}{ccccc}
[\pi^{pl_1}] & \bbf{b}_{12} & \bbf{b}_{13} & \dots & \bbf{b}_{1n} \\
& [\pi^{pl_2}] & \bbf{b}_{23} & & \bbf{b}_{2n} \\
& & \ddots & \ddots & \vdots \\
& & & [\pi^{pl_{n-1}}] & \bbf{b}_{n-1,n} \\
0 & & & & [\pi^{pl_n}] \\
\end{array}
\right)
$$
with entries in $W^f(R)$, satisfying
\begin{trivlist}
\itemm{1} $F(A)/A\ge 0$, $F(B)/B\ge 0$,
\itemm{2} $(pA-\cP\cU A)/B\ge 0$,
\end{trivlist}
such that $G$ is the kernel of an isogeny $\cE(A)\too \cE(B)$.
Moreover when $A$ is chosen, $B$ is unique up to the equivalence relation $\sim$.

If there exists $A$ and $B$ as above then $\cE(A)$ contains a finite and flat Kummer subgroup scheme $G$ and $\cE/G\simeq \cE(B)$.
\end{theo}

\begin{proo}
Let us first suppose that $(l_1,\dots,l_n)=(0,\dots,0)$. Then the
unique Kummer scheme of type $(0,\dots,0)$ is $\mu_{p^n}$. In this
case any matrix as in the statement is unipotent and so is equivalent
to the identity, since it is invertible. Therefore such a matrix gives
rise to the group scheme $\mu_{p^n}$. Conversely $A=B=Id$ satisfy the
conditions of the theorem.

We will now suppose that all the $l_i$ are all strictly positive.
The case with some $l_i$ equal to zero can be deduced form this
one using Remark \ref{rm:l_i=0}.

By definition of a Kummer subgroup, it is embedded in a filtered
group scheme $\cE(A)$, for some $A\in \sM_n$ of type $\bbf{l}$.
Moreover by \ref{result3} it is easy to see that
$$
pA-\cP\cU A=V_n\star_TB_n
$$
where the matrices $B_n,V_n$ are defined by induction:
$$
B_{n+1}=\left(
\begin{array}{cccc}
& & & \\
& B_n& & \bbf{u}^{n+1}\\
& & & \\
0 & \dots & 0 & [\pi^{pl_{n+1}}] \\
\end{array}
\right)
$$
and
$$
V_{n+1}=\left(
\begin{array}{cccc}
& & & \\
& V_n & & \bbf{v}^{n+1}\\
& & & \\
0 & \dots & 0 & p[\pi^{l_{n+1}}]/\pi^{pl_{n+1}} \\
\end{array}
\right)
$$
with $\bbf{u}^{n+1}$ and $\bbf{v}^{n+1}$ as in the statement of
Theorem \ref{result3}. It follows from the same theorem that
$\cE(A)/G\simeq\cE(B)$. Similar argument for the converse.
\end{proo}

\begin{rema}\label{rm:A_and_FA}
One significant difference between this theorem and
Theorem \ref{classification_BK} is that here we do not provide a normal form,
or distinguished choice, for a matrix~$A$ defining a Kummer group scheme.
Theorem \ref{classification_BK} suggests that maybe one could choose a pair
$(A,B)$ of the form $(A,F(A))$. This is true for instance for $n=2$ and at least in some cases for $n=3$, as we will
see in the next section.
\end{rema}

This theorem should be seen as the analogue in Sekiguchi-Suwa Theory
of Theorem~\ref{classification_BK} in Breuil-Kisin Theory.


\section{Computation of Kummer group schemes for $n=3$} \label{SSn3}

In this section, we apply the general theory to compute some
Kummer group schemes for $n=3$, that is to say, those models of
$\mu_{p^3}$ constructed using Sekiguchi-Suwa Theory. From the
start, we see that the complexity of the computations with Witt
vectors is a serious obstacle. In fact, the difficulty increases
with the number of nonzero coefficients of the vectors. It is
therefore interesting to know if in Theorem~\ref{result4} we can
choose matrices $A$ and $B$ with "short" Witt vector entries. The
results of \cite{To2} show that in the case $n=2$, any Kummer
group scheme may be described by matrices $A,B$ such that $A$ has
{\em Teichm\"uller entries}. We could not settle the question
whether this is possible for all $n$, but in our opinion it is not
very likely. It is much more plausible that $A$ may be chosen with
entries of bounded length; more precisely, it seems reasonable to
hope that each Kummer model of $\mu_{p^n}$ may be defined by a
matrix $A$ all whose entries are Witt vectors of length at most
$n-1$. However, it is probably impossible to make simple choices
{\em simultaneously for $A$ and~$B$}, which means that one should
not make a priori assumptions on $B$. Because of these remarks,
here we compute the Kummer groups defined by matrices $A,B\in
\sM_3$ such that $A$ has Teichm\"uller entries and $B$ is
arbitrary.

Even though this is not essential, it will simplify matters to
assume throughout that $p\ge 3$. Recall that $R$ is a discrete
valuation ring of characteristic $0$, residue characteristic $p$,
uniformizer~$\pi$,
valuation $v$ and $v(p)=e$. For $l\geq 1$, we denote again by $v$
the induced valuation on $R/\pi^lR$. Whenever the context does not
allow confusions, we keep the same notation for a Witt vector
$\bbf{a}=(a_0,a_1,\ldots,a_i\ldots)\in W(R)$ and its image in
$W(R/\pi^lR)$. Following the comments at the beginning
of~\ref{ss:matrices_for_Kummer} and the notation in~\ref{result4},
in this section we will write $\bbf{a}_{ij}$ the entries of the
matrices.

\subsection{Two lemmas}

We collect two easy lemmas for future reference.

\begin{lemm} \label{lm:kerF}
Let $l\ge 1$ be an integer. Then the following statements hold.
\begin{trivlist}
\itemn{1} If $pe\ge (p-1)l$, then
for any $\bbf{a}\in \ker(F:\hat W(R/\pi^lR)\endto)$ we have
$v(a_i)\ge l/p$ for all $i\ge 0$. \itemn{2} For all
$\bbf{a},\bbf{b}\in \hat W(R/\pi^lR)$ such that $v(a_i)\ge l/p$
and $v(b_i)\ge l/p$ for all $i\ge 0$, we have
$\bbf{a}+\bbf{b}=(a_0+b_0,a_1+b_1,\dots)$.
\end{trivlist}
\end{lemm}

\begin{proo}
This is Lemma~2.4 of \cite{To2}. More precisely, (1) and (2) are
proven in the proof of {\em loc. cit.} and the assertion of
Lemma~2.4 itself is a combination of these statements.
\end{proo}

\begin{lemm} \label{lemmar}
Let $\XX=(X_0,X_1,X_2,\dots)$ and $\YY=(Y_0,Y_1,Y_2,\dots)$ be
sequences of indeterminates and $S_0,S_1,S_2,\dots$ the
polynomials giving Witt vector addition (see~\ref{ss:p_adic_exp}).
\begin{trivlist}
\itemn{1} If the variables $X_i,Y_i$ are given the weight $p^i$,
then the polynomial $S_n(\XX,\YY)\in\ZZ[\XX,\YY]$ is homogeneous
of degree $p^n$. \itemn{2} We have $S_0(\XX,\YY)=X_0+Y_0$,
$$
S_1(\XX,\YY) =S_0(X_1,Y_1)+\sigma_1(X_0,Y_0) \ \mbox{ with } \
\sigma_1(X_0,Y_0) =\frac{X_0^p+Y_0^p-(X_0+Y_0)^p}{p}
$$
and
$$
S_2(\XX,\YY)=S_0(X_2,Y_2)+\sigma_1(X_1,Y_1)+
\sigma_1\Big(X_1+Y_1,\sigma_1(X_0,Y_0)\Big)+\sigma_2(X_0,Y_0)
$$
with
$$
\sigma_2(X_0,Y_0)=\frac{X_0^{p^2}+Y_0^{p^2}-(X_0+Y_0)^{p^2}
-p\sigma_1(X_0,Y_0)^p}{p^2}.
$$
\itemn{3} We have $\sigma_i(X,-Y)=\sigma_i(X,Y-X)$ and
$\sigma_i(X,-Y)=\sigma_i(X,Y-X)$. \itemn{4} For any $l\geq 1$,
$a,b\in R/\pi^lR$ and $i=1,2$ we have:
$$v(\sigma_i(a,b))\geq \min\Big((p^i-1)v(a)+v(b),(p^i-1)v(b)+v(a)\Big).$$
\itemn{5} If $p\ge 3$, then in the ring $W(\ZZ)$ we have
$$
p=(p,1-p^{p-1},\epsilon_2p^{p-1},\epsilon_3p^{p-1},\epsilon_4p^{p-1},\dots)
$$
where $\epsilon_2,\epsilon_3,\epsilon_4,\dots$ are principal
$p$-adic units.
\end{trivlist}
\end{lemm}

\begin{proo}
(1) is obvious and well-known, (2) is a simple computation, (3) is
proven in~\ref{lm:easy_formulas}, (4) follows from (2) with the
help of the binomial theorem, and (5) is proven for example
in Lemma~5.2.1 of~\cite{MRT}.
\end{proo}

\subsection{Computations for $n=2$}\label{Compn2}

As already said, the case $n=1$ is well-known (see for instance
\cite{MRT} Lemma 5.1.1). All the models of $\mu_{p,K}$ are given
by certain group schemes, called $G_{\pi^l,1}$ with $l\in \NN$
such that $\frac{e}{p-1}\ge l$. If $l=0$ we obtain the group
scheme $\mu_{p,R}\subset {\GG}_m$. One proves (see \cite{To2} \S~1)
that there exists a model map between $G_{\pi^l,1}$ and
$G_{\pi^m,1}$ if and only if $l \ge m$. Using
Lemma~\ref{lemma_inclusion_of_lattices}, this implies that for models of
$\mu_{p,K}$, the covariant equivalence of Kisin described
in \S~\ref{sec:Kisin equivalence} is given by $G_{\pi^l,1}\mapsto u^lk[[u]]$.
We recall that the matrix associated to the
Breuil-Kisin module $u^lk[[u]]$ is the $1\times 1$ matrix $(u^l)$.

We now consider the case $n=2$. First we recall the following
lemma.

\begin{lemm} \label{l_1 > l_2}
If $G$ is a  Kummer group scheme of type $(l_1,l_2)$
then $l_1\ge l_2$.
 \end{lemm}

 \begin{proo}
If $G$ is of type $(l_1,l_2)$, it is an extension of $G_{\pi^{l_1},1}$ by
$G_{\pi^{l_2},1}$ so the result follows from \cite{To2}, Lemma 3.2.
Another way to see this is to argue that by~\ref{k_i_p} there is a model
map $G_{\pi^{l_1},1}\to G_{\pi^{l_2},1}$ and this forces $l_1\ge l_2$, as recalled
above. Finally, in the case $l_1,l_2>0$ one
can also obtain the lemma using Theorem \ref{result4}.
 \end{proo}

Now fix $l_1,l_2>0$. We consider a matrix 
$$
A=\left(
\renewcommand{\arraystretch}{1.5}
\begin{array}{cc}
[\pi^{l_1}] & \bbf{a}_{12} \\
0 & [\pi^{l_2}] \\
\end{array}\right)
\in M_2(W^f(R)).
$$

\begin{lemm} \label{lm:belonging_to_sM_2}
The condition $A\in \sM_2$, that is to say $F(A)/A\ge 0$, is
equivalent to the congruence $F^{(l_1)}(\bbf{a}_{12})\equiv 0\mod
\pi^{l_2}$. Moreover let
$$
A'=\left(
\renewcommand{\arraystretch}{1.5}
\begin{array}{cc}
[\pi^{l_1}] & \bbf{a}'_{12} \\
0 & [\pi^{l_2}] \\
\end{array}\right)
\in M_2(W^f(R)).
$$
Then
 $\bbf{a}'_{12}\equiv \bbf{a}_{12}\mod \pi^{l_2}$ if and only if
 $A'\in \sM_2$ and $\cE(A)\simeq \cE(A')$.
\end{lemm}

\begin{proo}
By definition we have $F(A)/A\ge 0$ if and only if there exists a
positive matrix
$$
C=\left(
\renewcommand{\arraystretch}{1.5}
\begin{array}{cc}
[\pi^{(p-1)l_1}] & \bbf{c}_{12} \\
0 & [\pi^{(p-1)l_2}] \\
\end{array}\right)
$$
such that $F(A)=C\star_T A$. The equality of entries in position
$(1,2)$ gives $F^{(l_1)}(\bbf{a}_{12})=\pi^{l_2}.\bbf{c}_{12}$.
This is equivalent to $F^{(l_1)}(\bbf{a}_{12})=0$ in
$W(R/\pi^{l_2}R)$, which is the first assertion. In order to prove
the second assertion, let $A'\in M_2(W^f(R))$ be as in the statement. If
$\bbf{a}'_{12}\equiv \bbf{a}_{12}\mod \pi^{l_2}$, then:
\begin{trivlist}
\itemm{i} $F^{(l_1)}(\bbf{a}'_{12})\equiv F^{(l_1)}(\bbf{a}_{12})
\equiv 0\mod \pi^{l_2}$,
\itemm{ii} there exists $\bbf{r} \in W(R)$ such that
$\bbf{a}'_{12}=\bbf{a}_{12}+ \pi^{l_2}.\bbf{r}$, so
$A'=D\star_T A$ with $D={1 \ \bbf{r} \choose 0 \ 1}$.
\end{trivlist}
By (i) and the first assertion of the lemma we have $A'\in \sM_2$,
and by (ii) and Prop.~\ref{pp:pos_mat_and_model_maps} we have
$\cE(A')\simeq \cE(A)$. Conversely if $A'\in \sM_2$ and
$\cE(A)\simeq \cE(A')$, then by Prop.~\ref{pp:pos_mat_and_model_maps}
there exists a unitriangular matrix $D$ as above such that $A'=D\star_T A$.
It follows that $\bbf{a}'_{12}=\bbf{a}_{12}+ \pi^{l_2}.\bbf{r}$
and $\bbf{a}'_{12}\equiv \bbf{a}_{12}\mod \pi^{l_2}$.
\end{proo}

Now we use Theorem \ref{result4} in order to tell exactly when $A$
gives rise to a model of $\mu_{p^2}$, in the case $A$ is a matrix
with Teichm\"uller entries.

\begin{prop} \label{coroTo}
Let
$$
A:=\left(
\renewcommand{\arraystretch}{1.5}
\begin{array}{cc}
[\pi^{l_1}] & [{a}_{12}] \\
0 & [\pi^{l_2}] \\
\end{array}\right)
$$
be a matrix with Teichm\"uller entries and $l_1,l_2\ge 0$.
Then $A$ belongs to $\sM_2$ and $\cE(A)$ contains a finite flat
Kummer subgroup $G$ if and only if
\begin{trivlist}
\itemm{i} $\frac{e}{p-1}\ge l_1 \ge l_2$,
\itemm{ii} $a_{12}^p\equiv 0\mod \pi^{l_2}$, and
\itemm{iii} $pa_{12}-\pi^{l_1}-\frac{p}{\pi^{(p-1)l_1}}(a_{12})^p
\equiv 0 \mod\pi^{pl_2}$.
\end{trivlist}
In such a case, we have  $\cE(A)/G\simeq \cE(F(A))$. Moreover, if
we set $D_1(T)=\sum_{k=0}^{p-1}a_{12}^kT^k/k!$, then
$$
G=\Spec\left(R[T_1,T_2]\Big/\Big(\frac{(1+\pi^{l_1}T_1)^p-1}{\pi^{l_1p}},
\frac{(D_1(T_1)+\pi^{l_2}T_2)^p(1+\pi^{l_1}T_1)^{-1}-1}
{\pi^{pl_2}} \Big)\right).
$$
Finally $G$ depends only on the reduction of $a_{12}$ modulo
$\pi^{l_2}$.
\end{prop}

\begin{proo}
If $l_1=0$ or $l_2=0$, the result comes from the case $n=1$
(see \ref{rm:l_i=0}) and is easy so we suppose that $l_1,l_2>0$.
By Theorem \ref{result3} and Lemmas~\ref{l_1 > l_2} and
\ref{lm:belonging_to_sM_2}, the matrix $A$ gives  a
finite flat Kummer group $G$ if and only $\frac{e}{p-1}\ge
l_1\ge l_2$, $F^{(l_1)}([a_{12}])\equiv 0\mod \pi^{l_2}$ and there
exists a vector $\bbf{u}=\bbf{u}_{12}\in W^f(R)$ with reduction in
$\ker(F^{({pl_1})}:\hat{W}(R/\pi^{pl_2}R)\endto)$ such that:
\begin{equation}\label{eq:equation n=2}
p[a_{12}]-[\pi^{l_1}]=T_{p[\pi^{l_1}]/\pi^{pl_1}}(\bbf{u}) \mod
\pi^{pl_2}.
\end{equation}
Since $l_1\ge l_2$, the congruence $F^{(l_1)}([a_{12}])\equiv 0\mod \pi^{l_2}$
is equivalent to $a_{12}^p\equiv 0\mod \pi^{l_2}$.
It remains only to prove that \eqref{eq:equation n=2} and (iii)
are equivalent equations.

First, let us consider the left hand side of the congruence (iii).
Using the expression $p=(p,1-p^{p-1},\dots)\in W(R)$ recalled
in~\ref{lemmar}(5) and the minoration $v(a_{12})\geq l_2/p$, we
find:
$$
p[a_{12}]=\big(pa_{12},(1-p^{p-1})(a_{12})^p,0,0,\ldots\big)\in
W(R/\pi^{pl_2}R).
$$
Since $p\geq 3$, we have $v(p^{p-1})=(p-1)e\geq (p-1)^2l_2\ge
pl_2$. Thus in fact:
$$
p[a_{12}]=(pa_{12},(a_{12})^p,0,0,\ldots) \mod \pi^{pl_2}.
$$
Using Lemma \ref{lm:kerF}(2) 
we obtain \begin{equation}
p[a_{12}]-[\pi^{l_1}]=(pa_{12}-\pi^{l_1},(a_{12})^p,0,0,\ldots)
 \mod \pi^{pl_2}.
\end{equation}
Now let us turn to the right hand side. We set
$\omega:=\frac{p}{\pi^{(p-1)l_1}}$. In the same way as before, we
obtain
$$
\frac{p[\pi^{l_1}]}{\pi^{pl_1}} =\left(\omega,1,0,0,\dots\right)
\mod \pi^{pl_2}.
$$
It follows that $T_{p[\pi^{l_1}]/\pi^{pl_1}}\bbf{u}=
[\omega]\bbf{u}+V\bbf{u}$ in $W(R/\pi^{pl_2}R)$. Now note that in
$W(R/\pi^{pl_2}R)$ we have $\ker(F^{(pl_1)})=\ker(F)$, so it
follows from Lemma~\ref{lm:kerF}(1) that
$$
\bbf{u}\equiv 0 \mod \pi^{l_2},
$$
hence also
$$
[\omega]\bbf{u}\equiv V\bbf{u}\equiv 0 \mod \pi^{l_2}.
$$
Thus by Lemma~\ref{lm:kerF}(2), their sum in $W(R/\pi^{pl_2}R)$ is
computed componentwise:
$$
T_{p[\pi^{l_1}]/\pi^{pl_1}}\bbf{u} = \left(\omega
u_0,\omega^pu_1+u_0,\omega^{p^2}u_2+u_1,
\omega^{p^3}u_3+u_2,\dots\right).
$$
Since $\bbf{u}$ has finitely many nonzero coefficients, we may
call $u_k$ the last of them. It follows from the above that
$$
(pa_{12}-\pi^{l_1},(a_{12})^p,0,\dots)= \left(\omega
u_0,\omega^pu_1+u_0,\omega^{p^2}u_2+u_1,\dots
\omega^{p^3}u_k+u_{k-1},u_k\dots\right).
$$
This is possible only if $k=0$, hence
$\bbf{u}=(u_0,0,\dots)=[u_0]$ in $W(R/\pi^{pl_2}R)$. Now if we
identify
$$
T_{p[\pi^{l_1}]/\pi^{pl_1}}\bbf{u} =\left(\omega
u_0,u_0,0,0,\dots\right) =(pa_{12}-\pi^{l_1},(a_{12})^p,0,\dots),
$$
we obtain $u_0\equiv (a_{12})^p \mod \pi^{pl_2}$ and the
congruence
$$
pa_{12}-\pi^{l_1}- \frac{p}{\pi^{(p-1)l_1}}(a_{12})^p\equiv 0 \mod
\pi^{pl_2}.
$$
This finishes the proof of the first assertion of the proposition.
As a bonus, we see that $\bbf{u}=F([a_{12}])$ is a solution to
\eqref{eq:equation n=2}. Then Theorem \ref{result4} shows that
$\cE(A)/G\simeq \cE(F(A))$. The final expression for the
function ring of $G$ follows from the general theory, which of
course provides also the group law. The final statement follows
from the above lemma.
\end{proo}

\begin{remas}\label{rema:n=2 tous les modeles}
(1) In the set of conditions (i)-(ii)-(iii) of the proposition, the
inequality $l_1\ge l_2$ is a consequence of the rest. Indeed, if
we assume the three conditions satisfied except that $l_1<l_2$, it
is clear that (iii) has no solution.

\smallskip

\noindent (2) We recover, with essentially the same proof, all the
group schemes exhibited in Tossici's paper \cite{To2}. In {\em
loc. cit.} it is also proven with some more work that if
$(l_1,l_2,a_{12})$ and $(l_1,l_2, a'_{12})$ give rise to
isomorphic group schemes then $a_{12}\equiv a'_{12}\mod
\pi^{l_2}$, and that all models of $\mu_{p^2}$ are obtained in
this way. Moreover Proposition 3.34 in {\em loc. cit.} can be
complemented by saying that the existence of model maps
corresponds to the divisibility between their matrices (with
Teichm\"uller entries). All these things works also in
characteristic $2$.
\end{remas}

The following remark is the equivalent of \ref{rm:recap}.

\begin{rema} \label{consequences utiles}
The above results for $n=2$ have consequences for general $n$. Let
$G$ be a finite flat Kummer group scheme in a filtered group
scheme $\cE(A)$ of type $(l_1,\dots,l_n)\in \NN^n$. Then
\begin{itemize}
\item[(1)] $\frac{e}{p-1}\ge l_1$,
\item[(2)] $l_1\ge l_2 \ge \dots \ge l_n$,
\item[(3)] if
$\bbf{a}_{i,i+1}=(a_{i,i+1}^0,a_{i,i+1}^1,\dots, a_{i,i+1}^k,
\dots )$ then $v(a_{i,i+1}^k)\ge l_{i+1}/p$ for all $i,k$.
\end{itemize}
In fact point (1) is already known (first sentences of \ref{Compn2}).
For $n=2$, point (2) follows from Lemma~\ref{l_1 > l_2} and point (3)
follows from Lemmas~\ref{lm:belonging_to_sM_2} and \ref{lm:kerF}(1).
The statement for arbitrary~$n$ follows simply by considering the $n-1$
subquotients of $G$ of order $p^2$, whose matrices are the diagonal
blocks of size $(2,2)$ of the matrix $A$ (see
Proposition~\ref{pp:pos_mat_and_model_maps}).
\end{rema}

\begin{noth}
{\bf Comparison Sekiguchi-Suwa Theory / Breuil-Kisin Theory for
$n=2$.} \label{Comp_SST_BKT_n=2} The existence of a link between
these two theories was already known in \cite{To2}, Appendix~A but
in a less precise way. Explaining it in details using our
formalism will give an idea of what the problems are for $n>2$. We
shall construct an explicit bijection between the set of matrices
parametrizing models of $\mu_{p^2}$ viewed as Kummer group
schemes, and the set of matrices parametrizing Breuil-Kisin
lattices.

We recall the setting: $R$ is a complete discrete valuation ring
with perfect residue field~$k$, totally ramified over $W(k)$. We
fix a uniformizer $\pi\in R$ and we call $E(u)$ its minimal
polynomial over $K$, so that $u\mapsto \pi$ induces an isomorphism
$W(k)[u]/(E(u))\simeq R$. Note that since $E(u)$ is Eisenstein, we
have $E(u)\equiv u^e+p\,[E_1(u)]\mod p^2$ with $E_1(u)\in k[u]$,
$\deg(E_1(u))<e$ and $E_1(0)\ne 0$.

The central point in the dictionary between the two theories is
the map $(-)^*:k[[u]]\to R$ sending a power series
$c=\sum_{i=0}^{\infty}{c_i}u^i$ to
$c^*=\sum_{i=0}^{\infty}{[c_i]}\pi^i$. It is an isometry for the
$u$-adic distance on the domain and the $\pi$-adic distance on the
target, which means simply that $f\equiv g \mod u^l$ if and only
if $f^*\equiv g^*\mod \pi^l$, for all $l \ge 0$. Moreover, we have
the property $(u^nc)^*=\pi^n c^*$. For each $l\ge 1$, the map
$c\mapsto c^*$ induces a map $k[u]/u^lk[u]\simeq R/\pi^{l}R$ which
for $l\le e$ is an isomorphism of rings but is neither additive nor
multiplicative in general. Now, using $(-)^*$ we map any matrix
$$
A=\left(
\renewcommand{\arraystretch}{1.5}
\begin{array}{cc}
u^{l_1} & a_{12} \\
0 & u^{l_2} \\
\end{array}\right)
\in \sG_2((u))
$$
to the matrix
$$
A^*=\left(
\renewcommand{\arraystretch}{1.5}
\begin{array}{cc}
[\pi^{l_1}] & [a^*_{12}] \\
0 & [\pi^{l_2}] \\
\end{array}\right)
\in M_2(W(R)).
$$
We claim that $A$ is a $\mu$-matrix if and only if $A^*$ gives
rise to a model of $\mu_{p^2}$.
In order to prove this, we just have to check that the congruences
in the two columns correspond to each other:

\medskip

\begin{center}
\begin{tabular}{|c|c|}
\hline Breuil-Kisin ({\em cf} \ref{coro1})
& Sekiguchi-Suwa ({\em cf} \ref{coroTo}) \\
\hline
& \\
$\frac{e}{p-1}\ge l_1\ge l_2$ & $\frac{e}{p-1}\ge l_1 \ge l_2$ \\ & \\
$(a_{12})^p\equiv 0 \mod u^{l_2}$ & $(a^*_{12})^p\equiv 0\mod \pi^{l_2}$ \\ & \\
$u^e a_{12}+E_1(u)u^{l_1}-u^{e-(p-1)l_1}(a_{12})^p\equiv 0\mod
u^{pl_2}$ &
$pa^*_{12}-\pi^{l_1}-\frac{p}{\pi^{(p-1)l_1}}(a^*_{12})^p
\equiv 0 \mod\pi^{pl_2}$ \\ & \\
\hline
\end{tabular}
\end{center}

\medskip

Since $l_2\le e$, the equivalence between the congruences in the
second line comes from the isomorphism $k[u]/u^lk[u]\simeq
R/\pi^{l}R$. It remains only to prove the equivalence of the
congruences in the third line: this is not immediate since
$R/\pi^{pl_2}R$ is not isomorphic to $k[u]/u^{pl_2}$ if $pl_2>e$.
So we look at the image under $c\mapsto c^*$ of the Breuil-Kisin
congruence
\begin{equation}\label{eq:BK_cong}
u^e a_{12}+E_1(u)u^{l_1}\equiv u^{e-(p-1)l_1}(a_{12})^p\mod
u^{pl_2}.
\end{equation}
We compute the image of both sides.
Since, for $\alpha,\beta \in k$, the difference
$[\alpha+\beta]-[\alpha]-[\beta]\in W(k)$ is a multiple of $p$, one
sees that the difference between $(u^ea_{12}+E_1(u)u^{l_1})^*$ and
$\pi^e a^*_{12}+[E_1](\pi)\pi^{l_1}$ is a multiple of $p\pi^e$.
Hence using the fact that $2e>pl_2$ we see that
$$
(u^e a_{12}+E_1(u)u^{l_1})^*\equiv
\pi^e a^*_{12}+[E_1](\pi)\pi^{l_1} \mod \pi^{pl_2},
$$
where $[E_1](\pi)$ is the evaluation of the
polynomial $[E_1]$ at $u=\pi$.
Now using $\val_u(a_{12})\ge l_2/p$ and $\frac{e}{p-1}\ge l_2$, one
sees using the binomial theorem that $((a_{12})^p)^*=(a^*_{12})^p
\mod \pi^{pl_2}$. Putting things together, it follows that the image
of the congruence \eqref{eq:BK_cong} is:
$$
\pi^e
a^*_{12}+[E_1](\pi)\pi^{l_1}\equiv\pi^{e-(p-1)l_1}(a^*_{12})^p
\mod \pi^{pl_2}.
$$
Since $E$ vanishes at $u=\pi$, we have $\pi^e+p [E_1](\pi)\equiv
0\mod p^2$ (beware that $[E_1(u)]$ evaluated at $u=\pi$ is
$[E_1](\pi)$). Given that $p^2\equiv 0\mod \pi^{pl_2}$, we may
replace $\pi^e$ by $-p [E_1](\pi)$ in the previous congruence and obtain:
$$
p\,[E_1](\pi)
a^*_{12}-[E_1](\pi)\pi^{l_1}+\frac{p}{\pi^{(p-1)l_1}}
[E_1](\pi)(a^*_{12})^p \equiv 0\mod \pi^{pl_2}.
$$
Since $[E_1](\pi)$ is invertible mod $\pi^{pl_2}$, this is indeed
equivalent to the equation on the Sekiguchi-Suwa side in the third
line. Our claim is thus proved.

As we noticed in~\ref{rema:n=2 tous les modeles}(2), the results
of \cite{To2} imply that in fact the matrices with Teichm\"uller
coefficients that we are considering above on the Sekiguchi-Suwa
side are in one-one correspondence with the models of
$\mu_{p^2,K}$. Since the map $A\mapsto A^*$ also preserves
divisibility between matrices, it follows that we have set up a
covariant equivalence between the category of $\mu$-matrices and
the category of models of $\mu_{p^2,K}$. We have not proven that
this equivalence is the covariant equivalence constructed by
Kisin, but it seems natural to conjecture that they are indeed the
same.
\end{noth}

\subsection{Computations for $n=3$}

Fix $l_1,l_2,l_3>0$. We consider a matrix:
$$
A=\left(
\renewcommand{\arraystretch}{1.5}
\begin{array}{ccc}
[\pi^{l_1}] & \bbf{a}_{12} & \bbf{a}_{13} \\
0 & [\pi^{l_2}] & \bbf{a}_{23} \\
0 & 0 & [\pi^{l_3}] \\
\end{array}\right)\in \sH_3(W^f(R)).
$$

\begin{lemm} \label{lm:belonging_to_sM_3}
The condition $A\in\sM_3$, i.e. $F(A)/A\ge 0$, is equivalent to the
congruences:
\begin{trivlist}
\itemn{1} $F^{(l_1)}(\bbf{a}_{12})\equiv 0\mod \pi^{l_2}$,
$F^{(l_2)}(\bbf{a}_{23})\equiv 0\mod \pi^{l_3}$, and
\itemn{2}
$F^{(l_1)}(\bbf{a}_{13})\equiv
T_{\frac{F^{(l_1)}(\bbf{a}_{12})}{\pi^{l_2}}}(\bbf{a}_{23})
\mod\pi^{l_3}$.
\end{trivlist}
Moreover if $A'=({\bbf{a}'}_{\!ij})\in \sM_3$ then  $\cE(A')\simeq
\cE(A)$ as filtered group schemes if and only if
${\bbf{a}'}_{\!12}\equiv\bbf{a}_{12}\in W(R/\pi^{l_2}R)$,
${\bbf{a}'}_{\!23}\equiv\bbf{a}_{23}\in W(R/\pi^{l_3}R)$ and
\begin{equation}\label{eq:congruences for equivalence}
{\bbf{a}'}_{\!13}\equiv
\bbf{a}_{13}+T_{\frac{{\bbf{a}'}_{\!12}-\bbf{a}_{12}}{\pi^{l_2}}}(\bbf{a}_{23})\mod
\pi^{l_3}.
\end{equation}
If $A$ has Teichm\"uller entries $[a_{ij}]$ and $l_1\ge
l_2 \ge l_3$ then $F(A)/A\ge 0$ is equivalent to the congruences:
\begin{trivlist}
\itemn{1} $a_{12}^p\equiv 0\mod \pi^{l_2}$, $a_{23}^p\equiv 0\mod
\pi^{l_3}$, and \itemn{2} $\pi^{l_2}a_{13}^p\equiv a_{23}a_{12}^p
\mod\pi^{l_2+l_3}$.
\end{trivlist}
Moreover if $A'=([{{a}'}_{\!ij}])\in \sH_3(W^f(R))$ with
${{a}'}_{\!ii}=\pi^{l_i}$, then  $A'\in \sM_3$ and $\cE(A)\simeq \cE(A')$
as filtered group schemes if and only if
${{a}'}_{\!12}\equiv{a}_{12} \mod \pi^{l_2}$,
${{a}'}_{\!23}\equiv{a}_{23} \mod \pi^{l_3}$ and
\begin{equation}
[{{a}'}_{\!13}]-
[{a}_{13}]\equiv\big[\frac{({{a}'}_{\!12}-{a}_{12}){a}_{23}}{\pi^{l_2}}\big]\mod
\pi^{l_3}.
\end{equation}
\end{lemm}

\begin{rema}\label{rema: v(a13)>l3/p2}
Here is a remark for later use. Let us suppose that $pe\ge
(p-1)l_3$  and $l_1,l_2\ge l_3$. Let
$\bbf{a}_{i3}=(a_{i3}^0,a_{i3}^1,\dots, a_{i3}^k,\dots)$ for
$i=1,2$. Then, if $\bbf{a}_{13}$  and $\bbf{a}_{23}$ satisfy the
congruences of the first part of the above lemma then
$v(a_{13})\ge l_{3}/p^2$. To prove this we first observe
that since $F(\bbf{a}_{23})\equiv 0\mod \pi^{l_3}$ it follows from
Lemma \ref{lm:kerF}(1) that $v(a_{23}^k)\ge l_{3}/p$, for
any $k$. Hence from
$$
F^{(l_1)}(\bbf{a}_{13})\equiv
T_{\frac{F^{(l_1)}(\bbf{a}_{12})}{\pi^{l_2}}}(\bbf{a}_{23})
\mod\pi^{l_3}
$$
it follows  that all the components of $F(\bbf{a}_{13})$ have valuation
at least $l_3/p$. Again by Lemma \ref{lm:kerF}(1) we obtain
$v(a_{13}^k)\ge l_3/p^2$ for any $k$.
\end{rema}

\begin{proo}
We begin with the general case. By definition we have $F(A)/A\ge
0$ if and only if there exists a positive matrix
$$
C=\left(
\renewcommand{\arraystretch}{1.5}
\begin{array}{ccc}
[\pi^{(p-1)l_1}] & \bbf{c}_{12} & \bbf{c}_{13} \\
0 & [\pi^{(p-1)l_2}] & \bbf{c}_{23} \\
0 & 0 & [\pi^{(p-1)l_3}] \\
\end{array}\right)
$$
such that $F(A)=C\star_T A$. By the case $n=2$, this gives the
congruences in (i). The equality of entries in position $(1,3)$
gives:
$$
F^{(l_1)}(\bbf{a}_{13})=T_{\bbf{c}_{12}}(\bbf{a}_{23})+\pi^{l_3}.\bbf{c}_{13}.
$$
The coefficient $\bbf{c}_{12}$ is determined by the equality of entries
in position $(1,2)$, namely it is equal to
$\frac{F^{(l_1)}([a_{12}])}{\pi^{l_2}}.\bbf{c}_{12}$. This gives
the congruence (ii).

Let $A'$ be another matrix in $\sM_3$. Then by
\ref{pos_matrices_model_maps} we have $\cE(A)\simeq \cE(A')$
as filtered group schemes  if and only if there exist Witt vectors
$\bbf{r}',\bbf{s}',\bbf{t}'$ such that
$$
A':=\left(
\renewcommand{\arraystretch}{1.5}
\begin{array}{ccc}
[\pi^{l_1}] & {\bbf{a}'}_{12} & {\bbf{a}'}_{13} \\
0 & [\pi^{l_2}] & {\bbf{a}'}_{23} \\
0 & 0 & [\pi^{l_3}] \\
\end{array}\right)
= \left(
\renewcommand{\arraystretch}{1.5}
\begin{array}{ccc}
1 & \bbf{r}' & \bbf{t}' \\
0 & 1 & \bbf{s}' \\
0 & 0 & 1 \\
\end{array}\right)
\star_T \left(
\renewcommand{\arraystretch}{1.5}
\begin{array}{ccc}
[\pi^{l_1}] & \bbf{a}_{12} & \bbf{a}_{13} \\
0 & [\pi^{l_2}] & \bbf{a}_{23} \\
0 & 0 & [\pi^{l_3}] \\
\end{array}\right).
$$
It is easy to see that
$\bbf{r}'=\frac{{\bbf{a}'}_{\!12}-\bbf{a}_{12}}{\pi^{l_2}}$ and the
rest follows.

We now show how things simplify if one supposes that $A$ has
Teichm\"uller entries and $l_1\ge l_2\ge l_3$.

Formulas in (1) are immediate, under these hypothesis, from the
general case, and we already found them in the case $n=2$. Now let
us suppose that
$$
F^{(l_1)}([a_{13}])\equiv
T_{\frac{F^{(l_1)}([{a}_{12}])}{\pi^{l_2}}}([{a}_{23}])
\mod\pi^{l_3}.
$$
Since $(p-1)l_1\ge 2l_1\ge l_2+l_3$, this is equivalent to say
$$
a_{13}^p\equiv \frac{a_{12}}{\pi^{l_2}}a_{23} \mod \pi^{l_3},
$$
which is equivalent to (2).


We now  study when two matrices with Teichm\"uller entries are
equivalent.  The assertions about $a_{12}$ and $a_{23}$ clearly
come from the general case. Now let us take an upper triangular
matrix $A'=([a'_{ij}])\in \sM_3$. The condition
\eqref{eq:congruences for equivalence} reads, in  this case,
 $$
[a'_{13}]\equiv
[{a}_{13}]+T_{\frac{[a'_{12}]-[{a}_{12}]}{\pi^{l_2}}}([{a}_{23}])\mod
\pi^{l_3}.
$$
Since $a'_{12}=a_{12}+\pi^{l_2} r$ for some $r\in R$, $l_2\ge l_3$
and  $v(a_{12}),v(a_{23})\ge \frac{l_3}{p}$ we have
$$T_{\frac{[a'_{12}]-[{a}_{12}]}{\pi^{l_2}}}([a_{23}])\equiv
\frac{[(a'_{12}-a_{12})a_{23}]}{\pi^{l_2}}\mod \pi^{l_3}.
$$
 So we get
$$
[{a'}_{13}]\equiv
[{a}_{13}]+\frac{[({a'}_{12}-{a}_{12}){a}_{23}]}{\pi^{l_2}}\mod
\pi^{l_3},
$$
as desired.
\end{proo}

We now state our final result for $n=3$, and we provide some comments
after the statement.

\begin{theo} \label{coroToco}
Let $0\leq l_3\leq l_2\leq l_1\leq e/(p-1)$ be integers and
$a_{12}, a_{23}, a_{13}\in R$ elements satisfying the congruences:
$$
a_{12}^p\equiv 0\mod \pi^{l_2} \quad,\quad
a_{23}^p\equiv 0\mod \pi^{l_3} \quad,\quad
\pi^{l_2}a_{13}^p\equiv a_{23}a_{12}^p \mod\pi^{l_2+l_3}.
$$
Let $A=([a_{ij}])$ be the matrix with Teich\"muller entries of $\sM_3$ defined by these parameters
(Lemma~\ref{lm:belonging_to_sM_3}). Assume that $l_1\ge pl_3$.
Then the pre-Kummer subgroup $G\subset \cE(A)$ is finite flat if and
only if the following congruences are satisfied:
$$
\begin{array}{c}
\displaystyle
pa_{12}-\pi^{l_1}-\frac{p}{\pi^{(p-1)l_1}}a_{12}^p
\equiv 0\mod\pi^{pl_2}, \\ \\
\displaystyle
pa_{23}-\pi^{l_2}-\frac{p}{\pi^{(p-1)l_2}}a_{23}^p
\equiv 0\mod \pi^{pl_3}, \\ \\
\displaystyle
\frac{p}{\pi^{(p-1)l_1}}a_{13}^p \equiv pa_{13}-a_{12}-a_{23}^p
\,\frac{pa_{12}-\pi^{l_1}-\frac{p}{\pi^{(p-1)l_1}}a_{12}^p}{\pi^{pl_2}}
\mod \pi^{pl_3}.
\end{array}
$$
When this is the case, we have
$$
G=\Spec \left(R[T_1,T_2,T_3]\Bigg/ \left(
\begin{array}{l}
\frac{(1+\pi^{l_1}T_1)^p}{\pi^{l_1p}},
\frac{(D_1(T_1)+\pi^{l_2}T_2)^p(1+\pi^{l_1}T_1)^{-1}-1}{\pi^{pl_2}}, \\ \\
\frac{(D_2(T_1,T_2)+\pi^{l_3}T_3)^p
(D_1(T_1)+\pi^{l_2}T_2)^{-1}-1} {\pi^{pl_3}} \\
\end{array}
\right)\right)
$$
where $D_1(T)=E_p(a_{12}T)=\sum_{k=0}^{p-1}a_{12}^k\frac{T^k}{k!}$ and
$D_2(T_1,T_2)$ is a lifting of
$$
E_p(a_{13}T_1)E_p\Big(a_{23}\frac{T_2}{D_1(T_1)}\Big) \mod \pi^{l_3},
$$
which under the above congruences is a polynomial.
Finally  if $A'=([{{a}'}_{ij}])\in \sM_3$  then the finite and flat group scheme of $\cE(A')$ is isomorphic to $G$ if and only if
${{a}'}_{12}\equiv{a}_{12}\in R/\pi^{l_2}R$,
${{a}'}_{23}\equiv{a}_{23}\in R/\pi^{l_3}R$ and
\begin{equation}
[{{a}'}_{13}]-
[{a}_{13}]\equiv\big[\frac{({{a}'}_{12}-{a}_{12}){a}_{23}}{\pi^{l_2}}\big]\mod
\pi^{l_3}.
\end{equation}
\end{theo}



\begin{rema}
(1) This result says that we are able to describe completely the
congruences  satisfied by matrices with Teichm\"uller entries giving
rise to Kummer group schemes of order $p^3$, under the (light) assumption
that $l_1\ge pl_3$. Removing this assumption would require more work.
See the final remarks in~\ref{ss:conclusion} for more comments on the
case $l_1< pl_3$. However  we do not know if  Kummer group schemes of order
$p^3$ arising from matrices with Teichm\"uller entries provide all the
Kummer group schemes of order $p^3$, under the hypothesis $l_1\ge pl_3$.

\smallskip

\noindent (2) A consequence of the above statement is that in the
situation of~\ref{coroToco}, we may take $B=F(A)$ in Theorem~\ref{result4}.
See also Remark~\ref{rm:A_and_FA}.

\smallskip

\noindent (3) In the third congruence of the second set of congruences,
one may in fact remove the term $p a_{13}$ since $e\ge l_1\ge pl_3$.
But leaving it emphasizes the similarity with the congruences we
obtained with the Breuil-Kisin approach, as we will see in~\ref{Comp_SST_BKT_n=3}.

\smallskip

%
\end{rema}

\begin{proo}
The dependency on the parameters (in the end of the statement)
follows from Lemma~\ref{lm:belonging_to_sM_3}. The rest is proven
by the general theory, except for the precise shape of the
congruences. Proposition~\ref{coroTo} gives the congruences for
the subgroup and quotient of degree $p^2$ in $G$ to be finite
flat, and fills up the upper left and lower right matrices of size
$2$. More precisely, we have the congruences:
\begin{align*}
pa_{12}-\pi^{l_1}-\frac{p}{\pi^{(p-1)l_1}}(a_{12})^p
\equiv 0 \mod\pi^{pl_2}, \\
pa_{23}-\pi^{l_2}-\frac{p}{\pi^{(p-1)l_2}}(a_{23})^p
\equiv 0 \mod\pi^{pl_3}.
\end{align*}
With the notations of Theorem~\ref{result4}, we have:
$$
B=B_3=\left(
\renewcommand{\arraystretch}{1.5}
\begin{array}{ccc}
[\pi^{pl_1}] & [a_{12}^p] & \bbf{u}_{13} \\
0 & [\pi^{pl_2}] & [a_{23}^p] \\
0 & 0 & [\pi^{pl_3}] \\
\end{array}\right)
$$
and
$$
V_3=\left(
\renewcommand{\arraystretch}{1.5}
\begin{array}{ccc}
p[\pi^{l_1}]/\pi^{pl_1} & \bbf{v}_{12} & \bbf{v}_{13} \\
0 & p[\pi^{l_2}]/\pi^{pl_2} & \bbf{v}_{23} \\
0 & 0 & p[\pi^{l_3}]/\pi^{pl_3} \\
\end{array}\right)
$$
where $\bbf{v}_{12}=\bbf{v}_1^2$ and $\bbf{v}_{23}=\bbf{v}_2^3$ are
the following vectors of $W(R)$:
\begin{align*}
\bbf{v}_{12}=\frac{1}{\pi^{pl_2}}
\left(p[a_{12}]-[\pi^{l_1}]-T_{p[\pi^{l_1}]/\pi^{pl_1}}[a_{12}^p]\right), \\
\bbf{v}_{23}=\frac{1}{\pi^{pl_3}}
\left(p[a_{23}]-[\pi^{l_2}]-T_{p[\pi^{l_2}]/\pi^{pl_2}}[a_{23}^p]\right).
\end{align*}
Thus $G$ is finite flat if and only if the previous congruences
are satisfied as well as the following last one:
$$
p[a_{13}]-[a_{12}]-T_{p[\pi^{l_1}]/\pi^{pl_1}}\bbf{u}_{13}
-T_{\bbf{v}_{12}}[a_{23}^p]\equiv 0 \mod \pi^{pl_3}.
$$
It only remains to prove that this is equivalent to:
$$
\frac{p}{\pi^{(p-1)l_1}}a_{13}^p \equiv pa_{13}-a_{12}-a_{23}^p
\,\frac{pa_{12}-\pi^{l_1}-\frac{p}{\pi^{(p-1)l_1}}a_{12}^p}{\pi^{pl_2}}
\mod \pi^{pl_3}.
$$
This is done in Subsection~\ref{completion_of_the_proof}.
\end{proo}

\begin{noth}
{\bf Comparison Sekiguchi-Suwa Theory / Breuil-Kisin Theory for $n=3$.}
\label{Comp_SST_BKT_n=3}
We proceed as in~\ref{Comp_SST_BKT_n=2} to compare the two theories.
Since we conducted the computations only under the additional assumption
$l_1\ge pl_3$ on the Sekiguchi-Suwa side, we will stay in this restricted
setting. We consider again the map $(-)^*:k[[u]]\to R$,
$\sum_{i=0}^{\infty}{c_i}u^i \mapsto \sum_{i=0}^{\infty}{[c_i]}\pi^i$
and the induced map on matrices:
$$
A=\left(
\renewcommand{\arraystretch}{1.5}
\begin{array}{ccc}
u^{l_1} & a_{12} & a_{13} \\
   0   & u^{l_2} & a_{23} \\
   0   &   0     & u^{l_3} \\
\end{array}
\right)
\in\sG_3((u))
\quad \mto \quad
A^*=\left(
\renewcommand{\arraystretch}{1.5}
\begin{array}{ccc}
[\pi^{l_1}] & [a^*_{12}] & [a^*_{13}] \\
   0   & [\pi^{l_2}] & [a^*_{23}] \\
   0   &   0     & [\pi^{l_3}] \\
\end{array}
\right) \in M_3(W(R)).
$$
We want to check that $A$ is a $\mu$-matrix if and only if $A^*$
gives rise to a model of $\mu_{p^3}$. For this we compare the
congruences from Breuil-Kisin Theory ({\em cf} \ref{coro1}) on the
left, and the congruences from Sekiguchi-Suwa Theory ({\em cf}
\ref{coroTo}) on the right.

\begin{center}
\begin{tabular}{cc}
{\bf A}
\fbox{\begin{minipage}{10cm}
\begin{center}
$a_{12}^p\equiv 0\mod u^{l_2}$, $a_{23}^p\equiv 0\mod u^{l_3}$
\end{center}
\end{minipage}}
& \begin{minipage}{5cm} \ \end{minipage} \\
\end{tabular}
\end{center}

\vspace{-7mm}

\begin{center}
\begin{tabular}{cc}
\begin{minipage}{5cm} \ \end{minipage} &
\fbox{\begin{minipage}{10cm}
\begin{center}
${a_{12}^*}^p\equiv 0\mod \pi^{l_2}$, ${a_{23}^*}^p\equiv 0\mod
\pi^{l_3}$
\end{center}
\end{minipage}} \\
\end{tabular}
\end{center}

\begin{center}
\begin{tabular}{cc}
{\bf B}
\fbox{\begin{minipage}{10cm}
\begin{center}
$u^ea_{12}+u^{l_1}E_1-u^{e-(p-1)l_1}a_{12}^p \equiv 0 \mod u^{pl_2}$ \\
$u^ea_{23}+u^{l_2}E_1-u^{e-(p-1)l_2}a_{23}^p \equiv 0 \mod u^{pl_3}$
\end{center}
\end{minipage}}
& \begin{minipage}{5cm} \ \end{minipage} \\
\end{tabular}
\end{center}

\vspace{-7mm}

\begin{center}
\begin{tabular}{cc}
\begin{minipage}{5cm} \ \end{minipage} &
\fbox{\begin{minipage}{10cm}
\begin{center}
$p{a_{12}^*}-\pi^{l_1}-\frac{p}{\pi^{(p-1)l_1}}{a_{12}^*}^p\equiv 0\mod\pi^{pl_2}$ \\
$p{a_{23}^*}-\pi^{l_2}-\frac{p}{\pi^{(p-1)l_2}}{a_{23}^*}^p\equiv
0\mod \pi^{pl_3}$
\end{center}
\end{minipage}} \\
\end{tabular}
\end{center}

\begin{center}
\begin{tabular}{cc}
{\bf C}
\fbox{\begin{minipage}{10cm}
\begin{center}
$a_{12}-u^{l_1-l_2}a_{23}\equiv 0\mod u^{l_3}$
\end{center}
\end{minipage}}
& \begin{minipage}{5cm} \ \end{minipage} \\
\end{tabular}
\end{center}

\vspace{-7mm}

\begin{center}
\begin{tabular}{cc}
\begin{minipage}{5cm} \ \end{minipage} &
\fbox{\begin{minipage}{10cm}
\begin{center}
???
\end{center}
\end{minipage}} \\
\end{tabular}
\end{center}

\begin{center}
\begin{tabular}{cc}
{\bf D}
\fbox{\begin{minipage}{10cm}
\begin{center}
$u^{l_2}a_{13}^p-a_{12}^pa_{23}\equiv 0\mod u^{l_2+l_3}$
\end{center}
\end{minipage}}
& \begin{minipage}{5cm} \ \end{minipage} \\
\end{tabular}
\end{center}

\vspace{-7mm}

\begin{center}
\begin{tabular}{cc}
\begin{minipage}{5cm} \ \end{minipage} &
\fbox{\begin{minipage}{10cm}
\begin{center}
$\pi^{l_2}{a_{13}^*}^p\equiv {a_{23}^*}{a_{12}^*}^p
\mod\pi^{l_2+l_3}$
\end{center}
\end{minipage}} \\
\end{tabular}
\end{center}

\begin{center}
\begin{tabular}{cc}
{\bf E}
\fbox{\begin{minipage}{10cm}
\begin{center}
$\begin{array}[t]{l}
u^ea_{13}+a_{12}E_1+\SS_1(u^ea_{12},u^{l_1}E_1)+u^{l_1}E_2 \\
-u^{e-(p-1)l_1}a_{13}^p
-\frac{u^ea_{12}+u^{l_1}E_1-u^{e-(p-1)l_1}a_{12}^p}{u^{pl_2}}a_{23}^p
\equiv 0 \mod u^{pl_3}
\end{array}$
\end{center}
\end{minipage}}
& \begin{minipage}{5cm} \ \end{minipage} \\
\end{tabular}
\end{center}

\vspace{-7mm}

\begin{center}
\begin{tabular}{cc}
\begin{minipage}{5cm} \ \end{minipage} &
\fbox{\begin{minipage}{10.5cm}
\begin{center}
$\begin{array}[t]{l} \frac{p}{\pi^{(p-1)l_1}}{a_{13}^*}^p \equiv
p{a_{13}^*}-{a_{12}^*}-{a_{23}^*}^p
\,\frac{p{a_{12}^*}-\pi^{l_1}-\frac{p}{\pi^{(p-1)l_1}}{a_{12}^*}^p}{\pi^{pl_2}}
\mod \pi^{pl_3}
\end{array}$
\end{center}
\end{minipage}} \\
\end{tabular}
\end{center}
\end{noth}

One sees immediately that in the Breuil-Kisin side there is one more
equation. We will prove below that in fact this congruence is a consequence
of the others, so we do not bother considering it for the moment.

In fact, only conditions {\bf D} and {\bf E} need to be compared,
since the previous ones match by the case $n=2$. The equivalence between the
congruences in {\bf D} is immediate since the operator $(-)^*$ induces an
isomorphism on the truncations of level $l_2+l_3$, given that
$e\ge l_{2}+l_3$. We pass to {\bf E}. Taking into account our assumption
that $l_1\ge pl_3$, on the Breuil-Kisin side we have
$$
u^{l_1}E_2\equiv \SS_1(u^ea_{12},u^{l_1}E_1)\equiv 0 \mod u^{pl_3}.
$$
With what remains, we can see that the two equations are
equivalent in the same way as in \ref{Comp_SST_BKT_n=2}. It is
still true that $p^2\equiv 0\mod \pi^{pl_3}$. So we have $\pi^e+p
[E_1](\pi)\equiv 0\mod p^2$ (beware that $[E_1(u)]$ evaluated at
$u=\pi$ is $[E_1](\pi)$). With no more difficulty than in the case
$n=2$ one shows, using $l_1\ge pl_3$, that
$(a_{13}^p)^*\equiv(a_{13}^*)^p\mod \pi^{l_3}$,
$(a_{23}^p)^*\equiv(a_{23}^*)^p\mod \pi^{l_3}$ and

\begin{align*}
\big(u^ea_{13}+a_{12}E_1-&\frac{u^ea_{12}+u^{l_1}E_1-u^{e-(p-1)l_1}a_{12}^p}{u^{pl_2}}a_{23}^p\big)^*\equiv
\\
&\pi^ea_{13}^*+a_{12}^*[E_1](\pi)-\frac{\pi^ea_{12}^*+\pi^{l_1}[E_1](\pi)-\pi^{e-(p-1)l_1}{a_{12}^p}}{\pi^{pl_2}}(a_{23}^*)^p\mod
\pi^{pl_3}.
\end{align*}
This gives the result. 

We now prove that the congruence $a_{12}\equiv u^{l_1-l_2}a_{23}
\mod u^{l_3}$, in the Breuil-Kisin side, is implied by the others.
We first observe that by the Breuil Kisin congruence in {\bf E},
and since $a_{23}^p\equiv 0\mod u^{l_{3}}$, we have
$$
a_{12}E_1\equiv -u^{e-(p-1)l_1}a_{13}^p \mod u^{l_3}.
$$
So, using $a_{13}^p\equiv a_{23}\frac{a_{12}^p}{u^{l_2}}\mod
u^{l_3}$, it follows that
\begin{equation}\label{eq:a12 modulo pil3}
a_{12}E_1\equiv -u^{e-(p-1)l_1}a_{23}\frac{a_{12}^p}{u^{l_2}} \mod
u^{l_3}.
\end{equation}
But if we divide $u^ea_{12}-u^{l_1}E_1\equiv
u^{e-(p-1)l_{1}}a_{12}^p\mod u^{l_2p}$ by $u^{l_2}$, and we
consider what we obtain modulo $u^{l_3}$, we get
$$
-u^{l_1-l_2}E_1\equiv u^{e-(p-1)l_{1}}\frac{a_{12}^p}{u^{l_2}}\mod
u^{l_3}.
$$
Putting this congruence inside \eqref{eq:a12 modulo pil3} one gets
the claim. Finally the sets of congruences on the left is equivalent
to the set of congruences on the right.

\subsection{Five lemmas} \label{completion_of_the_proof}

The lemmas in this subsection complete the proof of Theorem~\ref{coroToco}.
We use all the notations introduced in the statement and the proof of
the theorem.

\begin{lemm} \label{lm:lm1}
Modulo $\pi^{pl_3}$, we have:
$$
\bbf{v}_{12}=\frac{1}{\pi^{pl_2}}
\left(pa_{12}-\pi^{l_1}-\frac{p}{\pi^{(p-1)l_1}}a_{12}^p,
\sigma_1(pa_{12},-\pi^{l_1})-\sigma_1\big(pa_{12}-\pi^{l_1},
\frac{p}{\pi^{(p-1)l_1}}a_{12}^p\big),0,\dots\right).$$
\end{lemm}

\begin{proo}
We will compute
$\pi^{pl_2}\bbf{v}_{12}=
p[a_{12}]-[\pi^{l_1}]-T_{p[\pi^{l_1}]/\pi^{pl_1}}[(a_{12})^p]$
modulo $\pi^{p(l_2+l_3)}$. Since $v(a_{12})\geq l_2/p$,
we obtain $v(S_i([a_{12}],[a_{12}]))\geq l_2p^{i-1}$ for each $i\in\NN$.
For $i\geq 3$, we have $l_2p^{i-1}\geq l_2p^2\geq pl_2+pl_3$, hence
$$
p[a_{12}]=
\left(pa_{12},(1-p^{p-1})(a_{12})^p,*p^{p-1}(a_{12})^{p^2},\dots\right)
\in \hat{W}(R/\pi^{p(l_2+l_3)}R), \mbox{ with } v(*)\geq 0.
$$
Note that
$$
v(p^{p-1}(a_{12})^{p^2})\ge e(p-1)+pl_2\ge l_3(p-1)^2+pl_2 \ge pl_2+pl_3
$$
so finally
$$
p[a_{12}]=\left(pa_{12},(1-p^{p-1})(a_{12})^p,0,\dots\right)
\in \hat{W}(R/\pi^{p(l_2+l_3)}R).
$$
We now compute:
$$
p[a_{12}]-[\pi^{l_1}]=
\left(pa_{12}-\pi^{l_1},(1-p^{p-1})(a_{12})^p+\sigma_1(pa_{12},-\pi^{l_1}),
S_2(p[a_{12}],-[\pi^{l_1}]),\dots\right).
$$
Using the minorations $v (pa_{12})\geq e$,
$v(1-p^{p-1})(a_{12})^p\geq l_2$, $v(\pi^{l_1})=l_1$,
we obtain $S_i(p[a_{12}],-[\pi^{l_1}])=0$ in $R/\pi^{p(l_2+l_3)}R$.
$$
v(\sigma_2(pa_{12},-\pi^{l_1}))
\ge (p^2-1)l_1+e
\ge (p^2-1)l_1+(p-1)l_1
= pl_1+(p^2-2)l_1
\ge pl_2+pl_3
$$
so finally
$$
p[a_{12}]-[\pi^{l_1}]=
\left(pa_{12}-\pi^{l_1},(1-p^{p-1})(a_{12})^p+\sigma_1(pa_{12},-\pi^{l_1}),0,\dots\right).
$$
The last term contributing to $\bbf{v}_{12}$ is
$$
T_{p[\pi^{l_1}]/\pi^{pl_1}}[(a_{12})^p]=
\left(\frac{p\pi^{l_1}}{\pi^{pl_1}}(a_{12})^p,(1-p^{p-1})(a_{12})^p,0\dots\right)\in \hat{W}(R/\pi^{p(l_2+l_3)}R).
$$
We add up and we obtain the lemma.
\end{proo}

\begin{lemm}
In $\widehat{W}(R/\pi^{pl_3}R)$, we have the equality:
$$
T_{\bbf{v}_{12}}[(a_{23})^p]=
\left(\frac{(a_{23})^p}{\pi^{pl_2}}
\left(pa_{12}-\pi^{l_1}-\frac{p}{\pi^{(p-1)l_1}}(a_{12})^p\right),
\frac{(a_{23})^p}{\pi^{pl_2}}\sigma_1\left(pa_{12},-\pi^{l_1}\right),0,\dots
\right).
$$
\end{lemm}

\begin{proo}
Now we can compute, in $\hat{W}(R/\pi^{pl_3}R)$ this time:
\begin{align*}
T_{\bbf{v}_{12}}[(a_{23})^p]
& =
\left(\frac{(a_{23})^p}{\pi^{pl_2}}
\left(pa_{12}-\pi^{l_1}-\frac{p}{\pi^{(p-1)l_1}}(a_{12})^p\right), \right. \\
& \quad \left.
\frac{(a_{23})^p}{\pi^{pl_2}}\sigma_1\left(pa_{12},-\pi^{l_1}\right)
+\frac{(a_{23})^p}{\pi^{pl_2}}\sigma_1\left(pa_{12}-\pi^{l_1},
\frac{p}{\pi^{(p-1)l_1}}(a_{12})^p\right),0,\dots
\right).
\end{align*}
We simplify a little bit. Using the identity
$\sigma_1(x,y)=\sigma_1(x,-y-x)$, we get
$$
\sigma_1\left(pa_{12}-\pi^{l_1},\frac{p}{\pi^{(p-1)l_1}}(a_{12})^p\right)
=\sigma_1\left(pa_{12}-\pi^{l_1},
-\frac{p}{\pi^{(p-1)l_1}}(a_{12})^p+\pi^{l_1}-pa_{12}\right).
$$
We use the inequality
$v(\sigma_1(a,b))\ge \min((p-1)v(a)+v(b),(p-1)v(b)+v(a))$
from Lemma \ref{lemmar}(4).
In our case $a=pa_{12}-\pi^{l_1}$ has valuation $l_1$ and
$b=-\frac{p}{\pi^{(p-1)l_1}}(a_{12})^p+\pi^{l_1}-pa_{12}$ has
valuation at least $pl_2$, and we find
$$
v\left(\frac{(a_{23})^p}{\pi^{pl_2}}\sigma_1(a,b)\right)
\ge l_3-pl_2+\min\big((p-1)l_1+pl_2,l_1+(p-1)pl_2\big) \ge pl_3
$$
so this term vanishes. Finally we obtain the lemma.
\end{proo}

\begin{lemm}\label{leminter}
We have the following equalities in $\widehat{W}(R/\pi^{pl_3}R)$:
\begin{align*}
& p[a_{13}]-[a_{12}]-T_{\bbf{v}_{12}}[(a_{23})^p] \\
& \qquad = \left(pa_{13}-a_{12},(a_{13})^p+\sigma_1(pa_{13},-a_{12}),0,\dots\right) \\
& \qquad \quad - \left(\frac{(a_{23})^p}{\pi^{pl_2}}
\left(pa_{12}-\pi^{l_1}-\frac{p}{\pi^{(p-1)l_1}}(a_{12})^p\right),
\frac{(a_{23})^p}{\pi^{pl_2}}\sigma_1\left(pa_{12},-\pi^{l_1}\right),0,\dots
\right) \\
& \qquad = (c_0,c_1,c_2,\dots)
\end{align*}
with
$$
\begin{array}{cl}
c_0 &= pa_{13}-a_{12}-\frac{(a_{23})^p}{\pi^{pl_2}}
\left(pa_{12}-\pi^{l_1}-\frac{p}{\pi^{(p-1)l_1}}(a_{12})^p\right) \\
& \\
c_1 &=
(a_{13})^p 
\end{array}
$$
and
$c_i=0$, $i\geq 2$.
\end{lemm}

\begin{proo}
We have
\begin{align*}
p[a_{13}]-[a_{12}]
& = (pa_{13},(a_{13})^p,0,\dots)-(a_{12},0,\dots) \\
& = (pa_{13}-a_{12},(a_{13})^p+\sigma_1(pa_{13},-a_{12}),
\sigma_2(pa_{13},-a_{12}),\dots)\in\hat{W}(R/\pi^{pl_3}R).
\end{align*}
Recall from Lemma~\ref{lemmar}(4) that
$v(\sigma_2(a,b))\ge \min\left((p^2-1)v(a)+v(b),(p^2-1)v(b)+v(a)\right)$.
Using this we see immediately that $\sigma_2(pa_{13},-a_{12})\equiv 0 \mod \pi^{pl_3}$
so that
$$
p[a_{13}]-[a_{12}]=
(pa_{13}-a_{12},(a_{13})^p+\sigma_1(pa_{13},-a_{12}),0,\dots)\in\hat{W}(R/\pi^{pl_3}R).
$$
%
So $c_0$ is as in the statement and
\begin{align*}
c_1 & =
(a_{13})^p+\sigma_1(pa_{13},-a_{12})
-\frac{(a_{23})^p}{\pi^{pl_2}}\sigma_1\left(pa_{12},-\pi^{l_1}\right) \\
& \quad +\sigma_1\left(pa_{13}-a_{12},-\frac{(a_{23})^p}{\pi^{pl_2}}
\left(pa_{12}-\pi^{l_1}-\frac{p}{\pi^{(p-1)l_1}}(a_{12})^p\right)\right).
\end{align*}
Then we have
$v(\sigma_1(pa_{13},-a_{12}))\geq e+(p-1)l_2/p\geq (p-1)l_1 \geq pl_3$,
$$
\begin{array}{l}
v\Big( \frac{(a_{23})^p}{\pi^{pl_2}} \sigma_1\left(pa_{12},-\pi^{l_1}\right)\Big)\geq l_3-pl_2+l_1(p-1)+e+l_2/p\geq pl_3, \bigskip \\
v\Big( \sigma_1\Big(pa_{13}-a_{12},-\frac{(a_{23})^p}{\pi^{pl_2}}
\left(pa_{12}-\pi^{l_1}-\frac{p}{\pi^{(p-1)l_1}}(a_{12})^p\Big)\right)\Big)
\geq l_1/p+(p-1)l_3\geq pl_3.
\end{array}
$$
\end{proo}

\begin{lemm}
For  $\bbf{u}_{13}=(u_0,u_1,\dots)$, the condition
$$
p[a_{13}]-[a_{12}]-T_{p[\pi^{l_1}]/\pi^{pl_1}}\bbf{u}_{13}
-T_{\bbf{v}_{12}}[a_{23}^p]\equiv 0 \mod \pi^{pl_3}.
$$
implies in $R/\pi^{pl_3}R$:
$$
\begin{array}{l}
c_0=\frac{p}{\pi^{(p-1)l_1}}u_0 \\ \\
c_1-u_0=\left(\frac{p}{\pi^{(p-1)l_1}}\right)^pu_1 \\ \\
\sigma_1(c_1,-u_0)-u_1=0. \\
\end{array}
$$
\end{lemm}
\begin{proo}
Since $\frac{p[\pi^{l_1}]}{\pi^{pl_1}}\equiv\left(\frac{p}{\pi^{(p-1)l_1}},1,0,\dots\right)$
modulo $\pi^{pl_3}$,  then
$$
T_{\frac{p[\pi^{l_1}]}{\pi^{pl_1}}}\bbf{u}_{13}\equiv \left[\frac{p}{\pi^{(p-1)l_1}}\right]\bbf{u}_{13}+V\bbf{u}_{13}\mod \pi^{pl_3}.
$$
Then by the condition in the statement it follows
$$
\left[\frac{p}{\pi^{(p-1)l_1}}\right]\bbf {u}_{13}
= -V\bbf{u}_{13}+(c_0,c_1,0,\dots)
= -(0,u_0,u_1,\dots)+(c_0,c_1,0,\dots).
$$
Remark that $v(c_0)\ge \frac{l_2}{p}$ and $v(c_1)\ge l_3/p$. Also note that
by Remark \ref{rema: v(a13)>l3/p2} it follows that $v(u_i)\ge
l_3/p$, for all $i\geq 0$. We deduce:
$$
\left[\frac{p}{\pi^{(p-1)l_1}}\right]\bbf{u}_{13}
=\Big(c_0,c_1-u_0,\sigma_1(c_1,-u_0)-u_1,-u_2,\dots,-u_k,0,\dots\Big)
\in \hat{W}(R/\pi^{pl_3}R)
$$
where $u_k$ is by definition the last nonzero term in $\bbf{u}_{13}$
and $-u_k$ occurs here at the $(k+1)$-th place.
On the other hand,
$$
\left[\frac{p}{\pi^{(p-1)l_1}}\right]\bbf{u}_{13}=
\left(\frac{p}{\pi^{(p-1)l_1}}u_0,\left(\frac{p}{\pi^{(p-1)l_1}}\right)^pu_1,
\left(\frac{p}{\pi^{(p-1)l_1}}\right)^{p^2}u_2,\dots,
\left(\frac{p}{\pi^{(p-1)l_1}}\right)^{p^k}u_k,0,\dots\right)
$$
where here the term involving $u_k$ occurs at the $k$-th place. This is
not possible if $k\ge 2$. Hence $k\le 1$ and
$\bbf{u}_{13}=(u_0,u_1,0,\dots)\in \hat{W}(R/\pi^{pl_3}R)$.
And we obtain the expected formulas.
\end{proo}

\begin{lemm} \label{lm:lm6}
We have $u_1=0$. Therefore $u_0=a_{13}^p$ and the condition of
the previous lemma is equivalent to the congruence:
$$
\frac{p}{\pi^{(p-1)l_1}}a_{13}^p \equiv pa_{13}-a_{12}-a_{23}^p
\,\frac{pa_{12}-\pi^{l_1}-\frac{p}{\pi^{(p-1)l_1}}a_{12}^p}{\pi^{pl_2}}
\mod \pi^{pl_3}.
$$
\end{lemm}

\begin{proo}
We have $u_0=c_1-\Big(\frac{p}{\pi^{(p-1)l_1}}\Big)^p\sigma_1(c_1,-u_0)$.
Since $c_1$ divides $\sigma_1(c_1,-u_0)\in R/\pi^{pl_3}R$, $u_0=c_1u_0'$
with $v(  u_0')\geq 0$. Write $\beta=c_1\Big(\frac{p}{\pi^{(p-1)l_1}}\Big)$ then
$$
c_0=\beta-\Big(\frac{p}{\pi^{(p-1)l_1}}\Big)\beta^p\sigma_1(1,u_0')
$$
and $v(\beta)\geq v(c_1)\geq l_3/p$. We obtain
\begin{align*}
\beta &=c_0+\Big(\frac{p}{\pi^{(p-1)l_1}}\Big)\sigma_1(1,u_0')
(c_0+\Big(\frac{p}{\pi^{(p-1)l_1}}\Big)\beta^p\sigma_1(1,u_0'))^p \\
& =c_0+\Big(\frac{p}{\pi^{(p-1)l_1}}\Big)\sigma_1(1,u_0')c_0^p\in R/\pi^{pl_3}R.
\end{align*}
Recall that
$c_0=pa_{13}-a_{12}-\frac{(a_{23})^p}{\pi^{pl_2}}
\left(pa_{12}-\pi^{l_1}-\frac{p}{\pi^{(p-1)l_1}}(a_{12})^p\right)$.
Hence $c_0^p=-a_{12}^p\in R/\pi^{pl_3}R$.
Now we come back to the congruence
$$
pa_{12}-\pi^{l_1}-\frac{p}{\pi^{(p-1)l_1}}(a_{12})^p\equiv 0 \mod \pi^{pl_2}.
$$
If $e+ v(a_{12})\leq e-l_1(p-1)+p v(a_{12})$ then $v(a_{12})\geq l_1$
and $a_{12}^p=0\in R/\pi^{pl_3}R$.\\
If $e+v(a_{12})> e-l_1(p-1)+p v(a_{12})$, then
$e-l_1(p-1)+p v(a_{12})\geq \min(pl_2,l_1)$ and
$\Big(\frac{p}{\pi^{(p-1)l_1}}\Big)a_{12}^p=0\in
R/\pi^{pl_3}R$.
Hence $\beta=c_0\in R/\pi^{pl_3}R$.
Since $\beta=\Big(\frac{p}{\pi^{(p-1)l_1}}\Big)c_1$, we obtain the lemma.
\end{proo}

\subsection{Conclusion} \label{ss:conclusion}

The case $l_1< pl_3$ excluded in Proposition~\ref{coroToco} shows
the complexity of the ramification. To achieve this case,
we should compute the integrality conditions with matrices with
non-Teichm\"uller entries (see the introduction of Section~\ref{SSn3}).
More generally, for $n\geq 4$, in order to compute the Breuil-Kisin
modules of Kummer groups, we need to define adapted liftings of parameters
to $R$. In view of Theorem \ref{classification_BK}, these choices should
be related to $E(u) \mod p^n$ in some way.

We have seen that any $\mu$-matrix has to satisfy the condition
$\cU A/\cL A\ge 0$. This condition is not present in the context of the
classification of Kummer group schemes. And in the case $n=2$ and
$n=3$ (with $l_1\ge pl_3$ and matrices with Teichm\"uller entries)
 we have seen that in fact this condition is
consequence of the others. We do not know  if we could remove this
condition in the classification of $\mu$-matrices.

Let us emphasize that the explicit formulas for Kummer subgroups
are relevant not only in the perpective of classification of Hopf
orders of rank $p^n$ (\cite{ChUn}, \cite{GrCh}) but also for the
computation of dimension and irreducible components
(\cite{Car},\cite{Im}) of Kisin's variety parametrizing some group
schemes over $\cO_K$ (\cite{Ki2}). At last, Kummer group
schemes could be useful to give an explicit form for
Breuil-Kisin's equivalence of categories between $(\Mod/\fS)$ and
the category of finite flat group schemes of $p$-power order.

\bigskip

\begin{flushleft}
A. M\'ezard, {\sc Institut de Math\'ematiques de Jussieu,
Universit{\'e} Pierre et Marie Curie,
4~place Jussieu,
75252 Paris Cedex 05, France}

{\em Email address:} {\tt mezard@math.jussieu.fr}

\bigskip

M. Romagny, {\sc 
Institut de Recherche Math\'ematique de Rennes,
Universit\'e de Rennes 1,
Campus de Beaulieu,
35042 Rennes Cedex,
France}

{\em Email address:} {\tt matthieu.romagny@univ-rennes.fr}

\bigskip

D. Tossici, {\sc Scuola Normale Superiore di Pisa,
Piazza dei Cavalieri~7, 56126 Pisa, Italy.}

{\em Email address:} {\tt dajano.tossici@gmail.com}
\end{flushleft}


\begin{thebibliography}{200}

\bibitem[AR]{AR} {\sc D. Abramovich, M. Romagny},
{\it Moduli of Galois covers in mixed characteristics}, to
appear in Algebra and Number Theory.

\bibitem[AS]{AS} {\sc A. Abbes, T. Saito}, {\it Ramification of local fields
with imperfect residue fields I}, Amer. J. Math. 124 (2002), 879--920.

\bibitem[Br1]{Br1} {\sc C. Breuil},
{\em Sch\'emas en groupes et corps des normes},
unpublished manuscript, September~1998.

\bibitem[Br2]{Br2} {\sc C. Breuil}, {\it Integral $p$-adic Hodge Theory},
Algebraic geometry 2000, Azumino (Hotaka),  51--80,
Adv. Stud. Pure Math. 36, Math. Soc. Japan (2002).

\bibitem[By]{By} {\sc N.P. Byott},
{\em Cleft extensions of Hopf algebras},
Proc. London Math. Soc. 67 (1993), 227--307.

\bibitem[Car]{Car} {\sc X. Caruso}, {\it Estimation des dimensions de
certaines vari\'et\'es de Kisin}, preprint, arXiv:1005.2394.

\bibitem[Ch]{Ch} {\sc L. N. Childs},
{\em Taming Wild Extensions: Hopf Algebras and Local Galois Module Theory},
American Mathematical Society, Mathematical Surveys and Monographs 80, 2000.

\bibitem[ChUn]{ChUn} {\sc L. N. Childs, R.G. Underwood}, {\it Cyclic Hopf orders
defined by isogenies of formal groups}, Amer. J. of Math., 125
(2003), 1295--1334.

\bibitem[Ei]{Ei} {\sc D. Eisenbud}, {\it Commutative algebra with a view
  toward algebraic geometry}, Graduate Texts in Math. 150, Springer-Verlag (1995).

\bibitem[Fa]{Fa} {\sc L. Fargues}, {\it La filtration de Harder-Narasimhan des
sch\'emas en groupes finis et plats}, J. Reine Angew. Math., 645 (2010), 1--39.

\bibitem[Fo]{Fo} {\sc J.-M. Fontaine}, {\it Repr\'esentations $p$-adiques des corps
locaux I}, The Grothendieck Festschrift, Vol. II, 249--309, Progr. Math., 87,
Birkha\"user (1990).

\bibitem[Gr]{Gr} {\sc C. Greither}, {\it Extensions of finite group schemes,
and Hopf Galois theory over a complete discrete valuation ring},
Math. Z. 210 (1992), 37--67.

\bibitem[GrCh]{GrCh} {\sc C. Greither, L. Childs}, {\it $p$-Elementary group
schemes-constructions and Raynaud's Theory}, in Hopf algebra,
Polynomial formal Groups and Raynaud Orders, Mem. Amer. Soc.  136
(1998), 91--118.

\bibitem[Im]{Im} {\sc N. Imai}, {\it On the connected components of moduli spaces
of finite flat models}, to appear in Amer. J. Math.

\bibitem[KM]{KM} {\sc N. Katz, B. Mazur}, {\it Arithmetic moduli of
elliptic curves}, Annals of Mathematics Studies 108, Princeton University
Press (1985).

\bibitem[Ki1]{Ki1} {\sc M. Kisin},
{\it Crystalline representations and $F$-crystals},
Algebraic geometry and number theory, 459--496, Progr. Math. 253,
Birkh\"auser (2006).

\bibitem[Ki2]{Ki2} {\sc M. Kisin},
{\it Moduli of finite flat group schemes, and modularity},
Ann. of Math. (2) 170 (2009), no. 3, 1085--1180.

\bibitem[Lar]{Lar} {\sc R. Larson},
{\em Hopf algebra orders determined by group valuations},
J. Algebra 38 (1976), no. 2, 414--452. 

\bibitem[Lau]{Lau} {\sc E. Lau}, {\it A relation between Dieudonn\'e displays
and crystalline Dieudonn\'e theory}, preprint, arXiv:1006.2720.

\bibitem[Li]{Li} {\sc T. Liu}, {\it The correspondence between Barsotti-Tate
groups and Kisin modules when $p=2$}, preprint (2011).

\bibitem[Ma]{Ma} {\sc Y. I. Manin}, {\it Cubic forms. Algebra, geometry,
arithmetic}, Second edition, North-Holland (1986).

\bibitem[MRT]{MRT} {\sc A. M{\'e}zard, M. Romagny, D. Tossici},
{\it Sekiguchi-Suwa Theory revisited}, preprint (2011).

\bibitem[Ro]{Ro} {\sc M. Romagny}, {\it Effective models of group
schemes}, to appear in the Journal of Algebraic Geometry.

\bibitem[Se]{Se} {\sc J.-P. Serre}, {\it Corps Locaux}, Hermann (1980).

\bibitem[Sm]{Sm} {\sc J. D. H. Smith}, {\it An introduction to
quasigroups and their representations}, Studies in Advanced Mathematics,
Chapman~\&~Hall (2007).

\bibitem[SOS]{SOS} {\sc T. Sekiguchi, F. Oort, N. Suwa},
{\it On the deformation of Artin-Schreier to Kummer},
Ann. Sci. \'Ecole Norm. Sup. (4) 22 (1989), no. 3, 345--375.

\bibitem[SS1]{SS1} {\sc T. Sekiguchi, N. Suwa},
{\it On the unified Kummer-Artin-Schreier-Witt Theory}, no.~111 in the
preprint series of the Laboratoire de Math\'ematiques Pures de Bordeaux (1999).

\bibitem[SS2]{SS2} {\sc T. Sekiguchi, N. Suwa},
{\it A note on extensions of algebraic and formal groups. IV.
Kummer-Artin-Schreier-Witt theory of degree $p^2$},
Tohoku Math. J. (2) 53 (2001), no. 2, 203--240.

\bibitem[TO]{TO} {\sc J. Tate, F. Oort},
{\em Group schemes of prime order},
Ann. Sci. Ec. Norm. Sup., 3 (1970), 1--21.

\bibitem[To1]{To1} {\sc D. Tossici}, {\it Effective models and extension
of torsors over a discrete valuation ring of unequal characteristic},
Int. Math. Res. Not. IMRN 2008, Art. ID rnn111, 68 pp.

\bibitem[To2]{To2} {\sc D. Tossici}, {\it Models of $\mu_{p^2,K}$
over a discrete valuation ring. With an appendix by Xavier Caruso},
J. Algebra 323 (2010), no. 7, 1908--1957.

\bibitem[Un]{Un} {\sc R. Underwood},
{\em $R$-Hopf algebra orders in $KC_p^2$},
J. Alg., 169 (1994), 418--440.

\bibitem[WW]{WW} {\sc W. Waterhouse, B. Weisfeiler},
{\it One-dimensional affine group schemes},
J. Algebra 66 (1980), no. 2, 550--568.


\end{thebibliography}
\end{document}